\documentclass[a4paper]{amsart}

\usepackage[utf8]{inputenc}			
\usepackage[T1]{fontenc}			
\usepackage[british]{babel} 		

\usepackage{
amsmath,			
amsthm,				
amssymb,			
latexsym,			
xfrac,				
faktor,				
mathtools,			
bbm,				
mathrsfs,			
enumitem,			
microtype,			
hyphenat,			
hyperref,			
cite,				
cleveref,			
}

\newtheorem{theo}{Theorem}[section]
\newtheorem{prop}[theo]{Proposition}
\newtheorem{lem}[theo]{Lemma}
\newtheorem{coro}[theo]{Corollary}
\newtheorem*{claim}{Claim}

\theoremstyle{definition}
\newtheorem{ex}[theo]{Example}
\newtheorem{rmk}[theo]{Remark}

\crefname{theo}{Theorem}{Theorems}
\crefname{prop}{Proposition}{Propositions}
\crefname{lem}{Lemma}{Lemmas}
\crefname{coro}{Corollary}{Corollaries}
\crefname{defi}{Definition}{Definitions}
\crefname{ex}{Example}{Examples}
\crefname{rmk}{Remark}{Remarks}
\crefname{enumi}{}{}

\setenumerate{itemsep=0pt, topsep=3pt}

\renewcommand\qedsymbol{\setbox0=\hbox{\ \ \ \footnotesize{\normalfont Q.E.D.}}\kern\wd0 \strut \hfill \kern-\wd0 \box0}
\newcommand*\claimqed{\renewcommand\qedsymbol{\setbox0=\hbox{\ \ \ {\normalfont \ensuremath{\square}}}\kern\wd0 \strut \hfill \kern-\wd0 \box0}\qedhere\renewcommand\qedsymbol{\setbox0=\hbox{\ \ \ \footnotesize{\normalfont Q.E.D.}}\kern\wd0 \strut \hfill \kern-\wd0 \box0}}

\newcommand*\N{\mathbb{N}}
\newcommand*\Z{\mathbb{Z}}
\newcommand*\R{\mathbb{R}}
\newcommand*\Q{\mathbb{Q}}

\newcommand*\sth{\mathrel{:}}

\catcode`@=11 \def\rfrac#1#2{\def\rf@num{#1}\def\rf@den{#2}\mathpalette\rfr@c\relax} \def\rfr@c#1#2{\raise.5ex\hbox{$\mathsurround=0pt #1\rf@num$}\kern-0.1em/\kern-0.1em\lower.25ex\hbox{$\mathsurround=0pt #1\rf@den$}}\catcode`@=12

\newcommand*\Ord{\mathbf{\mathbbm{O}\mathrm{n}}}
\newcommand*\Ima{\mathrm{Im}}
\newcommand*\tp{\mathrm{tp}}
\newcommand*\inc{\iota}
\newcommand*\id{\mathrm{id}}
\newcommand*\lang{\mathtt{L}}
\newcommand*\For{\mathrm{For}}
\newcommand*\Aut{\mathrm{Aut}}
\newcommand*\cf{\mathrm{cf}}

\newcommand*\proj{\mathrm{proj}}
\newcommand*\quot{\mathrm{quot}}
\newcommand*\Stab{\mathrm{Stab}}
\newcommand*\St{\mathrm{St}}
\newcommand*\Stone{\mathrm{S}}

\begin{document}

\title{On piecewise hyperdefinable groups}
\author{Arturo Rodr\'{\i}guez Fanlo}
\address{Mathematical Institute, University of Oxford\\ 
OX2 6GG, Oxford, the United Kingdom}
\email{Arturo.RodriguezFanlo@math.ox.ac.uk}
\email[current]{Arturo.Rodriguez@mail.huji.ac.il}

\begin{abstract} The aim of this paper is to generalise and improve two of the main model{\hyp}theoretic results of ``Stable group theory and approximate subgroups'' by E. Hrushovski to the context of piecewise hyperdefinable sets. The first one is the existence of Lie models. The second one is the Stabilizer Theorem. In the process, a systematic study of the structure of piecewise hyperdefinable sets is developed. In particular, we show the most significant properties of their logic topologies. 
\end{abstract}


\maketitle

\section*{Introduction}
Various enlightening results were found in \cite{hrushovski2011stable} with significant consequences on model theory and additive combinatorics. Two of them are particularly relevant: the Stabilizer Theorem and the existence of Lie models. \smallskip

The Stabilizer Theorem \cite[Theorem 3.5]{hrushovski2011stable}, subsequently improved in \cite[Theorem 2.12]{montenegro2018stabilizers}, was originally itself a generalisation of the classical Stabilizer Theorem for stable and simple groups (see for example \cite[Section 4.5]{wagner2010simple}) changing the stability and simplicity hypotheses by some kind of measure{\hyp}theoretic ones. Here, we extend that theorem to piecewise hyperdefinable groups in Section 3. Once one has the right definition of dividing and forking for piecewise hyperdefinable sets, the original proof of \cite{hrushovski2011stable}, and its improved version of \cite{montenegro2018stabilizers}, can be naturally adapted. However, we also manage to simplify the proof in such a way that we get a slightly stronger result.\smallskip

In \cite{hrushovski2011stable}, Hrushovski studied piecewise definable groups generated by near{\hyp}subgroups, i.e. approximate subgroups satisfying a kind of measure{\hyp}theoretic condition. In particular, he worked with ultraproducts of finite approximate subgroups. In that context, Hrushovski proved that there exist some Lie groups, named \emph{Lie models}\footnote{In \cite{hrushovski2011stable} they are simply called \emph{associated Lie groups}. The term \emph{Lie model} was later introduced in \cite{breuillard2012structure}.}, deeply connected to the model{\hyp}theoretic structure of the piecewise definable group \cite[Theorem 4.2]{hrushovski2011stable}. Furthermore, among all these Lie groups, Hrushovski focused on the minimal one, showing its uniqueness and its independence of expansions of the language. 

Here, in Section 2, we improve these results by defining the more general notion of \emph{Lie core}. Then, we prove the existence of Lie cores for any piecewise hyperdefinable group with a generic piece and the uniqueness of the minimal Lie core. In the process, we adapt the classical model{\hyp}theoretic components (with parameters) $G^0$, $G^{00}$ and $G^{000}$ for piecewise hyperdefinable groups --- our definitions extend some particular cases already studied (e.g. \cite{hrushovski2022amenability}). We also introduce a new component $G^{\mathrm{ap}}$; $G^{\mathrm{ap}}$ is the smallest possible kernel of a continuous projection to a locally compact topological group without non{\hyp}trivial compact normal subgroups. We use these components to show that the minimal Lie core is precisely $\rfrac{G^0}{G^{\mathrm{ap}}}$. Using this canonical presentation of the minimal Lie core, we conclude that the minimal Lie core is piecewise $0${\hyp}hyperdefinable and independent of expansions of the language. \smallskip

Hyperdefinable sets, originally introduced in \cite{hart2000coordinatisation}, are quotients of $\bigwedge${\hyp}definable (read infinitely{\hyp}definable or type{\hyp}definable) sets over $\bigwedge${\hyp}definable equivalence relations. Hyperdefinable sets have been already well studied by different authors (e.g. \cite{wagner2010simple} and \cite{kim2014simplicity}). Here we extend this study to piecewise hyperdefinable sets. 

A piecewise hyperdefinable set is a strict direct limit of hyperdefinable sets. We are interested in the piecewise hyperdefinable sets as elementary objects in themselves that inherit an underlying model{\hyp}theoretic structure. In particular, we are principally interested in studying the natural logic topologies that generalise the usual Stone topology and can be defined for any piecewise hyperdefinable set. 

We will discuss the general theory of piecewise hyperdefinable sets in Section 1, elaborating the basic theory that we need for the rest of the paper. Mainly, we study their basic topological properties such as compactness, local compactness, normality, Hausdorffness and metrizability, and also their relations with the quotient, product and subspace topologies. All this study leads naturally to the definition of \emph{locally hyperdefinable sets}, which is in fact one of the main notions of this document. \medskip

The original motivation of this paper comes from the study of rough approximate subgroups. Rough approximate subgroups are the natural generalisation of approximate subgroups where we also allow a small thickening of the cosets. This generalisation is particularly natural in the context of metric groups, where the thickenings are given by balls. In \cite{tao2008product} (see also \cite{tao2014metric}), Tao showed that a significant part of the basic theory about approximate subgroups can be extended to rough approximate subgroups of metric groups via a discretisation process. In \cite{gowers2020partial}, this generalisation to the context of metric groups has recently shown interesting applications. 

It is a well{\hyp}known fact in first{\hyp}order model theory that, when working with metric spaces, non{\hyp}standard real numbers appear as soon as we get saturation. To solve this issue, we have to deal with some kind of continuous logic. Piecewise hyperdefinable sets provide a natural way of doing it. Roughly speaking, the idea behind continuous logic is to restrict the universe to the bounded elements and to quotient out by infinitesimals via the standard part function, as can be done using piecewise hyperdefinable sets. 

Hrushovski already indicated in unpublished works that, using piecewise hyperdefinable groups, it should be possible to extend some of the results of \cite{hrushovski2011stable} to the context of metric groups. The aim of this paper is to give the abstract basis for that kind of results, to find in the end possibly interesting applications to combinatorics. We conclude the paper giving the natural generalisation of the Lie model Theorem for rough approximate subgroups.  Applications of this result to the case of metric groups will be studied in a future paper. \medskip

We study in general piecewise hyperdefinable sets in Section 1, focussing on the properties of their logic topologies. The most important results of this section are given after the introduction of locally hyperdefinable sets in Section 1.4. Section 2 is the core of this paper and is devoted to the general study of piecewise hyperdefinable groups. The first fundamental result of the section is \cref{t:generic set lemma}, in which we show that piecewise hyperdefinable groups satisfying a natural combinatorial condition are locally hyperdefinable. In Section 2.3, we define the model{\hyp}theoretic components for piecewise hyperdefinable groups, proving the existence of $G^{\mathrm{ap}}$ in \cref{t:the aperiodic component}. Finally, we focus on the study of Lie cores, proving their existence (\cref{t:logic existence of lie core}), the uniqueness of the minimal one (\cref{t:logic uniqueness of minimal lie core}), giving a canonical representation of the minimal one in terms of the model{\hyp}theoretic components  (\cref{t:canonical representation of the minimal lie core}) and showing its independence of expansions of the language (\cref{c:independence of expansions of the minimal lie core}). Section 3 is devoted to the Stabilizer Theorem for piecewise hyperdefinable groups, which is divided over \cref{t:stabilizer theorem 1,c:stabilizer theorem mos b2,t:stabilizer theorem 2} --- \cref{t:stabilizer theorem 2} being the standard statement of the Stabilizer Theorem. We conclude the paper stating the Rough Lie model \cref{t:rough lie model} which generalises Hrushovski's Lie Model Theorem to the case of rough approximate subgroups.

\subsection*{Notations and conventions:} From now on, unless otherwise stated, fix a many{\hyp}sort first order language, $\lang$, and an infinite $\kappa${\hyp}saturated and strongly $\kappa${\hyp}homogeneous $\lang${\hyp}structure, $\mathfrak{M}$, with $\kappa>|\lang|$ a strong limit cardinal. We say that a subset of $\mathfrak{M}$ is \emph{small} if its cardinality is smaller than $\kappa$. 

We would like to remark that our assumptions on $\mathfrak{M}$ are far stronger than needed. While saturation is fundamental, strong homogeneity is actually irrelevant. Also, assuming that $\kappa$ is a strong limit cardinal is much more than necessary. Most of the results and arguments work when we only assume $\kappa>|\lang|$. The only significant exceptions are sections 2.3 and 2.4, where we need to assume $\kappa> 2^{|A|+|\lang|}$, with $A$ the set of parameters we are using.\smallskip

Against the common habit, we consider the logic topologies on the models themselves, rather than on the spaces of types. This has the unpleasant consequence of handling with non $\mathrm{T}_0$ topological spaces. In particular, we need to study normality and local compactness without Hausdorffness. 

Anyway, whether on the model or on the space of types, these two topologies are in fact the ``same". Indeed, the spaces of types are obtained from the logic topologies on the model after saturation (i.e. compactification) as the quotient by the equivalence relation of being topologically indistinguishable (i.e. the \emph{Kolmogorov quotient}). In particular, the canonical projections $\tp_A:\ a\mapsto \tp(a/A)$ are continuous, open and closed maps\footnote{In fact, $X\mapsto\{\tp(x/A)\sth x\in X\}$ is an isomorphism (in the sense of lattices of sets) between the topologies. Also, $\tp^{-1}_A[\tp_A[X]]=X$ for any $A${\hyp}invariant set.} from the $A${\hyp}logic topologies to the spaces of types over $A$. Moreover, based on this, we will show that piecewise hyperdefinable sets with their logic topologies generalise type spaces, the latter being the particular case of $\mathrm{T}_0$ logic topologies. 

We choose to work with logic topologies over the models for three reasons. The first one is that we save one unneeded step (the Kolmogorov quotient). The second one is that the spaces of types of products are not the products of the spaces of types, while we want to be able to prove and naturally use \cref{p:product logic topology}, which says that the global logic topology of the product is the product topology of the global logic topologies. The last one, and related with the previous one, is that the space of types of a group does not preserve usually the group structure, which makes odd the statement of the Isomorphism \cref{t:isomorphism theorem} for piecewise hyperdefinable groups, fundamental in Section 2, expressed trough the spaces of types.\smallskip
\begin{enumerate}[label={\raisebox{0.4ex}{\rule{0.4em}{0.4em}}}, itemsep=2pt, wide]
\item From now on, except when otherwise stated, we use $a,b,\ldots$ and $A,B,\ldots$ to denote hyperimaginaries and small sets of hyperimaginaries respectively, while $a^*,a^{**},\ldots$ and $A^*,A^{**},\ldots$ denote associated representatives. If we use $a^*$ or $A^*$ without mention of $a$ or $A$, we mean real elements. To avoid technical issues, we assume that $A$ contains a subset $A_{\mathrm{re}}\subseteq A$ of real elements such that every element of $A$ is hyperimaginary over $A_{\mathrm{re}}$.
\item We use $x,y,z,\ldots$ to denote variables. Except when otherwise stated, variables are always finite tuples of (single) variables. Sometimes, to simplify the notation, we also allow infinite tuples of variable. We write $\overline{x},\overline{y},\overline{z},\ldots$ to denote infinite tuples of variables. If in some exceptional place we need to consider single variables, we will explicitly indicate it.
\item Except when otherwise stated, elements are finite tuples of (single) elements. We write $\overline{a},\overline{b},\ldots$ and $\overline{a}^*,\overline{b}^*,\ldots$ to indicate when we consider infinite tuples of elements. In this paper, we do not study piecewise hyperdefinable sets of infinite tuples. It seems that our results can be generalised in that sense by taking also inverse limits.
\item By a type we always mean a complete type. The Stone space of types contained in the set $X$ with parameters $A^*$ is denoted by $\Stone_{X}(A^*)$. The set of formulas of $\lang$ with parameters in $A^*$ and variables in the (possibly infinite) tuple $\overline{x}$ is denoted by $\For^{\bar{x}}(\lang(A^*))$. The cardinality of the language is the cardinality of its set of formulas; $|\lang|\coloneqq |\For(\lang)|$. In particular, $|\lang|$ is always at least $\aleph_0$.
\item Here, definable means definable with parameters. Also, $\bigwedge${\hyp}definable means infinitely{\hyp}definable (or type{\hyp}definable) over a small set of parameters. Similarly, $\bigvee${\hyp}definable means coinfinitely{\hyp}definable (or cotype{\hyp}definable) over a small set of parameters. To indicate that we use parameters from $A$, we write $A${\hyp}definable, $\bigwedge_A${\hyp}definable and $\bigvee_A${\hyp}definable. We also use cardinals and cardinal inequalities. In that case, the subscript should be read as an anonymous set of parameters whose size satisfies the indicated condition. For example, $\bigwedge_{<\omega}${\hyp}definable means $\bigwedge_A${\hyp}definable for some subset $A$ with $|A|<\omega$. The same notation will be naturally used for hyperdefinable, piecewise hyperdefinable and piecewise $\bigwedge${\hyp}definable sets.
\item Following the terminology of \cite{tent2012course}, a $\kappa${\hyp}homogeneous structure is a structure such that every partial elementary internal map defined on a subset with cardinality less than $\kappa$ can be extended to any further element, while a strongly $\kappa${\hyp}homogeneous structure is a structure such that every partial elementary internal map defined on a subset with cardinality less than $\kappa$ can be extended to an automorphism. 
\item Let $R\subseteq X\times Y$ be a set{\hyp}theoretic binary relation. For $x\in X$, we write $R(x)\coloneqq \{y\in Y\sth (x,y)\in R\}$. For a subset $V\subseteq X$, we write $R[V]\coloneqq\{y\in Y\sth \exists x\in V\ (x,y)\in R\}$. We also write $R^{-1}\coloneqq\{(y,x)\sth (x,y)\in R\}$, so $R^{-1}(y)\coloneqq \{x\in X\sth (x,y)\in R\}$ for $y\in Y$ and $R^{-1}[W]\coloneqq \{x\in X\sth \exists y\in V\ (x,y)\in R\}$ for $W\subseteq Y$. We denote the image and preimage functions between the power sets by $\Ima\ R:\  V\mapsto R[V]$ and $\Ima^{-1}R:\ W\mapsto R^{-1}[W]$. Most of the time, this notation is used for partial functions, which are always identified with their graphs --- note  that $f^{-1}$ is only a function when $f$ is invertible.
\item Cartesian projections are denoted by $\proj$, quotient maps are denoted by $\quot$, quotient homomorphisms in groups are denoted by $\pi$, inclusion maps are denoted by $\inc$ and identity maps are denoted by $\id$.
\item A lattice of sets is a family of sets closed under finite unions and intersections. A complete algebra on a set $X$ is a family of subsets of $X$ closed under complements and arbitrary unions. 
\item The class of ordinals is denoted by $\Ord$. The cardinal of a set $X$ is $|X|$.
\item We use product notation for groups. Also, unless otherwise stated, we consider the group acting on itself on the left. In particular, by a coset we mean a left coset. A subset $X$ of a group is called symmetric if $1\in X=X^{-1}$. For subsets $X$ and $Y$ of a group, we write $XY$ for the set of pairwise products, and abbreviate $X^n\coloneqq XX^{n-1}$ and $X^{-n}\coloneqq (X^{-1})^n$ for $n\in\N$. We say that $X$ normalises $Y$ if $x^{-1}Yx\subseteq Y$ for every $x\in X$.
\item By a Lie group we always mean here a finite{\hyp}dimensional real Lie group.
\end{enumerate}

\section*{Acknowledgements} I want to specially thank my supervisor, Prof. Ehud Hrushovski, for all his help, ideas and suggestions: without him it would not have been possible to do this. I also want to thank my PhD mate, Alex Chevalier, for his helpful reviews and comments. Finally, I want to thank the anonymous referee of the JML for the thoughtful review and many comments.

\section{Piecewise hyperdefinable sets}
\subsection{Hyperdefinable sets}
Let $A^*$ be a small set of parameters. An \emph{$A^*${\hyp}hyperdefinable} set is a quotient $P=\rfrac{X}{E}$ where $X$ is a non{\hyp}empty $\bigwedge_{A^*}${\hyp}definable set and $E$ is an $\bigwedge_{A^*}${\hyp}definable equivalence relation. If we do not indicate the set of parameters, we mean that it is hyperdefinable for some small set of parameters. Write $\quot_P:\ X\rightarrow P$ for the \emph{quotient map} given by $x\mapsto [x]_E\coloneqq \rfrac{x}{E}\coloneqq E(x)=\{x'\in X\sth (x,x')\in E\}$. The elements of $A${\hyp}hyperdefinable sets are called \emph{hyperimaginaries over $A$}. Given a hyperimaginary element $a$, a \emph{representative} of $a$ is an element $a^*\in a$ of the structure such that $\quot_P(a^*)=a$. The elements of the structure will be called \emph{real}. 

An \emph{$\bigwedge_{A^*}${\hyp}definable} subset, $V\subseteq P$, is a subset such that $\quot^{-1}_P[V]$ is $\bigwedge_{A^*}${\hyp}definable in $X$. We will say that a partial type defines $V\subseteq P$ if it defines $\quot^{-1}_P[V]$. If $V\subseteq P$ is a non{\hyp}empty $\bigwedge${\hyp}definable set, after declaring the parameters, we will write $\underline{V}$ to denote a partial type defining $\quot^{-1}_P[V]$. The following basic proposition is the starting point to study hyperdefinable sets.\smallskip

\begin{lem}[Correspondence Lemma] \label{l:correspondence lemma} Let $P=\rfrac{X}{E}$ be an $A^*${\hyp}hyperdefinable set. Then, the image by $\quot_P$ of any $\bigwedge_{A^*}${\hyp}definable subset of $X$ is an $\bigwedge_{A^*}${\hyp}definable subset of $P$. Moreover, the preimage function $\Ima^{-1}\quot_P$ is an isomorphism, whose inverse is $\Ima\ \quot_P$, between the lattice of $\bigwedge_{A^*}${\hyp}definable subsets of $P$ and the lattice of $\bigwedge_{A^*}${\hyp}definable subsets of $X$ closed under $E$.
\end{lem}
The main part of this lemma can be further generalised to \cref{l:infinite definable functions with real parameters}. For this, note firstly that, given two $A^*${\hyp}hyperdefinable sets $P=\rfrac{X}{E}$ and $Q=\rfrac{Y}{F}$, the Cartesian product $P\times Q$ is canonically identified with the hyperdefinable set $\rfrac{X\times Y}{E\hat{\times} F}$ via $([x]_E,[y]_F)\mapsto [x,y]_{E\hat{\times} F}$, where $E\hat{\times} F\coloneqq\{((x,y),(x',y'))\sth (x,x')\in E$, $(y,y')\in F\}$. Then, we can talk about $\bigwedge_{A^*}${\hyp}definable relations and partial functions. 
\begin{ex} The inclusion $\inc:\ V\rightarrow P$ of an $\bigwedge_{A^*}${\hyp}definable set $V$ is $\bigwedge_{A^*}${\hyp}definable. The Cartesian projections $\proj_P:\ P\times Q\rightarrow P$ and $\proj_Q:\ P\times Q\rightarrow Q$ are $\bigwedge_{A^*}${\hyp}definable. Also, the quotient map $\quot_P:\ X\rightarrow P$ is $\bigwedge_{A^*}${\hyp}definable. 
\end{ex}
\begin{lem} \label{l:infinite definable functions with real parameters} 
Let $P=\rfrac{X}{E}$ and $Q=\rfrac{Y}{F}$ be two $A^*${\hyp}hyperdefinable sets and $f$ an $\bigwedge_{A^*}${\hyp}definable partial function from $P$ to $Q$. Then, for any $\bigwedge_{A^*}${\hyp}definable sets $V\subseteq P$ and $W\subseteq Q$, $f[V]$ and $f^{-1}[W]$ are $\bigwedge_{A^*}${\hyp}definable. 
\begin{proof} It is enough to check that the projection maps satisfy the proposition. It is trivial that $\proj^{-1}_P[V]=V\times Q$ is $\bigwedge_{A^*}${\hyp}definable for any $\bigwedge_{A^*}${\hyp}definable subset $V$. On the other hand, if $V\subseteq P\times Q$ is $\bigwedge_{A^*}${\hyp}definable, by compactness, $\Sigma(x)=\{\exists y\mathrel{} \bigwedge \Delta(x,y)\sth \Delta\subseteq \underline{V}\mbox{ finite}\}$ defines $\quot^{-1}_P[\proj_P[V]]$. 
\end{proof}
\end{lem}
\begin{rmk} Let $f:\ P\rightarrow Q$ and $g:\ Q\rightarrow R$ be functions. As subsets of $P\times R$, we have $g\circ f=\proj_{P\times R}[(f\times R)\cap (P\times g)]$, where $\proj_{P\times R}:\ P\times Q\times R\rightarrow P\times R$ is the natural projection. Thus, compositions of $\bigwedge_{A^*}${\hyp}definable partial functions are also $\bigwedge_{A^*}${\hyp}definable partial functions. 
\end{rmk}

A \emph{type over $A^*$, or $A^*${\hyp}type}, in $P$ is an $\bigwedge_{A^*}${\hyp}definable subset of $P$ which is $\subset${\hyp}minimal in the family of non{\hyp}empty $\bigwedge_{A^*}${\hyp}definable subsets. For  $a\in P$, we write $\tp(a/A^*)$ for the type over $A^*$ containing $a$. As the lattice of $\bigwedge_{A^*}${\hyp}definable sets is closed under arbitrary intersections, the type of a hyperimaginary element always exists.  

As usual, for infinite tuples $\overline{a}=(a_t)_{t\in T}$ and $\overline{b}=(b_t)_{t\in T}$ of hyperimaginaries over $A^*$, we write $\tp(\overline{a}/A^*)=\tp(\overline{b}/A^*)$ to mean that $\tp(\overline{a}_{\mid T_0}/A^*)=\tp(\overline{b}_{\mid T_0}/A^*)$ for any $T_0\subseteq T$ finite.
\begin{lem} \label{l:types of hyperimaginaries} 
Let $P$ be an $A^*${\hyp}hyperdefinable set, $a\in P$ and $a^*\in a$. Then, $\tp(a/A^*)=\quot_P[\tp(a^*/A^*)]$. In particular, $b\in\tp(a/A^*)$ if and only if there is $b^*\in b$ such that $\tp(a^*/A^*)=\tp(b^*/A^*)$. In other words, $\tp(a/A^*)$ is the orbit of $a$ under the action of $\Aut(\mathfrak{M}/A^*)$ on $P$.
\begin{proof} By the Correspondence \cref{l:correspondence lemma} and minimality. 
\end{proof}
\end{lem}
Let $P$ be $A^*${\hyp}hyperdefinable. A subset $V\subseteq P$ is \emph{$A^*${\hyp}invariant} if $\tp(a/A^*)\subseteq V$ for any $a\in V$. As immediate corollaries of \cref{l:types of hyperimaginaries} we get the following results:
\begin{coro} \label{c:invariance} 
Let $P$ be $A^*${\hyp}hyperdefinable and $V\subseteq P$. Then, $V$ is $A^*${\hyp}invariant if and only if it is setwise invariant under the action of $\Aut(\mathfrak{M}/A^*)$ on $P$, if and only if $\quot^{-1}[V]$ is $A^*${\hyp}invariant.
\begin{proof} Obvious by \cref{l:types of hyperimaginaries}. 
\end{proof}
\end{coro}
\begin{coro} \label{c:parameters of definition} 
Let $P$ be $A^*${\hyp}hyperdefinable and $V\subseteq P$ an $\bigwedge${\hyp}definable subset. Then, $V$ is $\bigwedge_{A^*}${\hyp}definable if and only if it is $A^*${\hyp}invariant.
\end{coro} 
We also want to be able to use hyperimaginary parameters. To do that we need to redefine the previous notions. Let $A$ be a set of hyperimaginaries. We start by defining types over $A$ and $A${\hyp}invariance for real elements. Then, we can say what is an $A${\hyp}hyperdefinable set, an $A${\hyp}invariant set, an $\bigwedge_A${\hyp}definable set and a type over $A$ for hyperimaginaries.

Since we assume that every element in $A$ is hyperimaginary over $A_{\mathrm{re}}$, the group $\Aut(\mathfrak{M}/A_{\mathrm{re}})$ naturally acts on the elements of $A$, i.e. $\sigma(a)$ makes sense for $a\in A$ and $\sigma\in\Aut(\mathfrak{M}/A_{\mathrm{re}})$. The \emph{group of automorphisms $\Aut(\mathfrak{M}/A)$ pointwise fixing $A$} is defined as the subgroup of elements of $\Aut(\mathfrak{M}/A_{\mathrm{re}})$ fixing each element of $A$. Note that, obviously, $\Aut(\mathfrak{M}/A^*)\leq \Aut(\mathfrak{M}/A)$ for any set of representatives $A^*$ of $A$. 
 
Let $a^*$ be a tuple of real elements. The type of $a^*$ over $A$ is the set $\tp(a^*/A)\coloneqq\{b^*\sth \tp(a^*,A)=\tp(b^*,A)\}$. Equivalently, by \cref{l:types of hyperimaginaries}, $\tp(a^*/A)$ is the orbit of $a^*$ under $\Aut(\mathfrak{M}/A)$. A set of real elements $X$ is \emph{$A${\hyp}invariant} if we have $\tp(a^*/A)\subseteq X$ for any $a^*\in X$. Equivalently, $X$ is $A${\hyp}invariant if and only if it is setwise invariant under the action of $\Aut(\mathfrak{M}/A)$. We say that a hyperdefinable set $P=\rfrac{X}{E}$ is \emph{$A${\hyp}hyperdefinable} if $X$ and $E$ are $A${\hyp}invariant. Note that, as $\Aut(\mathfrak{M}/A^*)\leq \Aut(\mathfrak{M}/A)$, an $A${\hyp}hyperdefinable set is, in particular, $A^*${\hyp}hyperdefinable for any set of representatives $A^*$ of $A$. If $P$ is $A${\hyp}hyperdefinable, it follows that $\Aut(\mathfrak{M}/A)$ naturally acts on $P$.

A subset of an $A${\hyp}hyperdefinable set is \emph{$A${\hyp}invariant} if its preimage by the quotient map is $A${\hyp}invariant. Equivalently, it is $A${\hyp}invariant if and only if it is setwise invariant under the action of $\Aut(\mathfrak{M}/A)$. By \cref{c:invariance}, as $\Aut(\mathfrak{M}/A^*)\leq \Aut(\mathfrak{M}/A)$, every $A${\hyp}invariant set is in particular $A^*${\hyp}invariant for any set of representatives. An \emph{$\bigwedge_A${\hyp}definable} subset is an $A${\hyp}invariant $\bigwedge${\hyp}definable subset. By \cref{c:parameters of definition}, every $\bigwedge_A${\hyp}definable set is $\bigwedge_{A^*}${\hyp}definable for any set of representatives. Note that the $A${\hyp}invariant subsets of $P$ form a complete algebra of subsets of $P$. Thus, $\bigwedge_A${\hyp}definable sets are a lattice of sets closed under arbitrary intersections. 
\begin{lem} \label{l:infinite definable functions} 
Let $P=\rfrac{X}{E}$ and $Q=\rfrac{Y}{F}$ be two $A${\hyp}hyperdefinable sets and $f$ an $\bigwedge_A${\hyp}definable partial function. Then, for any $\bigwedge_A${\hyp}definable sets $V\subseteq P$ and $W\subseteq Q$, $f[V]$ and $f^{-1}[W]$ are $\bigwedge_A${\hyp}definable. 
\begin{proof} By \cref{l:infinite definable functions with real parameters}, $f[V]$ and $f^{-1}[W]$ are $\bigwedge${\hyp}definable. As $f$ is $A${\hyp}invariant, $f(\sigma(a))=\sigma(f(a))$ for $a\in P$ and $\sigma\in\Aut(\mathfrak{M}/A)$, so $f[V]$ and $f^{-1}[W]$ are $A${\hyp}invariant. 
\end{proof}
\end{lem}

A \emph{type over $A$, or $A${\hyp}type}, is a $\subset${\hyp}minimal non{\hyp}empty $\bigwedge_A${\hyp}definable subset. For $a\in P$, we write $\tp(a/A)$ for the type over $A$ containing $a$. As the lattice of $\bigwedge_{A}${\hyp}definable sets is closed under arbitrary intersections, the type of a hyperimaginary element always exists. As usual, for infinite tuples $\overline{a}=(a_t)_{t\in T}$ and $\overline{b}=(b_t)_{t\in T}$, we write $\tp(\overline{a}/A)=\tp(\overline{b}/A)$ to mean that $\tp(\overline{a}_{\mid T_0}/A)=\tp(\overline{b}_{\mid T_0}/A)$ for any $T_0\subseteq T$ finite.
\begin{lem} \label{l:types and infinite definable functions}
Let $P,Q$ be $A${\hyp}hyperdefinable sets, $f:\ P\rightarrow Q$ an $\bigwedge_{A}${\hyp}definable function and $a\in P$. Then, $f[\tp(a/A)]=\tp(f(a)/A)$. 
\begin{proof} By \cref{l:infinite definable functions}, we have $f[\tp(a/A)]$ and $f^{-1}[\tp(f(a)/A)]$ are $\bigwedge_A${\hyp}definable. Then, $\tp(a/A)\subseteq f^{-1}[\tp(f(a)/A)]$ and $\tp(f(a)/A)\subseteq f[\tp(a/A)]$ by minimality. Thus, $f[\tp(a/A)]=\tp(f(a)/A)$. 
\end{proof}
\end{lem}
As an immediate corollary we get the following result:
\begin{coro}\label{c:types over hyperimaginaries}
Let $P$ be an $A${\hyp}hyperdefinable set, $a\in P$ and $a^*\in a$. Then, $\tp(a/A)=\quot_P[\tp(a^*/A)]$. In particular, $\tp(a/A)$ is the orbit of $a$ under the action of $\Aut(\mathfrak{M}/A)$ on $P$. In other words, for any $A^*$ representatives of $A$, we have $b\in\tp(a/A)$ if and only if there are $b^{**}\in b$ and $A^{**}$ representatives of $A$ such that $\tp(a^*,A^*)=\tp(b^{**},A^{**})$. 

Consequently, $V$ is $A${\hyp}invariant if and only if $\tp(a/A)\subseteq V$ for any $a\in V$. 
\begin{proof} Clear by \cref{l:types and infinite definable functions,l:types of hyperimaginaries}. 
\end{proof}
\end{coro}


We explain now how to substitute hyperimaginary parameters. Let $P$ be $A${\hyp}hyperdefinable, $\overline{b}$ a small set of hyperimaginaries over $A$ and $V\subseteq P$ an $\bigwedge_{\overline{b}}${\hyp}definable subset. Let $\overline{c}$ be such that $\tp(\overline{b}/A)=\tp(\overline{c}/A)$. The set \emph{$V(\overline{c})$ given from $V$ by replacing $\overline{b}$ by $\overline{c}$} is the set $\sigma[V]$ for $\sigma\in\Aut(\mathfrak{M}/A)$ such that $\sigma(\overline{b})=\overline{c}$. Note that this does not depend on the choice of $\sigma\in\Aut(\mathfrak{M}/A)$.

Alternatively, we present a more explicit methodology using \emph{uniform definitions} that does not require the use of automorphisms. Let $P=\rfrac{X}{E}$ be an $A${\hyp}hyperdefinable set. We say that a partial type $\Sigma(x,A^*)$ is \emph{weak uniform on $P$ over $A$} if $\Sigma(x,A^{**})$ defines the same set on $P$ for any set of representatives $A^{**}$ of $A$, i.e. $\rfrac{\Sigma(\mathfrak{M},A^*)}{E}=\rfrac{\Sigma(\mathfrak{M},A^{**})}{E}$. We say that $\Sigma(x,A^*)$ is \emph{uniform on $P$ over $A$} if $\Sigma(x,A^*)\cap \For^x(\lang(B^*))$ is weak uniform on $P$ over $B$ for any $B^*\subseteq A^*$ such that $P$ is still $B${\hyp}hyperdefinable and $A_{\mathrm{re}}\subseteq B$. Let $V\subseteq P$ be $\bigwedge_A${\hyp}definable. A \emph{uniform definition of $V$ over $A$} is a partial type $\underline{V}$ uniform on $P$ over $A$ that defines $V$. 
\begin{lem} \label{l:uniform definition} 
Let $P=\rfrac{X}{E}$ be an $A${\hyp}hyperdefinable set and $V\subseteq P$ a non{\hyp}empty $\bigwedge_A${\hyp}definable set. Then, there is a uniform definition of $V$ over $A$.  
\begin{proof} Say $A=\{a_i\}_{i\in I}$ with $a_i\in Q_i=\rfrac{Y_i}{F_i}$. Pick $A^*$ representatives of $A$ and $\Sigma(x,A^*)$ partial type defining $V$ over $A^*$. Write $\underline{F}=\bigwedge_i \underline{F}_i$ and $\Gamma=\tp(A^*)$. Using saturation, take the partial type $\underline{V}(x,A^*)$ expressing $\exists \overline{y}\mathrel{} \Sigma(x,\overline{y})\wedge \underline{F}(\overline{y},A^*)\wedge \Gamma(\overline{y})$. By \cref{c:types over hyperimaginaries}, it is a uniform definition of $V$ over $A$. 
\end{proof}
\end{lem}

\begin{lem}\label{l:substitution and uniform definitions}
Let $P$ be an $A${\hyp}hyperdefinable set and $V\subseteq P$ be a non{\hyp}empty $\bigwedge_{\overline{b}}${\hyp}definable subset with $\overline{b}$ a small set of hyperimaginaries over $A$ and $A\subseteq \overline{b}$. Let $\underline{V}(x,\overline{b}^*)$ be a uniform definition of $V$. Let $\overline{c}$ be such that $\tp(\overline{b}/A)=\tp(\overline{c}/A)$ and $\overline{c}^*$ be representatives. Then, $\underline{V}(x,\overline{c}^*)$ is a uniform definition of $V(\overline{c})$. 
\begin{proof} Take $\sigma\in\Aut(\mathfrak{M}/A)$ with $\sigma(\overline{b})=\overline{c}$, so $\overline{b}^{**}=\sigma^{-1}(\overline{c}^*)$ is a representative of $\overline{b}$. As $\underline{V}(x,\overline{b}^*)$ is a uniform definition of $V$ on $P$ over $b$, we have that $\underline{V}(x,\overline{b}^{**})$ also defines $V$. Consequently, $V(\overline{c})$ is defined by $\underline{V}(x,\overline{c}^*)$. Now, take $\overline{c}_0\subseteq\overline{c}$ such that $P$ is still $ \overline{c}_0${\hyp}hyperdefinable and say $\overline{c}_0=\sigma(\overline{b}_0)$. Then, $P$ is still $\overline{b}_0${\hyp}hyperdefinable and, as $\underline{V}(x,\overline{b}^*)$ is uniform on $P$ over $b$, we have that $\underline{W}(x,\overline{b}^{**}_0)=\underline{V}(x,\overline{b}^{**})\cap \For(\lang(\overline{b}^{**}_0))$ is weakly uniform on $P$ over $\overline{b}_0$. Therefore, as $\tp(c^*_0/A)=\tp(b^{**}_0/A)$, $\underline{W}(x,\overline{c}^{*}_0)$ is weakly uniform on $P$ over $\overline{c}_0$. 
\end{proof}
\end{lem}
\begin{lem} \label{l:uniform definition finite} 
Let $P$ and $Q$ be $A${\hyp}hyperdefinable sets and $b\in Q$. Then, for any $\bigwedge_{A,b}${\hyp}definable subset $V\subseteq P$ there is an $\bigwedge_A${\hyp}definable set $W\subseteq Q\times P$ such that $V(c)=W(c)\coloneqq \proj_P[W\cap (\{c\}\times P)]$ for any $c\in \tp(b/A)$.
\begin{proof} Take a uniform definition $\Sigma(A^*,b^*,x)$ of $V$ on $P$ over $A$. By \cref{l:substitution and uniform definitions}, we get that $\Sigma(A^*,y,x)$ defines an $\bigwedge_A${\hyp}definable subset $W$ of $Q\times P$ such that $V(c)=W(c)\coloneqq \proj_P[W\cap (\{c\}\times P)]$ for any $c\in\tp(b/A)$. 
\end{proof}
\end{lem}
\begin{lem} \label{l:equivalence relation of equality of types} 
Let $P$ be an $A${\hyp}hyperdefinable set. Then,  
\[\Delta_P(A)=\{(a,b)\in P\times P\sth \tp(a/A)=\tp(b/A)\}\]
is an $\bigwedge_{A}${\hyp}definable equivalence relation. Furthermore, it has a uniform definition $\underline{\Delta_P}(A^*)$ such that $\underline{\Delta_P}(A^*)\cap \For(\lang(B^*))$ defines $\Delta_P(B)$ for any subset $B\subseteq A$ such that $P$ is $B${\hyp}hyperdefinable.
\end{lem}
\subsection{The logic topologies of hyperdefinable sets}
Let $P=\rfrac{X}{E}$ be an $A${\hyp}hyperdefinable subset. The \emph{$A${\hyp}logic topology} of $P$ is the one given by taking as closed sets the $\bigwedge_A${\hyp}definable subsets of $P$. In particular, by \cref{l:infinite definable functions}, the $A${\hyp}logic topology of $P$ is the quotient topology of the $A${\hyp}logic topology of $X$. 
\begin{prop}
Let $P$ and $Q$ be $A${\hyp}hyperdefinable sets and $f:\ P\rightarrow Q$ an $\bigwedge_A${\hyp}definable function. Then, $f$ is a continuous and closed function between the $A${\hyp}logic topologies.
\begin{proof} By \cref{l:infinite definable functions}. 
\end{proof}
\end{prop}
\begin{prop}\label{p:topology hyperdefinables} Let $P$ be an $A${\hyp}hyperdefinable set. Then, the $A${\hyp}logic topology of $P$ is compact and normal (i.e. any two disjoint closed sets can be separated by open sets). The closure of a point $a\in P$ is $\tp(a/A)$, so the properties $\mathrm{T}_0$, $\mathrm{T}_1$ and $\mathrm{T}_2$ are equivalent to $\tp(a/A)=\{a\}$ for all $a\in P$.
\begin{proof} Compactness follows from saturation. For normality we should note that the image of a normal space by a continuous closed function is always normal. Therefore, using $\quot_P$, from normality of the $A^*${\hyp}logic topology of $X$ we conclude the normality of the $A^*${\hyp}logic topology of $P$. Using that $a\mapsto \rfrac{a}{\Delta_P(A)}$ is continuous and closed from the $A^*${\hyp}logic topology of $P$ to the $A${\hyp}logic topology of $\rfrac{P}{\Delta_P(A)}$, we get normality of the latter. From there, we trivially conclude normality of the $A${\hyp}logic topology of $P$. 

Finally, as the closure of a point is its type by definition, $\mathrm{T}_1$ is equivalent to $\tp(a/A)=\{a\}$ for every $a\in P$. By \cref{c:types over hyperimaginaries}, $a\in\tp(b/A)$ if and only if $\tp(a/A)=\tp(b/A)$, so $b\in \overline{\{a\}}$ if and only if $a\in \overline{\{b\}}$ --- this topological property is sometimes called $\mathrm{R}_0$. Therefore, $\mathrm{T}_0$ is also equivalent to $\tp(a/A)=\{a\}$ for every $a\in P$. By normality, $\mathrm{T}_1$ and $\mathrm{T}_2$ are equivalent. 
\end{proof}
\end{prop} 
\begin{prop} \label{p:uniqueness of hausdorff logic topologies} 
Let $P$ be an $A${\hyp}hyperdefinable set and $A\subseteq A'$. If the $A${\hyp}logic topology of $P$ is Hausdorff, then the $A${\hyp}logic topology and the $A'${\hyp}logic topology are equal. Thus, there is at most one Hausdorff logic topology and, if it exists, it will be called the \emph{global logic topology}. 
\end{prop}
\begin{rmk} Furthermore, the global logic topologies are preserved by expansion of the language. Indeed, let $P$ be a hyperdefinable set whose $A${\hyp}logic topology is Hausdorff and $\mathfrak{M}'$ a $\kappa${\hyp}saturated and strongly $\kappa${\hyp}homogeneous expansion of $\mathfrak{M}$. Then, obviously, $\mathrm{id}:\ P\rightarrow P$ is a continuous bijection from the $A${\hyp}logic topology in $\mathfrak{M}'$ to the $A${\hyp}logic topology in $\mathfrak{M}$. As both are compact and the second one is Hausdorff, we conclude that both are the same.
\end{rmk}
Let $P$ be an $A${\hyp}hyperdefinable set. Assuming that $\kappa> 2^{|A|+|\lang|}$, either $|P|\geq \kappa$ or $|P|\leq 2^{|A|+|\lang|}$. Then, $P$ has a global logic topology if and only if $P$ is small, if and only if $|P|\leq 2^{|A|+|\lang|}$. Furthermore, write $\mathrm{bdd}(A)$ for the set of \emph{bounded hyperimaginaries over $A$}, i.e. hyperimaginaries over $A$ such that $|\tp(a/A)|<\kappa$. Note that $|\mathrm{bdd}(A)|\leq 2^{|A|+|\lang|}$. It follows that, for any $A${\hyp}hyperdefinable set $P$, it has a global logic topology if and only if the $\mathrm{bdd}(A)${\hyp}logic topology is the global logic topology.  
\begin{prop} \label{p:product logic topology} 
Let $P$ and $Q$ be $A${\hyp}hyperdefinable sets. Then, the $A${\hyp}logic topology in $P\times Q$ is at least as fine as the product topology of the $A${\hyp}logic topologies. Furthermore, $P\times Q$ has a global logic topology if and only if $P$ and $Q$ have so, and then the global logic topology of $P\times Q$ is the product topology of the global logic topologies of $P$ and $Q$.
\end{prop}
\subsection{Piecewise hyperdefinable sets}
A \emph{piecewise $A${\hyp}hyperdefinable} set is a strict direct limit of $A${\hyp}hyperdefinable sets with $\bigwedge_A${\hyp}definable inclusions. In other words, a piecewise $A${\hyp}hyperdefinable set is a direct limit 
\[P\coloneqq \underrightarrow{\lim}_{(I,\prec)}(P_i,\varphi_{ji})_{j\succeq i}=\rfrac{\bigsqcup_I P_i}{\sim_P}\]
of a direct system of $A${\hyp}hyperdefinable sets and $1${\hyp}to{\hyp}$1$ $\bigwedge_A${\hyp}definable functions $\varphi_{ji}:\ P_i\rightarrow P_j$.  

Recall that $(I,\prec)$ is a direct ordered set and, for each $i,j,k\in I$ with $i\preceq j\preceq k$, $\varphi_{kj}\circ\varphi_{ji}=\varphi_{kj}$ and $\varphi_{ii}=\mathrm{id}$. Also, recall that $\sim_{P}$ is the equivalence relation on $\bigsqcup_I P_i$ defined by $x\sim_{P}y$ for $x\in P_i$ and $y\in P_j$ if and only if there is some $k\in I$ with $i\preceq k$ and $j\preceq k$ and $\varphi_{ki}(x)=\varphi_{kj}(y)$. Note that, in fact, as we consider only direct systems where the functions $\varphi_{ji}$ are $1${\hyp}to{\hyp}$1$, the equivalence relation is given by $x\sim_{P}y$ with $x\in P_i$ and $y\in P_j$ if and only if $\varphi_{ki}(x)=\varphi_{kj}(y)$ for any $k\in I$ such that $i\preceq k$ and $j\preceq k$. 

The \emph{pieces} of $P$ are the subsets $\rfrac{P_i}{\sim_{P}}$. The \emph{canonical inclusions} are the maps $\inc_{P_i}:\ P_i\rightarrow \rfrac{P_i}{\sim_P}\subseteq P$ given by $a\mapsto [a]_{\sim_P}$ for $a\in P_i$.   

The \emph{cofinality} $\cf(P)$ of $P$ is the cofinality of $I$, which is the minimal ordinal $\alpha$ from which there is a function $f:\ \alpha\rightarrow I$ such that for every $i\in I$ there is $\xi\in \alpha$ with $i\preceq f(\xi)$. We say that $P$ is \emph{countably} piecewise hyperdefinable if it has countable cofinality. From now on, we always assume that $\cf(P)<\kappa$.

A \emph{piecewise $\bigwedge_{A}${\hyp}definable} subset of $P$ is a subset $V\subseteq P$ such that $\inc^{\ -1}_{P_i}[V]$ is $\bigwedge_{A}${\hyp}definable in $P_i$ for each $i\in I$. An \emph{$\bigwedge_A${\hyp}definable} subset of $P$ is a piecewise $\bigwedge_A${\hyp}definable subset contained in some piece. If $V\subseteq P$ is a non{\hyp}empty $\bigwedge${\hyp}definable set, after fixing a piece $\rfrac{P_i}{\sim_P}$ containing it and declaring the parameters, we will write $\underline{V}$ to denote a partial type defining it in $P_i$. In that case, we will also say that $\underline{V}$ defines $V$. If $a\in P$ is an element, after fixing a piece $\rfrac{P_i}{\sim_P}$ containing it, we will say that $a^*$ is a representative of $a$ if $\inc_{P_i}(\quot_{P_i}(a^*))=a$, where $\quot_{P_i}$ is the quotient map of $P_i$ as hyperdefinable set. 

\begin{rmk} The set of piecewise $\bigwedge_A${\hyp}definable subsets is a lattice of sets closed under arbitrary intersections, and the collection of $\bigwedge_A${\hyp}definable subsets is the ideal of that lattice generated by the pieces.\end{rmk}
 
Note that, in the previous definitions, piecewise hyperdefinable sets are not only sets but sets together with a particular structure. This structure is given by the lattices of piecewise $\bigwedge${\hyp}definable subsets and the ideals of $\bigwedge${\hyp}definable subsets. It is very important to remember this as the same set could be represented as a piecewise hyperdefinable set in several different ways --- see \cref{e:example 1}. \medskip

By the strictness condition we get the following fundamental lemma. 
\begin{lem}[Correspondence Lemma]\label{l:correspondence lemma in pieces} 
Let $P$ be piecewise $A${\hyp}hyperdefinable and $\rfrac{P_i}{\sim_P}$ a piece of $P$. Then, $\inc_{P_i}:\ P_i\rightarrow \rfrac{P_i}{\sim_P}$ is a bijection. Furthermore,  $V\subseteq \rfrac{P_i}{\sim_P}$ is $\bigwedge_A${\hyp}definable as subset of $P$ if and only if $\inc_{P_i}^{\ -1}[V]$ is $\bigwedge_A${\hyp}definable as subset of $P_i$.
\begin{proof} Obviously $\inc_{P_i}:\ P_i\rightarrow\rfrac{P_i}{\sim_P}$ is a bijection by the strictness condition. Say $V\subseteq \rfrac{P_i}{\sim_P}$. Using again the strictness condition, for any $j,k\in I$ with $i,j\preceq k$, we have that $\inc^{\ -1}_{P_j}[V]=\varphi_{kj}^{-1}[\inc^{\ -1}_{P_k}[V]]=\varphi_{kj}^{-1}[\varphi_{ki}[\inc^{\ -1}_{P_i}[V]]]$. Therefore, by $\bigwedge_A${\hyp}definability of the maps and \cref{l:infinite definable functions}, we conclude that $V$ is $\bigwedge_A${\hyp}definable if and only if $\inc_{P_i}^{\ -1}[V]$ is $\bigwedge_A${\hyp}definable as subset of $P_i$. 
\end{proof}
\end{lem}
The Correspondence \cref{l:correspondence lemma in pieces} says that $\inc_{P_i}$ is a true identification between $P_i$ and $\rfrac{P_i}{\sim_P}$ in terms of the model theoretic structure. In other words, it says that pieces of piecewise hyperdefinable sets are, indeed, hyperdefinable. From now on, slightly abusing of the notation, we make no distinction between $P_i$ and $\rfrac{P_i}{\sim_P}$, i.e. $P_i\coloneqq\rfrac{P_i}{\sim_P}$ and $\quot_{P_i}\coloneqq \inc_{P_i}\circ \quot_{P_i}$. 
\begin{coro} \label{c:infinite definable in piecewise hyperdefinable} 
Let $P=\underrightarrow{\lim}_I P_i$ be piecewise $A${\hyp}hyperdefinable. Then, a subset $V\subseteq P$ is piecewise $\bigwedge_A${\hyp}definable if and only if $V\cap P_i$ is $\bigwedge_A${\hyp}definable as subset of $P$ for each $i\in I$.
\end{coro}
\begin{rmk}\label{o:remark ideal infinite definable sets} As every piece is $\bigwedge${\hyp}definable, the lattice of piecewise $\bigwedge_A${\hyp}definable subsets can be recovered from the ideal of $\bigwedge_A${\hyp}definable subsets. In other words, the structure of $P$ is completely determined by the ideals of $\bigwedge${\hyp}definable sets. However, in general, the lattice of piecewise $\bigwedge_A${\hyp}definable subsets does not determine the ideal of $\bigwedge_A${\hyp}definable subsets --- see \cref{e:example 1}.
\end{rmk}

A \emph{type over $A$, or $A${\hyp}type}, is a $\subset${\hyp}minimal non{\hyp}empty piecewise $\bigwedge_A${\hyp}definable subset. The \emph{type of $a\in P$ over $A$}, $\tp(a/A)$, is the minimal piecewise $\bigwedge_A${\hyp}definable subset of $P$ containing $a$. Since $a\in P_i$ for some $i\in I$ and pieces are $\bigwedge_A${\hyp}definable, $\tp(a/A)$ is actually $\bigwedge_A${\hyp}definable.\smallskip 

Let $R\subseteq P\times Q$ be a binary relation between two piecewise hyperdefinable sets. We say that $R$ is \emph{piecewise bounded (or piecewise continuous)} if the image of any piece of $P$ is contained in some piece of $Q$. We say that it is \emph{piecewise proper} if the preimage of any piece of $Q$ is contained in some piece of $P$. We use this terminology in particular for partial functions. To simplify the terminology, we often omit reiterative uses of ``piecewise'' when they happen. So, for example, we say ``a piecewise bounded and proper $\bigwedge${\hyp}definable function'' instead of ``a piecewise bounded and piecewise proper piecewise $\bigwedge${\hyp}definable function''. 

Given two piecewise $A${\hyp}hyperdefinable sets $P=\underrightarrow{\lim}_{I}P_i$ and $Q=\underrightarrow{\lim}_{J}Q_j$, the Cartesian product $P\times Q$ is canonically identified with $\underrightarrow{\lim}_{I\times J}P_i\times Q_j$ via the map $([x]_{P},[y]_{Q})\mapsto [x,y]_{P\times Q}$, where $(x,y)\sim_{P\times Q}(x',y')$ if and only if $x\sim_Px'$ and $y\sim_Q y'$. Thus, we say that a binary relation $R$ between $P$ and $Q$ is piecewise $\bigwedge_A${\hyp}definable if it is so as subset of the Cartesian product. 
\begin{lem} \label{l:types and piecewise infinite definable functions}
Let $P$ and $Q$ be piecewise $A${\hyp}hyperdefinable sets, $f:\ P\rightarrow Q$ a piecewise $\bigwedge_{A}${\hyp}definable function and $a\in P$. Then, $f[\tp(a/A)]=\tp(f(a)/A)$.
\begin{proof} Clear from \cref{l:types and infinite definable functions}. 
\end{proof}
\end{lem}
\begin{prop} \label{p:piecewise infinite definable functions} 
Let $P$ and $Q$ be two piecewise $A${\hyp}hyperdefinable sets and $f$ a piecewise $\bigwedge_A${\hyp}definable partial function from $P$ to $Q$. Then:
\begin{enumerate}[label={\rm{(\arabic*)}}, ref={\rm{\arabic*}}, wide]
\item \label{itm:piecewise infinite definable functions:bounded} If $f$ is piecewise bounded, images of $\bigwedge_A${\hyp}definable sets are $\bigwedge_A${\hyp}definable, and preimages of piecewise $\bigwedge_A${\hyp}definable sets are piecewise $\bigwedge_A${\hyp}definable.
\item \label{itm:piecewise infinite definable functions:proper} If $f$ is piecewise proper, images of piecewise $\bigwedge_A${\hyp}definable sets are piecewise $\bigwedge_A${\hyp}definable, and preimages of $\bigwedge_A${\hyp}definable sets are $\bigwedge_A${\hyp}definable.
\end{enumerate}
\begin{proof} Both are quite similar, so let us show only (\ref*{itm:piecewise infinite definable functions:bounded}). For $i\in I$ and $j\in J$ with $f[P_i]\subseteq Q_j$, consider $f_{ji}:\ P_i\rightarrow Q_j$ given by the restriction of $f$. As $f$ is piecewise $\bigwedge_A${\hyp}definable, each $f_{ji}$ is $\bigwedge_A${\hyp}definable. Given an $\bigwedge_A${\hyp}definable subset $V\subseteq P_i$, we have that $f[V]=f_{ji}[V]$, concluding that the image of $V$ is $\bigwedge_A${\hyp}definable in $Q$ by \cref{l:infinite definable functions} and the Correspondence \cref{l:correspondence lemma in pieces}. On the other hand, given a piecewise $\bigwedge_A${\hyp}definable subset $W\subseteq Q$, we have that $f^{-1}[W]\cap P_i=f^{-1}_{ji}[W\cap Q_j]$, so $f^{-1}[W]$ is piecewise $\bigwedge_A${\hyp}definable in $P$ by \cref{l:infinite definable functions}. 
\end{proof}
\end{prop}
An \emph{isomorphism of piecewise $A${\hyp}hyperdefinable sets} is a piecewise bounded and proper $\bigwedge_A${\hyp}definable bijection. In that case, we will say that $P$ and $Q$ are \emph{isomorphic over $A$}. 

\subsection{The logic topologies of piecewise hyperdefinable sets}
The \emph{$A${\hyp}logic topology} of $P$ is the respective direct limit topology. In other words, a subset of $P$ is closed if and only if it is piecewise $\bigwedge_A${\hyp}definable. By the Correspondence \cref{l:correspondence lemma in pieces}, each piece is compact and, further, every $\bigwedge_A${\hyp}definable subset is compact. 

As in the case of hyperdefinable sets,  $\overline{\{a\}}=\tp(a/A)$ for any $a\in P$. Thus, the properties $\mathrm{T}_0$ and $\mathrm{T}_1$ are equivalent. If the $A${\hyp}logic topology is $\mathrm{T}_1$, for any other small set of hyperimaginary parameters $A'$ containing $A$, the logic topologies over $A$ and over $A'$ are the same. Thus, there is at most one $\mathrm{T}_1$ logic topology on $P$ and, if it exists, it is called the \emph{global logic topology}. It follows that $P$ has a global logic topology if and only if every piece has size at most $2^{|A|+|\lang|}$, if and only if the $\mathrm{bdd}(A)${\hyp}logic topology is the global logic topology. In particular, if $P$ has a global logic topology, then $|P|\leq 2^{|A|+|\lang|}+\cf(P)$. Thus, assuming $2^{|A|+|\lang|}+\cf(P)<\kappa$, we conclude that $P$ has a global logic topology if and only if it is small.

There are still some topological properties that we want to extend from hyperdefinable sets to piecewise hyperdefinable sets. Ideally, we would like to show that these topologies are locally compact, normal and satisfy $\mathrm{T}_1\Leftrightarrow \mathrm{T}_2$. Also, we would like to show that they are closed under taking finite products. In general, these properties may fail --- see \cref{e:example 4,e:example 6,e:example no topological group} for some counterexamples. The rest of the subsection is dedicated to give sufficient natural conditions for these properties.  

\

In general, for a topological space $X$, a \emph{covering} of $X$ is a set $\mathcal{C}\subseteq \mathcal{P}(X)$ such that $\bigcup\mathcal{C}=X$, whose elements are called \emph{pieces}. A covering $\mathcal{C}$ is \emph{coherent} when, for every $U\subseteq X$, $U$ is open in $X$ if and only if $U\cap P$ is open in $P$ for each $P\in\mathcal{C}$. Equivalently, $\mathcal{C}$ is coherent when, for every $V\subseteq X$, $V$ is closed in $X$ if and only if $V\cap P$ is closed in $P$ for each $P\in\mathcal{C}$. For example, $X$ is \emph{compactly generated} if the family of all the compact subsets is a coherent covering. 

We say that a covering is \emph{local} if for every point of $X$ there is a piece that is a neighbourhood of it. The following topological results are straightforward:

\begin{lem}[Local coverings] \label{l:local coherent} 
Let $X$ be a topological space and $\mathcal{C}$ a local covering. Then, $\mathcal{C}$ is coherent.  
\end{lem} 

\begin{lem} \label{l:lemma normal coverings}
Let $X$ be a topological space and $\{P_n\}_{n\in\N}$ a closed coherent covering with $P_n\subseteq P_{n+1}$ for each $n\in\N$. Suppose that $P_n$ is normal for each $n\in\N$. Then, $X$ is normal. 
\begin{proof} Let $V$ and $W$ be closed disjoint subsets of $X$. Take the continuous map $g:\ V\sqcup W\rightarrow\{-1,1\}$ given by $g_{\mid V}=1$ and $g_{\mid W}=-1$. By recursion, using Tietze's Theorem \cite[Theorem 35.1]{munkres1999topology}, we construct a chain $(f_n)_{n\in\N}$ of continuous maps $f_n:\ P_n\rightarrow [-1,1]$, each one extending $g_{\mid P_n}$. Taking $f=\bigcup f_n$, we get a continuous map separating $V$ and $W$. 
\end{proof}
\end{lem}
In the case of a piecewise hyperdefinable set $P$, we have a directed closed and compact coherent covering $\{P_i\}_{i\in I}$. In particular, it follows that logic topologies are always compactly generated. If $P$ is countably piecewise hyperdefinable, then it is $\sigma${\hyp}compact (i.e. a countable union of compact sets). Furthermore, by \cref{l:lemma normal coverings}, we get normality:
\begin{prop} \label{p:normal countably piecewise hyperdefinable} Let $P$ be a countably piecewise $A${\hyp}hyperdefinable set. Then, $P$ is normal with the $A${\hyp}logic topology. In particular, global logic topologies of countably piecewise hyperdefinable sets are Hausdorff.
\end{prop}
Say that a piecewise hyperdefinable set $P$ is \emph{locally $A${\hyp}hyperdefinable} if its covering is local in the $A${\hyp}logic topology, i.e. if for every point of $P$ there is an $\bigwedge_A${\hyp}definable set which is a neighbourhood of it. Say that $P$ is locally hyperdefinable if it is so for some small set of parameters. 
\begin{rmk} \emph{Piecewise definable} sets are the special case of piecewise hyperdefinable sets when we only consider strict direct limits of definable sets with definable maps. Definable sets are always open in the logic topology so, following the terminology of this paper, every piecewise definable set is trivially locally (hyper)definable. The distinction between these two notions only appears in the general context of piecewise hyperdefinable sets. In particular, note that our terminology is consistent with the typical use of the terms ``piecewise definable'' and ``locally definable'' as synonyms in the literature.
\end{rmk}
\begin{prop} \label{p:mapping locally hyperdefinable sets}
Let $P$ and $Q$ be piecewise $A${\hyp}hyperdefinable sets. Assume that $P$ is locally $A${\hyp}hyperdefinable. Let $f:\ P\rightarrow Q$ be a piecewise bounded and proper $\bigwedge_A${\hyp}definable onto map. Then, $Q$ is locally $A${\hyp}hyperdefinable. 
\begin{proof} Pick $y\in Q$. By \cref{p:piecewise infinite definable functions}, $f^{-1}[\tp(y/A)]$ is $\bigwedge_A${\hyp}definable, so compact in the $A${\hyp}logic topology of $P$. As $P$ is locally hyperdefinable, we find an $\bigwedge_A${\hyp}definable set $V$ such that $f^{-1}[\tp(y/A)]\subseteq U\subseteq V$, where $P\setminus U$ is piecewise $\bigwedge_A${\hyp}definable. By \cref{p:piecewise infinite definable functions}, $f[V]$ is $\bigwedge_A${\hyp}definable, and $f[P\setminus U]$ is piecewise $\bigwedge_A${\hyp}definable. As $f^{-1}[\tp(y/A)]\subseteq U$, $\tp(y/A)\cap f[P\setminus U]=\emptyset$. Also, as $f$ is onto, $Q=f[V\cup (P\setminus U)]= f[V]\cup f[P\setminus U]$. Therefore, $y\in \tp(y/A)\subseteq Q\setminus f[P\setminus U]\subseteq f[V]$, concluding that $f[V]$ is an $\bigwedge_A${\hyp}definable neighbourhood of $y$ in the $A${\hyp}logic topology. As $y$ is arbitrary, we conclude that $Q$ is locally $A${\hyp}hyperdefinable. 
\end{proof}
\end{prop}
For locally hyperdefinable sets we have a really good control of the logic topology. 
\begin{prop}\label{p:locally compact locally hyperdefinable}
Let $P$ be a locally $A${\hyp}hyperdefinable set. Then, $P$ is locally closed compact in the $A${\hyp}logic topology, i.e. every point has a local base of closed compact neighbourhoods.
\begin{proof} Say $P=\underrightarrow{\lim}\, P_i$. Pick $x\in P$ and $U$ open neighbourhood of $x$ in the $A${\hyp}logic topology. Take a piece $P_i$ and $U_0$ such that $x\in U_0\subseteq P_i$ with $U_0$ open in the $A${\hyp}logic topology of $P$. Then, $U_1\coloneqq U\cap U_0$ is an open neighbourhood of $x$ in the $A${\hyp}logic topology of $P$. Note that $\tp(x/A)\subseteq U_1\subseteq P_i$, so $\tp(x/A)$ and $P_i\setminus U_1$ are disjoint closed subsets of $P_i$. By \cref{p:topology hyperdefinables}, $P_i$ is normal, so there are $\tp(x/A)\subseteq U'$ and $P_i\setminus U_1\subseteq P_i\setminus V'$ with $U'\cap (P_i\setminus V')=\emptyset$ such that $P_i\setminus U'$ and $V'$ are $\bigwedge_A${\hyp}definable in $P_i$. Therefore, $x\in\tp(x/A)\subseteq U'\subseteq V'\subseteq U_1\subseteq P_i$. Now, since $U_1$ is open in the $A${\hyp}logic topology of $P$ and $U'\subseteq U_1$ is open in the subspace topology of $U_1$, we conclude that $U'$ is open in the $A${\hyp}logic topology of $P$. Hence, $V'$ is an $\bigwedge_A${\hyp}definable neighbourhood of $x$ in the $A${\hyp}logic topology contained in $U$. As $x$ and $U$ are arbitrary, we conclude that $P$ is locally closed compact. 
\end{proof}
\end{prop}
\begin{prop} \label{p:compact locally hyperdefinable} 
Let $P$ be a locally $A${\hyp}hyperdefinable set. Then, every compact subset of $P$ in the $A${\hyp}logic topology is contained in the interior of some piece.
\begin{proof} Clear from the definition. 
\end{proof}
\end{prop}
\begin{prop} \label{p:hausdorff locally hyperdefinable} 
Let $P$ be a locally $A${\hyp}hyperdefinable set. Then, any two closed compact disjoint subsets in the $A${\hyp}logic topology are separated by open sets. In particular, it is $\mathrm{T}_1$ if and only if it is $\mathrm{T}_2$.
\begin{proof} Say $P=\underrightarrow{\lim}\, P_i$. Let $K_1$ and $K_2$ be two disjoint closed compact subsets of $P$ in the $A${\hyp}logic topology. By \cref{p:compact locally hyperdefinable}, there is a piece $P_i$ and an open subset $U$ of $P$ such that $K_1,K_2\subseteq U\subseteq P_i$. By normality of $P_i$, there are disjoint open subsets $U_1$ and $U_2$ in the $A${\hyp}logic topology of $P_i$ such that $K_1\subseteq U_1$ and $K_2\subseteq U_2$. Thus, $U\cap U_1$ and $U\cap U_2$ are disjoint open subsets in the subspace topology of $U$ separating $K_1$ and $K_2$. As $U$ is open in the $A${\hyp}logic topology of $P$, we conclude that $U\cap U_1$ and $U\cap U_2$ are disjoint open subsets in the $A${\hyp}logic topology of $P$ separating $K_1$ and $K_2$. 
\end{proof}
\end{prop}
Recall that a function between topological spaces is \emph{proper} if the preimage of every compact set is compact. For instance, every closed function with compact fibres is proper \cite[Theorem 3.7.2]{engelking1989general}. 
\begin{prop} \label{p:functions logic topology} 
Let $P$ and $Q$ be piecewise $A${\hyp}hyperdefinable sets and $f:\ P\rightarrow Q$ a piecewise $\bigwedge_A${\hyp}definable function.
\begin{enumerate}[label={\rm{(\arabic*)}}, ref={\arabic*}, wide]
\item \label{itm:functions logic topology:continuous} If $f$ is piecewise bounded, then $f$ is continuous between the $A${\hyp}logic topologies.
\item \label{itm:functions logic topology:proper} If $f$ is piecewise proper, then $f$ is closed and has compact fibres between the $A${\hyp}logic topologies. In particular, it is proper.
\item \label{itm:functions logic topology:homeomorphism} If $f$ is an isomorphism of piecewise $A${\hyp}hyperdefinable sets, then $f$ is a homeomorphism between the $A${\hyp}logic topologies.
\item \label{itm:functions logic topology:continuous iff piece bounded} If $Q$ is locally $A${\hyp}hyperdefinable, then $f$ is continuous between the $A${\hyp}logic topologies if and only if it is piecewise bounded.
\item \label{itm:functions logic topology:proper iff piece proper} If $P$ is locally $A${\hyp}hyperdefinable, then $f$ is proper between the $A${\hyp}logic topologies if and only if it is piecewise proper.
\item \label{itm:functions logic topology:homeomorphism iff piece isomorphism} If $P$ and $Q$ are locally $A${\hyp}hyperdefinable, then $f$ is an isomorphism of piecewise $A${\hyp}hyperdefinable sets if and only if it is a homeomorphism between the $A${\hyp}logic topologies.
\end{enumerate}
\begin{proof} Point (\ref*{itm:functions logic topology:continuous}) is given by \cref{p:piecewise infinite definable functions}(\ref{itm:piecewise infinite definable functions:bounded}). 

Point (\ref*{itm:functions logic topology:proper}): closedness is given by \cref{p:piecewise infinite definable functions}(\ref{itm:piecewise infinite definable functions:proper}). On the other hand, for any point $a\in Q$, $f^{-1}(a)\subseteq f^{-1}[\tp(a/A)]$ and $f^{-1}[\tp(a/A)]$ is $\bigwedge_A${\hyp}definable, so compact. As $f$ is $A${\hyp}invariant and $\tp(a/A)$ is the orbit of $a$ under $\Aut(\mathfrak{M}/A)$ by \cref{c:types over hyperimaginaries}, it follows that $f^{-1}[\tp(a/A)]$ is the orbit of $f^{-1}(a)$ under $\Aut(\mathfrak{M}/A)$, i.e. the smallest $A${\hyp}invariant set containing $f^{-1}(a)$. Consequently, any open covering of $f^{-1}(a)$ in the $A${\hyp}logic topology is also a covering of $f^{-1}[\tp(a/A)]$, so $f^{-1}(a)$ is compact too. 

Point (\ref*{itm:functions logic topology:homeomorphism}) is given by points (\ref*{itm:functions logic topology:continuous}) and (\ref*{itm:functions logic topology:proper}). Point (\ref*{itm:functions logic topology:continuous iff piece bounded}) is given by (\ref*{itm:functions logic topology:continuous}) and \cref{p:compact locally hyperdefinable}. Point (\ref*{itm:functions logic topology:proper iff piece proper}) is given by (\ref*{itm:functions logic topology:proper}) and \cref{p:compact locally hyperdefinable}. Point (\ref*{itm:functions logic topology:homeomorphism iff piece isomorphism}) is given by points (\ref*{itm:functions logic topology:continuous iff piece bounded}) and (\ref*{itm:functions logic topology:proper iff piece proper}).
\end{proof}
\end{prop}
\begin{prop} \label{p:functions logic topology 2}
Let $P$ be a piecewise $A${\hyp}hyperdefinable set and  $Q$ a locally $A${\hyp}hyperdefinable set whose $A${\hyp}logic topology is its global logic topology. Let $f:\ P\rightarrow Q$ be a function. Then, $f$ is continuous between the $A${\hyp}logic topologies if and only if $f$ is a piecewise bounded $\bigwedge_A${\hyp}definable function.
\begin{proof} By the Closed Graph Theorem \cite[Exercise 8, Section 26]{munkres1999topology} and \cref{p:functions logic topology,p:compact locally hyperdefinable}. 
\end{proof}
\end{prop}
By \cref{p:functions logic topology}(\ref{itm:functions logic topology:continuous}), Cartesian projection maps are continuous between the logic topologies. Therefore, the $A${\hyp}logic topology of a finite product of piecewise $A${\hyp}hyperdefinable sets is at least as fine as the product topology of the $A${\hyp}logic topologies. In the case of local hyperdefinable sets with global logic topologies, we have that they coincide.
\begin{prop} \label{p:product of locally hyperdefinable sets}
Let $P$ and $Q$ be locally hyperdefinable sets with global logic topologies. Then, $P\times Q$ is a locally hyperdefinable set and its product topology is its global logic topology.
\begin{proof} It suffices to note that, for any two topological spaces $X$ and $Y$ with respective local coverings $\mathcal{C}$ and $\mathcal{D}$, $\mathcal{C}\times \mathcal{D}\coloneqq \{P\times Q\sth P\in\mathcal{C},\ Q\in\mathcal{D}\}$ is a local covering of the product topology. 
\end{proof}
\end{prop}
Similarly, in the case of countably piecewise hyperdefinable sets with global logic topologies, we also conclude that they coincide.
\begin{prop} \label{p:product of countably piecewise hyperdefinable sets} Let $P$ and $Q$ be two countably piecewise hyperdefinable sets with global logic topologies. Then, $P\times Q$ is a countably piecewise hyperdefinable set and its product topology is its global logic topology.
\begin{proof} Say $P=\underrightarrow{\lim}\, P_n$ and $Q=\underrightarrow{\lim}\, Q_n$. Let $\Gamma\subseteq P\times Q$ be closed in the global logic topology and $(a,b)\notin \Gamma$ arbitrary. We recursively define a sequence $(U_n,V_n)_{n\in\N}$ such that $a\in U_0$, $b\in V_0$, $U_n\subseteq U_{n+1}$, $V_n\subseteq V_{n+1}$, $U_n\subseteq P_n$ is open in $P_n$, $V_n\subseteq Q_n$ is open in $Q_n$ and $(\overline{U}_n\times \overline{V}_n)\cap \Gamma=\emptyset$. 

For any $n\in\N$, note that the product topology and the global logic topology in $P_n\times Q_n$ coincide. As $(a,b)$ and $\Gamma\cap (P_0\times Q_0)$ are disjoint and closed, by normality in $P_0\times Q_0$, we can find $U_0\subseteq P_0$ open in $P_0$ and $V_0\subseteq Q_0$ open in $Q_0$ such that $(a,b)\in U_0\times V_0$ and $(\overline{U}_0\times\overline{V}_0)\cap \Gamma=\emptyset$. Now, suppose $U_n$ and $V_n$ are defined. By normality of $P_{n+1}\times Q_{n+1}$, for any $(x,y)\in \overline{U}_n\times \overline{V}_n$, we can find $U^{xy}_{n+1}\subseteq P_{n+1}$ open in $P_{n+1}$ and $V^{xy}_{n+1}\subseteq Q_{n+1}$ open in $Q_{n+1}$ such that $x\in U^{xy}_{n+1}$, $y\in V^{xy}_{n+1}$ and $(\overline{U}^{xy}_{n+1}\times \overline{V}^{xy}_{n+1})\cap \Gamma=\emptyset$. As $\overline{V}_n$ is compact, for each $x\in \overline{U}_n$ there is $F_x\subseteq \overline{V}_n$ finite such that $\overline{V}_n\subseteq \bigcup_{y\in F_x}V^{xy}_{n+1}$. Take $V^x_{n+1}=\bigcup_{y\in F_x}V^{xy}_{n+1}$ and $U^x_{n+1}=\bigcap_{y\in F_x}U^{xy}_n$. Then, $U^x_{n+1}\subseteq P_{n+1}$ is open in $P_{n+1}$, $V^x_{n+1}\subseteq Q_{n+1}$ is open in $Q_{n+1}$, $x\in U^x_{n+1}$, $\overline{V}_n\subseteq V^x_{n+1}$ and $(\overline{U}^x_{n+1}\times \overline{V}^x_{n+1})\cap \Gamma=\emptyset$. As $\overline{U}_n$ is compact, there is $F\subseteq \overline{U}_n$ finite such that $\overline{U}_n\subseteq \bigcup_{x\in F}U^x_{n+1}$. Take $U_{n+1}=\bigcup_{x\in F}U^x_{n+1}$ and $V_{n+1}=\bigcap_{x\in F} V^x_{n+1}$. Then, $U_{n+1}\subseteq P_{n+1}$ is open in $P_{n+1}$, $V_{n+1}\subseteq Q_{n+1}$ is open in $Q_{n+1}$, $\overline{U}_n\subseteq U_{n+1}$, $\overline{V}_n\subseteq V_{n+1}$ and $(\overline{U}_{n+1}\times \overline{V}_{n+1})\cap\Gamma=\emptyset$.

Let $U=\bigcup_{n\in\N} U_n$ and $V=\bigcup_{n\in\N} V_n$. For $n\in\N$, we have $U\cap P_n=\bigcup_{m\geq n} (U_m\cap P_n)$ and $V\cap Q_n=\bigcup_{m\geq n}(V_m\cap Q_n)$. As $P_n\subseteq P_m$ and $Q_n\subseteq Q_m$ for $m>n$, we get that $U_m\cap P_n$ is open in $P_n$ and $V_m\cap Q_n$ is open in $Q_n$, so $U\cap P_n$ is open in $P_n$ and $V\cap Q_n$ is open in $Q_n$ for each $n\in\N$. Therefore, $U\times V$ is open in the product topology and $(a,b)\in U\times V$ with $(U\times V)\cap \Gamma=\emptyset$. As $(a,b)\notin \Gamma$ is arbitrary, $\Gamma$ is closed in the product topology. As $\Gamma$ is arbitrary, we conclude that the global logic topology in $P\times Q$ is the same as the product topology. 
\end{proof}
\end{prop}
Let $P=\underrightarrow{\lim}\, P_i$ be piecewise $A${\hyp}hyperdefinable set and $V\subseteq P$ a non{\hyp}empty piecewise $\bigwedge_A${\hyp}definable subset. Note that the subspace topology of $V$ inherited from the $A${\hyp}logic topology of $P$ is the $A${\hyp}logic topology of $V$ given as the piecewise $A${\hyp}hyperdefinable set $\underrightarrow{\lim}\ V\cap P_i$.  We conclude showing that local hyperdefinability is hereditary.
\begin{prop} \label{p:subspace locally hyperdefinable} 
Let $P=\underrightarrow{\lim}\, P_i$ be a locally $A${\hyp}hyperdefinable set and $V\subseteq P$ a piecewise $\bigwedge_A${\hyp}definable subset. Then, $V$ with the induced piecewise hyperdefinable substructure $\underrightarrow{\lim}\ V\cap P_i$ is a locally $A${\hyp}hyperdefinable set.
\begin{proof} It suffices to note that, for any topological space $X$ with local covering $\mathcal{C}$ and any $Y\subseteq X$, $\mathcal{C}_{\mid Y}=\{C\cap Y\sth C\in\mathcal{C}\}$ is a local covering of the subspace topology. 
\end{proof}
\end{prop}
\subsection{Spaces of types}
Let $P=\rfrac{X}{E}$ be an $A${\hyp}hyperdefinable set and $F\subseteq P\times P$ an $\bigwedge_A${\hyp}definable equivalence relation in $P$. Then, $\rfrac{P}{F}$ is actually identified with the $A${\hyp}hyperdefinable set $\rfrac{X}{\quot^{-1}_{P\times P}[F]}$ via the canonical bijection 
\[\eta:\ [x]_{\quot^{-1}_{P\times P}[F]}\mapsto [[x]_P]_F.\]
Furthermore, by the definition of the quotient topologies, this map is a homeomorphism.  

Let $P=\underrightarrow{\lim}_{(I,\prec)}(P_i,\varphi_{ji})_{j\succeq i}$ be a piecewise $A${\hyp}hyperdefinable set. Write $P_i=\rfrac{X_i}{E_i}$. Let $F\subseteq P\times P$ be a piecewise $\bigwedge_A${\hyp}definable equivalence relation. Write $F_i\coloneqq F_{\mid P_i\times P_i}$. Define $\widetilde{\varphi}_{ji}:\ \rfrac{P_i}{F_i}\rightarrow \rfrac{P_j}{F_j}$ by $\widetilde{\varphi}_{ji}([x]_{F_i})=[\varphi_{ji}(x)]_{F_j}$. Clearly, they are well{\hyp}defined $\bigwedge_A${\hyp}definable and $1${\hyp}to{\hyp}$1$ functions and $\widetilde{\varphi}_{ii}=\mathrm{id}$ and $\widetilde{\varphi}_{ki}=\widetilde{\varphi}_{kj}\circ\widetilde{\varphi}_{ji}$. Then, we have canonically the piecewise $A${\hyp}hyperdefinable structure in $\rfrac{P}{F}$ given by
\[\rfrac{P}{F}\coloneqq \underrightarrow{\lim}_{(I,\prec)}(\rfrac{P_i}{F_i},\widetilde{\varphi}_{ji})_{j\succeq i},\]
via the map $\eta:\ [\inc_{P_i}(x)]_F\mapsto \inc_{P_i/F_i}([x]_{F_i})$ for $i$ such that $x\in P_i$. Furthermore, topologically, by the definitions of the quotient and the direct limit topologies, this bijection is clearly a homeomorphism. \medskip

Let $P=\rfrac{X}{E}$ be an $A^*${\hyp}hyperdefinable set. Consider the space of types $\Stone_X(A^*)=\{\tp(x/A^*)\sth x\in X\}$ with the usual topology. We define the equivalence relation $\sim_E$ in $\Stone_X(A^*)$ as $p\sim_Eq$ if and only if there are $E${\hyp}equivalent realisations of $p(x)$ and $q(y)$. The \emph{space of types} of $P$ over $A^*$ is the space $\Stone_P(A^*)\coloneqq \rfrac{\Stone_X(A^*)}{\sim_E}$ with the quotient topology. 

On the other hand, by \cref{l:equivalence relation of equality of types}, $\rfrac{P}{\Delta_P(A^*)}=\{\tp(a/A^*)\sth a\in P\}$ is an $A^*${\hyp}hyperdefinable set and has its $A^*${\hyp}logic topology. 
\begin{prop} \label{p:space of types} 
Let $P=\rfrac{X}{E}$ be an $A^*${\hyp}hyperdefinable set. Then, $\rfrac{P}{\Delta_P(A^*)}$ and $\Stone_P(A^*)$ are homeomorphic.
\begin{proof}
Consider the map
\[\begin{array}{lccc} f:& \Stone_{P}(A^*)&\rightarrow & \rfrac{P}{\Delta_P(A^*)}\\ & \left[\tp(a^*/A^*)\right]_{\sim_{E}}&\mapsto& \tp([a^*]_{E}/A^*).\end{array}\]
By saturation, it is well{\hyp}defined. Clearly, it is onto. It is $1${\hyp}to{\hyp}$1$ by \cref{l:types of hyperimaginaries}. By the definition of the quotient topology, it is clear that $f$ is continuous. As $\Stone_{X}(A^*)$ is a compact topological space, $\Stone_{P}(A^*)$ is also compact. On the other hand, $\rfrac{P}{\Delta_P(A^*)}$ is a Hausdorff space. Then, as the domain is compact, the image is Hausdorff and $f$ is a continuous bijection, we conclude that $f$ is a homeomorphism. 
\end{proof}
\end{prop}
Therefore, for hyperimaginary parameters, we define the space of types of $P$ over $A$, $\Stone_P(A)$, as $\rfrac{P}{\Delta_P(A)}$ with its global logic topology.
\begin{rmk} By \cref{l:infinite definable functions,l:types and infinite definable functions}, any $\bigwedge_A${\hyp}definable function $f:\ P\rightarrow Q$ induces a continuous closed map between the spaces of types.
\end{rmk}
Similarly, let $P=\underrightarrow{\lim}_{I}P_i$ be a piecewise $A^*${\hyp}hyperdefinable set. For each $i\in I$, we have the space of types $\Stone_{P_i}(A^*)=\{\tp(a/A^*)\sth a\in P_i\}$. Given $i,j\in I$ with $i\preceq j$, the map $\varphi_{ji}:\ P_i\rightarrow P_j$ induces a continuous and closed $1${\hyp}to{\hyp}$1$ function $\Ima\ \varphi_{ji}:\ \Stone_{P_i}(A^*)\rightarrow \Stone_{P_j}(A^*)$ given by $\tp(a/A^*)\mapsto \tp(\varphi_{ji}(a)/A^*)$. Clearly, $\Ima\ \varphi_{kj}\circ\Ima\ \varphi_{ji}=\Ima\ \varphi_{ki}$ and $\Ima\ \varphi_{ii}=\mathrm{id}$. Therefore, we have a topological direct sequence $(\Stone_{P_i}(A^*),\Ima\ \varphi_{ji})_{j\succeq i}$. The \emph{space of types} of $P$ over $A^*$ is then the direct limit topological space 
\[\Stone_{P}(A^*)\coloneqq \underrightarrow{\lim}_I\Stone_{P_i}(A^*).\]
By \cref{p:space of types}, we have 
\[\rfrac{P}{\Delta_{P}(A^*)}=\underrightarrow{\lim}\,\rfrac{P_i}{\Delta_{P_i}(A^*)}\cong\underrightarrow{\lim}\,\Stone_{P_i}(A^*)=\Stone_{P}(A^*).\]

Therefore, for a piecewise hyperdefinable set $P$, like in the case of hyperdefinable sets, we also define the \emph{space of types} of $P$ over hyperimaginary parameters $A$, $\Stone_P(A)$, as $\rfrac{P}{\Delta_{P}(A)}$ with its global logic topology. Even more, we say that a piecewise hyperdefinable set $S$ is a \emph{space of types} if it has a global logic topology.\medskip

Let $X$ be a topological space. Recall that two points are \emph{topologically indistinguishable} if they have the same neighbourhoods. A \emph{Kolmogorov map} is an onto continuous and closed map $k:\ X\rightarrow Y$ between topological spaces such that $k(a)$ and $k(b)$ are topologically indistinguishable if and only if $a$ and $b$ are topologically indistinguishable. We will use the following basic characterisation --- see \cite{pirttimaki2021survey} for more details: 
\begin{lem} \label{l:kolmogorov maps} 
Let $X$ and $Y$ be topological spaces and $k:\ X\rightarrow Y$ a function. Then, the following are equivalent:\smallskip
\begin{enumerate}[label={\rm{(\alph*)}}, wide]
\item $k$ is a Kolmogorov map.
\item $\Ima\, k$ is a lattice isomorphism between the topologies with inverse $\Ima^{-1}k$.
\end{enumerate}
\end{lem}
\begin{rmk} When $Y$ is $\mathrm{T}_0$, we call a Kolmogorov map $k:\ X\rightarrow Y$ the \emph{Kolmogorov quotient} of $X$. In that case, up to a homeomorphism, $Y$ is the quotient space $\rfrac{X}{\sim}$, where $\sim$ is the topologically indistinguishable equivalence relation, and $k$ is the respective quotient map.
\end{rmk}
As it is discussed in the introduction, it is useful to note that the quotient map $\tp_A:\ P\rightarrow \Stone_{P}(A)$ given by $\tp_A:\ a\mapsto \tp(a/A)=\rfrac{a}{\Delta_P(A)}$ is the Kolmogorov quotient between the $A${\hyp}logic topologies. By definition, it also satisfies that $\tp^{-1}_A[\tp_A[V]]=V$ for any $A${\hyp}invariant set $V$. Thus, when studying $P$, it is typically possible to map the discussion via $\tp_A$, argue in $\Stone_P(A)$ and then lift the conclusions to $P$ via $\tp^{-1}_A$. Following this procedure, one can usually assume without loss of generality that $P$ is $\mathrm{T}_0$. This technique will be illustrated in the following subsection --- see \cref{t:uniform structure} and the Metrisation \cref{t:metrisation theorem}.
\begin{rmk} By \cref{p:functions logic topology}(\ref{itm:functions logic topology:continuous}) and \cref{l:types and piecewise infinite definable functions}, any piecewise bounded $\bigwedge_A${\hyp}definable function induces a continuous map between the spaces of $A${\hyp}types. 
\end{rmk}
\subsection{Metrisation results}\hspace{-.5em}\footnote{This section is inspired on results from \cite{ben2005uncountable}.}
A \emph{uniform space} is a pair $(X,\Phi)$ formed by a non{\hyp}empty set $X$ and a filter $\Phi$ of binary relations in $X$ satisfying that, for every $U\in\Phi$, 
\[\begin{array}{ll} 
\mathbf{(i)}& \Delta\subseteq U,\\
\mathbf{(ii)}& U^{-1}\in \Phi \mathrm{\ and\ }\\
\mathbf{(iii)}& \mathrm{there\ is\ }V\in \Phi\mathrm{\ such\ that\ }V\circ V\subseteq U;
\end{array}\]
where $\Delta\coloneqq\{(x,x)\sth x\in X\}$ is the diagonal (equality) relation, $U^{-1}\coloneqq\{(y,x)\sth (x,y)\in U\}$ and $W\circ V\coloneqq\{(x,z)\sth \exists y\mathrel{} (x,y)\in V,\ (y,z)\in W\}$. The filter $\Phi$ is called the \emph{uniform structure} of the uniform space $(X,\Phi)$.

A \emph{uniformity base} of $(X,\Phi)$ is a filter base of $\Phi$. Note that a filter base $\mathcal{B}$ of reflexive binary relations on $X$ is a uniformity base of some uniform structure if and only if, for any $U\in\mathcal{B}$, there are $V,W\in \mathcal{B}$ such that $V\circ V\subseteq U$ and $W\subseteq U^{-1}$.

Uniform spaces generalise (pseudo{\hyp})metric spaces. In other words, every pseudo{\hyp}metric $\rho$ induces a uniform space given by the uniformity base $\{\rho^{-1}[0,\varepsilon)\sth \varepsilon\in\R_{>0}\}$. We say then that a uniform structure $\Phi$ is \emph{pseudo{\hyp}metrisable} if it arises in this way; that means there is a pseudo{\hyp}metric $\rho$ such that $\{\rho^{-1}[0,\varepsilon)\sth \varepsilon\in\R_{>0}\}$ is a uniformity base of $\Phi$. 

As in the case of metric spaces, every uniform space has a topology given by the system of local bases of neighbourhoods $\{U(a)\sth U\in \Phi\}_{a\in X}$, where $U(a)\coloneqq\{b\sth (a,b)\in U\}$. We say that a topology $\mathcal{T}$ \emph{admits a uniform structure} if there is a uniform structure $\Phi$ on $X$ with uniform topology $\mathcal{T}$. Note that a topological space could admit many different uniform structures. Uniform structures inducing the same topology are called \emph{equivalent}.
\begin{rmk} Uniform structures are the natural abstract context to study uniform continuity, uniform convergence, Cauchy sequences or completeness. It is important to note that two equivalent uniform structures may differ on these aspects. For example, a uniformly continuous function with respect to one uniform structure might not be uniformly continuous in another equivalent uniform structure.
\end{rmk} 
Recall that a topological space $X$ is \emph{functionally regular} if, for any point $x\in X$ and any closed set $V\subseteq X$ such that $x\notin V$, there is a continuous function $f:\ X\rightarrow [0,1]$ such that $f(x)=0$ and $f_{\mid V}=1$. It is a well{\hyp}known fact from the theory of uniform spaces that a topological space admits a uniform structure if and only if it is functionally regular --- see \cite[Theorem 38.2]{willard1970general} for a proof. Hence, we get the following general result for countably piecewise $A${\hyp}hyperdefinable sets.
\begin{prop} \label{p:uniformity of countably piecewise hyperdefinable sets} 
Every countably piecewise $A${\hyp}hyperdefinable set $P$ with its $A${\hyp}logic topology admits a uniform structure.
\begin{proof} We know that $P$ is normal by \cref{p:normal countably piecewise hyperdefinable}. Let $V$ be piecewise $\bigwedge_A${\hyp}definable and $a\in P\setminus V$. As $\overline{\{a\}}=\tp(a/A)$, we have that $\overline{\{a\}}$ and $V$ are disjoint. By normality, using Urysohn's Lemma \cite[Theorem 33.1]{munkres1999topology}, we conclude that $P$ is functionally regular. Thus, we conclude that it admits a uniform structure. 
\end{proof}
\end{prop}
Our aim now is to give a better description of the uniform structure in each piece. In other words, we want to give an actual uniformity base for the logic topology.

Let $P$ be an $A${\hyp}hyperdefinable set. We say that an $\bigwedge_A${\hyp}definable binary relation $\varepsilon\subseteq P\times P$ is \emph{positive} if $\Delta_P(A)\subseteq \mathring{\varepsilon}$, where the interior is taken in the product topology of the $A${\hyp}logic topology. Write $\mathcal{E}_P(A)$ for the set of all positive $\bigwedge_A${\hyp}definable binary relations of $P$. It could seem odd the fact that we take the interior in the product topology. The following lemma explains why.

\begin{lem} Let $P$ be an $A${\hyp}hyperdefinable set and $\pi:\ P\times P\rightarrow \Stone_P(A)\times\Stone_P(A)$ the quotient map $(a,b)\mapsto (\tp(a/A),\tp(b/A))$. Then, for any $\varepsilon\in\mathcal{E}_P(A)$, $\pi[\varepsilon]\in\mathcal{E}_{\Stone_P(A)}(A)$, and for any $\varepsilon'\in \mathcal{E}_{\Stone_P(A)}(A)$, $\pi^{-1}[\varepsilon']\in \mathcal{E}_P(A)$. 
\end{lem}
\begin{rmk} When the $A${\hyp}logic topology is the global topology, $\Delta_P(A)$ is precisely the diagonal $\Delta$ and $\mathcal{E}_P(A)$ is the set of closed neighbourhoods of $\Delta$ in the global logic topology of $P\times P$, which coincides with the product topology by \cref{p:product logic topology}. This in particular applies to $\Stone_P(A)$.
\end{rmk}
We now prove the main result of this subsection:
\begin{theo} \label{t:uniform structure} 
Let $P$ be an $A${\hyp}hyperdefinable set. Then, $P$ with the $A${\hyp}logic topology admits a unique uniform structure and $\mathcal{E}_P(A)$ is a uniformity base of it.
\begin{proof} We start by proving the proposition for $\Stone_P(A)$ rather than $P$, so consider the set $\mathcal{E}_{\Stone_P(A)}(A)$ of positive $\bigwedge_A${\hyp}definable binary relations on $\Stone_P(A)$. By \cite[Theorem 36.19]{willard1970general}, since $\Stone_P(A)$ with the global logic topology is compact and Hausdorff, it admits one and only one uniform structure, which is precisely the filter of neighbourhoods of the diagonal in the product topology. Since the global logic topology and the product topology on $\Stone_P(A)\times\Stone_P(A)$ coincide by \cref{p:product logic topology}, $\mathcal{E}_{\Stone_P(A)}(A)$ is the collection of all closed neighbourhoods of the diagonal. Thus, by normality of $\Stone_P(A)\times\Stone_P(A)$ and closedness of $\Delta$ (which follows from Hausdorffness of $\Stone_P(A)$), the filter generated by $\mathcal{E}_{\Stone_P(A)}(A)$ is exactly the collection of all neighbourhoods of $\Delta$, which is precisely the unique uniform structure admitted by $\Stone_P(A)$. Hence, $\mathcal{E}_{\Stone_P(A)}(A)$ is a base of the unique uniform structure of $\Stone_P(A)$.

We now use the quotient maps $\tp_A:\ P\rightarrow \Stone_P(A)$ and $\tau:\ P\times P\rightarrow\Stone_P(A)\times\Stone_P(A)$ given by $\tp_A(a)=\tp(a/A)$ and $\tau(a,b)=(\tp(a/A),\tp(b/A))$ to extend the result to $P$. 
Let $\Phi$ be any uniform structure on $P$ inducing the $A${\hyp}logic topology. Note that at least one $\Phi$ exists by \cref{p:uniformity of countably piecewise hyperdefinable sets}. Consider $\tau\Phi\coloneqq \{\tau[U]\sth U\in \Phi\}$. Obviously, $\tau\Phi$ is a filter. Also, $\tau[U^{-1}]=\tau[U]^{-1}$ and $\tau[U]\circ\tau[U]\subseteq \tau[U\circ U\circ U]$ for any $U\in\Phi$, so $\tau\Phi$ is a uniform structure on $\Stone_P(A)$. On the other hand, for any $a\in P$, $\tp_A[U(a)]\subseteq\tau[U](\tp(a/A))\subseteq \tp_A[U\circ U(a)]$. Thus, by the properties of $\tp_A$, $\tau\Phi$ induces the global logic topology on $\Stone_P(A)$. As $\Stone_P(A)$ with the global logic topology only admits one uniform structure, it follows that $\mathcal{E}_{\Stone_P(A)}(A)$ is a uniformity base of $\tau\Phi$. In particular, we have $\tau^{-1}\mathcal{E}_{\Stone_P(A)}(A)\subseteq \Phi$, where $\tau^{-1}\mathcal{E}_{\Stone_P(A)}(A)\coloneqq \{\tau^{-1}[\varepsilon]\sth \varepsilon\in\mathcal{E}_{\Stone_P(A)}(A)\}$. On the other hand, for $U\in\Phi$, take $V\in\Phi$ such that $V\circ V\circ V\subseteq U$ and find $\varepsilon\in\mathcal{E}_{\Stone_P(A)}(A)$ such that $\varepsilon\subseteq \tau[V]$. Therefore, $\tau^{-1}[\varepsilon]\subseteq \tau^{-1}[\tau[V]]=\Delta_P(A)\circ V\circ\Delta_P(A)\subseteq V\circ V\circ V\subseteq U$. Hence, $\tau^{-1}\mathcal{E}_{\Stone_P(A)}(A)$ is a uniformity base of $\Phi$, concluding uniqueness. 

We now show that $\mathcal{E}_P(A)$ is a uniformity base of $\Phi$. Since $\tau^{-1}[\varepsilon]\in\mathcal{E}_P(A)$ for any $\varepsilon\in \mathcal{E}_{\Stone_P(A)}(A)$, it is enough to show that for any $\varepsilon\in \mathcal{E}_P(A)$ there is $\varepsilon'\in\mathcal{E}_{\Stone_P(A)}(A)$ such that $\tau^{-1}[\varepsilon']\subseteq \varepsilon$. Take $\varepsilon''\in \mathcal{E}_P(A)$ such that $\varepsilon''\circ\varepsilon''\circ\varepsilon''\subseteq \varepsilon$ and set $\varepsilon'=\tau[\varepsilon'']$. Then, $\tau^{-1}[\varepsilon']=\Delta_P(A)\circ\varepsilon''\circ\Delta_P(A)\subseteq \varepsilon$. 
\end{proof}
\end{theo}
Using compactness, we can get an even smaller uniformity base:
\begin{lem} \label{l:size uniformity base} Let $P$ be an $A${\hyp}hyperdefinable set. There is a family $\{\varepsilon_i\}_{i<|A|+|\lang|}$ of $\bigwedge_A${\hyp}definable positive binary relations on $P$ which is a uniformity base of $P$ with the $A${\hyp}logic topology.
\begin{proof} Say $\underline{\Delta}_P(A)=\bigwedge_{i< |A|+|\lang|} \varphi_i$ with $\varphi_i\in\For(\lang(A^*))$. For each $i<|A|+|\lang|$, write $E_i\coloneqq \quot_{P\times P}[\varphi_i(\mathfrak{M})]$ and $U_i\coloneqq P\times P\setminus \quot_{P\times P}[\neg\varphi_i(\mathfrak{M})]$. Since $\quot_{P\times P}^{-1}[\Delta_P(A)]\cap \neg\varphi_i(\mathfrak{M})=\emptyset$, $\Delta_P(A)\cap \quot_{P\times P}[\neg\varphi_i(\mathfrak{M})]=\emptyset$, so $\Delta_P(A)\subseteq U_i\subseteq E_i$ where $E_i$ is closed and $U_i$ is open in the $A${\hyp}logic topology. By compactness of $P\times P$ in the $A${\hyp}logic topology, we can take $\varepsilon_i\in \mathcal{E}_P(A)$ such that $\varepsilon_i\subseteq U_i\subseteq E_i$. It follows that $\{\varepsilon_i\}_{i<|A|+|\lang|}$ is a sequence of $\bigwedge_A${\hyp}definable positive binary relations on $P$ such that $\Delta_P(A)=\bigcap_{i<|A|+|\lang|}\varepsilon_i$. Now we show that $\{\varepsilon_i\}_{i<|A|+|\lang|}$ is a uniformity base. Take $\varepsilon\in \mathcal{E}_P(A)$ and $U\subseteq \varepsilon$ open in the product topology of $P\times P$ such that $\bigcap_{i<|A|+|\lang|} \varepsilon_i=\Delta_P(A)\subseteq U\subseteq \varepsilon$. Since $U$ is open in the product topology, it is also open in the $A${\hyp}logic topology. Then, by compactness, there is $i<|A|+|\lang|$ such that $\varepsilon_i\subseteq U\subseteq \varepsilon$, concluding that $\{\varepsilon_i\}_{i<|A|+|\lang|}$ is a uniformity base of $P$ with the $A${\hyp}logic topology.
\end{proof}
\end{lem}
Recall that a uniform structure is pseudo{\hyp}metrisable if and only if it has a countable uniformity base --- see \cite[Theorem 38.3]{willard1970general}. Therefore, we get the following metrisation results:
\begin{coro} Let $P$ be an $A${\hyp}hyperdefinable set. Then, $P$ with its $A${\hyp}logic topology is pseudo{\hyp}metrisable if and only if there is a countable family $\{\varepsilon_n\}_{n\in\N}$ of $\bigwedge_A${\hyp}definable positive binary relations on $P$ such that $\Delta_P(A)=\bigcap_{n\in \N} \varepsilon_n$. In particular, $P$ with the $A${\hyp}logic topology is pseudo{\hyp}metrisable if $\lang$ and $A$ are countable.
\end{coro}
\begin{theo}[Metrisation Theorem] \label{t:metrisation theorem} Let $P$ be a locally $A${\hyp}hyperdefinable set of countable cofinality. Then, $P$ with the $A${\hyp}logic topology is pseudo{\hyp}metrisable if and only if each piece is pseudo{\hyp}metrisable. In particular, $P$ with the $A${\hyp}logic topology is pseudo{\hyp}metrisable if $\lang$ and $A$ are countable.  
\begin{proof} Assume that each piece is pseudo{\hyp}metrisable. Taking the quotient map $a\mapsto \tp(a/A)$, we get that $\Stone_P(A)$ is locally metrisable and $\sigma${\hyp}compact. By $\sigma${\hyp}compactness, $\Stone_P(A)$ is trivially Lindel\"of (i.e. every open cover has a countable subcover). By \cref{p:normal countably piecewise hyperdefinable}, $\Stone_P(A)$ is normal and Hausdorff, so it is in particular regular (i.e. any closed set and any point outside it can be separated by open sets). Therefore, by \cite[Theorem 41.5]{munkres1999topology}, $\Stone_P(A)$ is paracompact (i.e. every open cover has a locally finite open refinement). By Smirnov's Metrisation Theorem \cite[Theorem 42.1]{munkres1999topology}, we conclude that it is metrisable. Taking the composition with the quotient map $a\mapsto \tp(a/A)$, we conclude that $P$ with the $A${\hyp}logic topology is pseudo{\hyp}metrisable.  
\end{proof}
\end{theo}
Using uniformities, it is also easy to find small dense subsets: 
\begin{coro} Let $P$ be a piecewise $A${\hyp}hyperdefinable set. Then, there is a subset $D\subseteq P$ with $|D|\leq |A|+|\lang|+\cf(P)$ which is dense in the $A${\hyp}logic topology. In particular, when $A$ and $\lang$ are countable, every countably piecewise $A${\hyp}hyperdefinable set is separable (i.e. has a countable dense subset) with the $A${\hyp}logic topology.
\begin{proof} Say $P=\underrightarrow{\lim}_{i<\cf(P)}\, P_i$. By \cref{l:size uniformity base}, for each $i$, there is a uniformity base $\mathcal{B}_i\subseteq \mathcal{E}_{P_i}(A)$ with $|\mathcal{B}_i|\leq |A|+|\lang|$. By compactness, for each $\varepsilon\in \mathcal{B}_i$, there is $D_{\varepsilon,i}\subseteq P_i$ finite such that $P_i\subseteq\bigcup_{a\in D_{\varepsilon,i}} \varepsilon(a)$. Take $D=\bigcup_{i<\cf(P)} \bigcup_{\varepsilon\in\mathcal{B}_i} D_{\varepsilon,i}$, so $|D|\leq |A|+|\lang|+\cf(P)$. Let $U$ be open. Then, $\varepsilon(a)\subseteq U\cap P_i$ for some $a\in P$, $i<\cf(P)$ and $\varepsilon\in\mathcal{B}_i$. Find $\varepsilon_0\in\mathcal{B}_i$ with $\varepsilon_0^{-1}\subseteq \varepsilon$ and $d\in D_{i,\varepsilon_0}\subseteq D$ with $a\in\varepsilon_0(d)$, so $d\in\varepsilon(a)\subseteq U$, concluding $D\cap U\neq \emptyset$. As $U$ is arbitrary, $D$ is dense.  
\end{proof}
\end{coro}
Arguing similarly for global logic topologies, we can reduce the number of parameters from $2^{|A|+|\lang|}$ to $|A|+|\lang|+\cf(P)$.
\begin{prop}\label{p:reduce number of parameters} Let $P$ be a piecewise $A${\hyp}hyperdefinable set with a global logic topology. Then, there is $B\subseteq P\cup A$ with $|B|\leq |A|+|\lang|+\cf(P)$ such that the $B${\hyp}logic topology of $P$ is its global logic topology.
\begin{proof} Say $P=\underrightarrow{\lim}_{i\in I}\, P_i$ and $P_i=\rfrac{X_i}{E_i}$ with $|I|=\cf(P)$. Write $\lambda=|\lang|+|A|+\cf(P)$. For $i\in I$, $\Delta_i\coloneqq \{(x,x)\sth x\in P_i\}$ is $\bigwedge_{A^*}${\hyp}definable as $\underline{\Delta}_i=\underline{E}_i$. Say $\underline{\Delta}_i=\bigwedge_{j\in  \lambda} \varphi_j(x,y)$ with $\varphi_j(x,y)\in\For(\lang(A^*))$. Write $V_j\coloneqq \quot_{P_i\times P_i}[\varphi_j(\mathfrak{M})]$ and $U_j\coloneqq P_i\times P_i\setminus \quot_{P_i\times P_i}[\neg\varphi_j(\mathfrak{M})]$. Then, $\Delta_i\subseteq U_j\subseteq V_j$ where $V_j$ is closed and $U_j$ is open in the global logic topology of $P_i\times P_i$. By compactness of $P_i\times P_i$, there is $\varepsilon_j$ positive $\bigwedge${\hyp}definable binary relation such that $\varepsilon_j\subseteq U_j\subseteq V_i$. It follows that $\mathcal{B}_i\coloneqq \{\varepsilon_j\}_{j\in \lambda}$ forms a uniformity base of $P_i$ with the global logic topology. 

By compactness of $P_i$, for each $\varepsilon\in \mathcal{B}_i$, there is $D_{\varepsilon,i}\subseteq P_i$ finite such that $P_i\subseteq\bigcup_{a\in D_{\varepsilon,i}} \varepsilon(a)$. Take $D_i=\bigcup_{\varepsilon\in\mathcal{B}_i} D_{\varepsilon,i}$ and $D=\bigcup_{i\in I}D_i$, so $|D|\leq \lambda$. It follows that $D_i$ is dense in the global logic topology of $P_i$. 

Take $B=D\cup A\subseteq \mathrm{bdd}(A)$, so $|B|\leq \lambda$. We claim that the $B${\hyp}logic topology of $P$ is its global logic topology. Indeed, take $a\in P$ and $\sigma\in \Aut(\mathfrak{M}/B)$ arbitrary and, aiming a contradiction, suppose $\sigma(a)\neq a$, so $\sigma^{-1}(a)\neq a$. Pick $i\in I$ such that $a\in P_i$. Note that $\sigma[P_i]=P_i$. By Hausdorffness of $P_i$ in the global logic topology, there are $U$ and $V$ such that $a\in U\subseteq V\subseteq P_i$ with $\sigma^{-1}(a)\notin V$, where $U$ is open and $V$ is closed in the global logic topology of $P_i$. Now, $D_i\cap U=\sigma[D_i\cap U]\subseteq \sigma[V]$. As $U$ is open and $D_i$ is dense in $P_i$, $U\cap D_i$ is dense in $U$. As $\sigma[V]$ is closed in the global logic topology, we conclude $a\in U\subseteq \sigma[V]$. Therefore, $\sigma^{-1}(a)\in V$, getting a contradiction. 
\end{proof}
\end{prop}
\subsection{Examples}
Most of the examples of locally hyperdefinable sets come from the following two basic remarks, which were, in fact, already known. First, note that any piecewise $A${\hyp}definable set is trivially locally $A${\hyp}definable. Secondly, note that if $P$ is locally $A${\hyp}hyperdefinable and $E$ is a piecewise bounded $\bigwedge_A${\hyp}definable equivalence relation on $P$, then $\rfrac{P}{E}$ is locally hyperdefinable by \cref{p:mapping locally hyperdefinable sets}. \smallskip

\begin{ex}\label{e:example 1} A classical example is the field of real numbers $\R$ with its usual topology, which is a locally hyperdefinable set of countable cofinality in the theory of real closed fields. It is explicitly given as $\rfrac{O(1)}{o(1)}$ with $O(1)=\bigcup_{n\in\N} [-n,n]$ and $o(1)=\bigcap_{n\in\N} \left[-\sfrac{1}{n},\sfrac{1}{n}\right]$. Furthermore, up to isomorphism of piecewise hyperdefinable sets, $\rfrac{O(1)}{o(1)}$ is the unique representation of $\R$ as a locally hyperdefinable set. Indeed, by \cref{p:functions logic topology}(\ref{itm:functions logic topology:homeomorphism iff piece isomorphism}), any two locally hyperdefinable sets homeomorphic with the logic topologies are isomorphic as piecewise hyperdefinable sets. 

However, note that the real numbers with the usual topology can be represented as a piecewise hyperdefinable set (non{\hyp}locally hyperdefinable) in other non{\hyp}isomorphic ways. For instance, consider the direct system of all compact subsets of $\rfrac{O(1)}{o(1)}$ with empty interior in the global logic topology with the natural inclusion maps. Using that the topology is first{\hyp}countable, it is easy to note that this is a coherent covering. Therefore, the resulting direct limit with the global logic topology is a piecewise hyperdefinable set homeomorphic to $\R$ with the usual topology. However, it is not locally hyperdefinable, so it is not isomorphic to $\rfrac{O(1)}{o(1)}$.
\end{ex}


\begin{ex} \label{e:example 2}
More generally, any topological manifold $X$ (i.e. a locally euclidean Hausdorff topological space) is a locally hyperdefinable set in the usual theory of real closed fields. Indeed, for any $m$, $\R^m$ is locally hyperdefinable, so every compact subset of $\R^m$ is hyperdefinable, concluding that the compact charts of $X$ are hyperdefinable. Using now \cref{p:functions logic topology 2}, the chart changing maps are $\bigwedge${\hyp}definable, so the finite unions of compact chart neighbourhoods are hyperdefinable. It follows then that the whole manifold is locally hyperdefinable. If it is also second countable, then it has countable cofinality too. \end{ex}

\begin{ex} \label{e:example 3}
As $\Q$ with the usual metric topology is not locally compact, it cannot be given as a locally hyperdefinable set. However, it is possible to give it as a piecewise hyperdefinable set in the theory of real closed fields. Indeed, using that $\Q$ is first countable, it is clear that $\Q$ is compactly generated. Every compact subset of $\Q$ is a hyperdefinable subset being $\bigwedge${\hyp}definable in $\R=\rfrac{O(1)}{o(1)}$. In other words, $\Q$ is the direct limit of all the $\bigwedge${\hyp}definable subsets of $\R=\rfrac{O(1)}{o(1)}$ contained in $\Q$ with the standard inclusion maps. Note that it has uncountable cofinality.
\end{ex}

\begin{ex} \label{e:example 4}
More generally, any first countable Hausdorff topological space $X$ can be given as a piecewise hyperdefinable set with a global logic topology in the theory of real closed fields. Indeed, let $\mathcal{C}$ be the family of countable compact subsets of $X$. As $X$ is first countable, $\mathcal{C}$ is coherent. Now, by Mazurkiewicz{\hyp}Sierpi\'{n}ski Theorem \cite[Theorem 4]{milliet2011remark}, every countable compact Hausdorff topological space is homeomorphic to a countable successor ordinal with the order topology. On the other hand, by induction, we easily see that every countable ordinal with the order topology is homeomorphic to a subset of $\Q$. Therefore, for every $A\in \mathcal{C}$, there is a compact subset $P_A\subseteq \Q$ such that $A$ with the subspace topology is homeomorphic to $P_A$ with the subspace topology. For each $A\in\mathcal{C}$, pick $\eta_A: A\rightarrow P_A$ a homeomorphism. For $A,B\in \mathcal{C}$ with $A\subseteq B$, take $\varphi_{BA}=\eta_B\circ \eta_A^{-1}:\ P_{A}\rightarrow P_{B}$. Now, as noted in the previous example, every compact subset of $\Q$ is homeomorphic to a hyperdefinable set with a global logic topology in the theory of real closed fields. Also, for $A,B\in\mathcal{C}$ with $A\subseteq B$, $\varphi_{BA}$ is continuous, so it is $\bigwedge${\hyp}definable by \cref{p:functions logic topology 2}. Then, we conclude that $X$ is homeomorphic to $\underrightarrow{\lim}_{A\in \mathcal{C}}\, P_A$ with the global logic topology.
\end{ex}

\begin{ex} \label{e:example 5} For a countably piecewise hyperdefinable set that is not locally hyperdefinable, consider the infinite countable rose, i.e. the infinite countable bouquet of circles. This is $\rfrac{\R}{\sim}$ with the equivalence relation $x\sim y\Leftrightarrow x,y\in\Z\vee x=y$. Note that $\R$ is countably piecewise hyperdefinable and the relation $\sim$ is piecewise $\bigwedge${\hyp}definable in the theory of real closed fields. Thus, the infinite countable rose is a countably piecewise hyperdefinable set. It is not locally hyperdefinable and not first countable, so it is not pseudo{\hyp}metrisable.
\end{ex}

\begin{ex} \label{e:example 6}
For a piecewise hyperdefinable set which is not normal, consider $\R$ with the rational sequence topology. For $i\in\R$, pick a sequence $(i_n)_{n\in\N}$ of rational numbers converging to $i$. When $i\in \Q$, take $i_n=i$ for every $n\in\N$. Take $P_i=\{i_n\sth n\in\N\}\cup\{i\}$ for $i\in\R$. For a finite subset $I\subseteq \R$, take $P_I=\bigcup_{i\in I} P_i$. Note that, for any finite $I\subseteq \R$, $P_I$ is compact in $\R$ with the usual topology, so each $P_I$ is hyperdefinable in the language of real closed fields. Take $P=\underrightarrow{\lim}\ P_I$ with the usual inclusion maps. We now check that $P$ with the global logic topology is homeomorphic to $\R$ with the rational sequence topology, i.e. the topology given by the local bases of open neighbourhoods $U_n(i)\coloneqq \{i_k\sth k\geq n\}\cup \{i\}$ for $i\in\R$. 

Note first that $U_n(i)$ is open in $P$ for each $i\in\R$ and $n\in\N$. Obviously, $U_n(i)\cap P_i$ is open in $P_i$. For $j\neq i$, as $(i_n)_{n\in\N}$ converges to $i$ and $(j_n)_{n\in\N}$ converges to $j$, we conclude that there is $m\in\N$ such that $U_n(i)\cap P_j\subseteq \{i_n\sth n\leq m\}$. Thus, $U_n(i)\cap P_j$ is open in $P_j$. 

On the other hand, suppose $U\subseteq P$ is open in $P$. Take $i\in U$. As $U\cap P_i$ is open in $P_i$, there is $n\in\N$ such that $U_n(i)\subseteq U\cap P_i$, so $U_n(i)\subseteq U$. As $i\in U$ is arbitrary, we conclude that $U$ is open in the rational sequence topology. 

By Jone's Lemma \cite[Lemma 15.2]{willard1970general}, $P$ is not normal. Indeed, $\Q$ is dense, $\R\setminus \Q$ is discrete and closed and $|\R\setminus\Q|\geq 2^{|\Q|}$.
\end{ex}

For a counterexample where the product topology of global logic topologies is not the global logic topology on the product, see \cref{e:example no topological group}. We have no counterexample of a piecewise hyperdefinable set with a non{\hyp}Hausdorff global logic topology. We have no counterexample of a countably piecewise $A${\hyp}hyperdefinable set that is not locally $A${\hyp}hyperdefinable but has a locally compact $A${\hyp}logic topology.

\section{Piecewise hyperdefinable groups}
In this section we study the particular case of piecewise hyperdefinable groups. Our main aim is to find sufficient and necessary conditions to conclude when they are locally compact topological groups with the logic topology. Then, we will discuss how to extend the classical Gleason{\hyp}Yamabe Theorem and some related results to piecewise hyperdefinable groups. We start with an introduction recalling some fundamental facts about topological groups.

\subsection{Preliminaries on topological groups} 
Recall that a \emph{topological group} is a Hausdorff topological space with a group structure whose group operations are continuous. 
\begin{rmk} \label{o:remarks topological groups} Let $G$ be a topological group. The following are the basic fundamental facts that we need:
\begin{enumerate}[label={\rm{(\arabic*)}}, ref={\rm{\arabic*}}, wide]
\item \label{itm:remarks topological groups:quotient} Let $H\trianglelefteq G$ be closed. Then, $\rfrac{G}{H}$ is a topological group and $\pi_{G/H}:\ G\rightarrow \rfrac{G}{H}$ is a continuous and open surjective homomorphism. Furthermore, if $H$ is compact, then $\pi_{G/H}$ is also a closed map and has compact fibres. In particular, it is proper by \cite[Theorem 3.7.2]{engelking1989general}.
\item \label{itm:remarks topological groups:open subgroup} Let $H\leq G$ be an open subgroup. Then, $H$ is also closed.
\item \label{itm:remarks topological groups:connected component} The connected component $G^0$ of the identity is a normal closed subgroup of $G$. If $G$ is locally connected (e.g. a Lie group), then $G^0$ is also open.
\item \label{itm:remarks topological groups:closed isomorphism theorem} {(Closed Isomorphism Theorem)} Let $f:\ G\rightarrow H$ be a continuous and closed surjective homomorphism between topological groups. Then, for any closed subgroup $S\trianglelefteq K\coloneqq \ker(f)$ with $S\trianglelefteq G$, the map $f_S:\ \rfrac{G}{S}\rightarrow H$ defined by $f=f_S\circ\pi_{G/S}$ is a continuous, closed and open homomorphism. In particular, $f_K:\ \rfrac{G}{K}\rightarrow H$ is an isomorphism of topological groups and $f$ is an open map.
\end{enumerate}
\end{rmk}
In the theory of topological groups, a \emph{Yamabe pair} of a topological group $G$ is a pair $(K,H)$ with $K\trianglelefteq H\leq G$ such that $K$ is compact, $H$ is open and $L=\rfrac{H}{K}$ is a finite dimensional Lie group. We say that $H$ is the \emph{domain}, $K$ is the \emph{kernel} and $L$ is the \emph{Lie core}. We write $\pi_{H/K}:\ H\rightarrow L$ for the quotient map. A Lie group is a \emph{Lie core} of $G$ if it is isomorphic, as topological group, to the Lie core of some Yamabe pair of $G$.
\begin{rmk} Let $G$ be a topological group and suppose that it has a Yamabe pair $(K,H)$ with Lie core $L$. By \cref{o:remarks topological groups}(\ref{itm:remarks topological groups:quotient}), $\pi_{H/K}:\ H\rightarrow L$ is a continuous, open, closed and proper surjective group homomorphism. In particular, as $L$ is locally compact, $H$ must be a locally compact topological group too. Since $H$ is open in $G$, we conclude that $G$ is locally compact as well.
\end{rmk} 
The following celebrated theorem, claiming that every locally compact topological group has Lie cores, is usually considered the solution to Hilbert's fifth problem.
\begin{theo}[Gleason{\hyp}Yamabe] \label{t:gleason yamabe} 
Let $G$ be a locally compact topological group and $U\subseteq G$ a neighbourhood of the identity. Then, $G$ has a Yamabe pair $(K,H)$ with $K\subseteq U$. In particular, a topological group has a Lie core if and only if it is locally compact.
\begin{proof} The original papers are \cite{gleason1952groups,montgomery1952small,yamabe1953conjecture,yamabe1953generalization}. A complete proofs can also be found in \cite{montgomery1955topological,tao2014hilbert}. Model{\hyp}theoretic treatments can be found in \cite{hirschfeld1990nonstandard,van2015hilbert}. 
\end{proof}
\end{theo}
In this paper we mainly use this classical version of Gleason{\hyp}Yamabe \cref{t:gleason yamabe}. Alternatively, we can use the following variation proved in \cite[Theorem 1.25]{carolino2015structure} which provides some extra control over some parameters. Recall that two subsets of a group are \emph{$k${\hyp}commensurable} if $k$ left translates of each one cover the other. Recall that a \emph{$k${\hyp}approximate subgroup} is a symmetric subset $X$ which is $k${\hyp}commensurable with its set of pairwise products $X^2$. 
\begin{theo}[Gleason{\hyp}Yamabe{\hyp}Kreitlon Theorem]  \label{t:gleason yamabe carolino}
There are functions $c:\ \N\rightarrow \N$ and $d:\ \N\rightarrow\N$ such that the following holds: 

Let $G$ be a locally compact topological group and $U\subseteq G$ an open precompact $k${\hyp}approximate subgroup for some $k\in\N$. Then, $G$ has a Yamabe pair $(K,H)$ with $K\subseteq U^4$ such that $\rfrac{H}{K}$ is a Lie group of dimension at most $d(k)$ and $H\cap U^4$ generates $H$ and is $c(k)${\hyp}commensurable to $U$.
\end{theo}
A Yamabe pair $(K',H')$ is \emph{smaller than or equal to} $(K,H)$ if $K\trianglelefteq K'\trianglelefteq H'\leq H$. A Yamabe pair is \emph{minimal} if it has no smaller ones. A Lie core is \emph{minimal} if it is the Lie core of some minimal Yamabe pair. In other words, let us define an \emph{aperiodic} topological group to be a topological group without non{\hyp}trivial compact normal subgroups. Then, by \cref{o:remarks topological groups}(\ref{itm:remarks topological groups:quotient}), a Lie core is minimal if and only if it is an aperiodic connected Lie core. The following basic proposition implies that every Yamabe pair has a smaller or equal minimal Yamabe pair.
\begin{lem} \label{l:maximal compact normal subgroup}
Every connected Lie group has a unique maximal compact normal subgroup.  
\begin{proof} By Cartan{\hyp}Iwasawa{\hyp}Malcev Theorem \cite[Chapter XV Theorem 3.1]{hochschild1965structure}, there is a maximal compact subgroup $T$, and every compact subgroup is contained in a conjugate of it. Hence, $\bigcap_{g\in G} gTg^{-1}$ is the unique maximal compact normal subgroup. 
\end{proof}
\end{lem}
\begin{rmk} A different proof of the previous result, without using Cartan{\hyp}Iwasawa{\hyp}Malcev Theorem, was explained in \cite{hrushovski2011stable}.
\end{rmk}
\begin{coro} \label{c:minimal yamabe pairs} Let $G$ be a topological group and $(K_1,H_1)$ a Yamabe pair. Then, there is a minimal Yamabe pair $(K,H)$ smaller than or equal to $(K_1,H_1)$. Furthermore, for any clopen subset $U$ containing $K_1$, we have $H\subseteq U^2$.    
\begin{proof} Write $\pi_1\coloneqq \pi_{H_1/K_1}:\ H_1\rightarrow L_1$ for the quotient map to the Lie core of $(K_1,H_1)$. Let $\widetilde{L}\subseteq L_1$ be the topological connected component of the identity. As Lie groups are locally connected, $\widetilde{L}$ is open by \cref{o:remarks topological groups}(\ref{itm:remarks topological groups:connected component}). Let $\widetilde{K}\trianglelefteq \widetilde{L}$ be its maximal compact normal subgroup, given by \cref{l:maximal compact normal subgroup}. Take $H\coloneqq \pi^{-1}_1[\widetilde{L}]$ and $K\coloneqq \pi^{-1}_1[\widetilde{K}]$. Then, by \cref{o:remarks topological groups}(\ref{itm:remarks topological groups:quotient}) and the Closed Isomorphism Theorem (\cref{o:remarks topological groups}(\ref{itm:remarks topological groups:closed isomorphism theorem})), $(K,H)$ is a minimal Yamabe pair of $G$ smaller than or equal to $(K_1,H_1)$. Finally, if $U$ is clopen and $K_1\subseteq U$, $\pi_1[U]$ is clopen with $1\in \pi_1[U]$ as $\pi_1$ is open and closed by \cref{o:remarks topological groups}(\ref{itm:remarks topological groups:quotient}). Thus, $\widetilde{L}\subseteq \pi_1[U]$ as $\widetilde{L}$ is connected, concluding that $H\subseteq UK_1\subseteq U^2$. 
\end{proof}
\end{coro}
Two Yamabe pairs $(K,H)$ and $(K',H')$ of $G$ with Lie cores $\pi\coloneqq \pi_{H/K}:\ H\rightarrow L$ and $\pi'\coloneqq \pi_{H'/K'}:\ H'\rightarrow L'$ are \emph{equivalent} if the map $\eta:\ \pi(h)\mapsto \pi'(h)$ for $h\in H\cap H'$ is a well{\hyp}defined isomorphism of topological groups between $L$ and $L'$. Equivalently, by the Closed Isomorphism Theorem (\cref{o:remarks topological groups}(\ref{itm:remarks topological groups:closed isomorphism theorem})), $(K,H)$ and $(K',H')$ are equivalent if and only if $H\cap K'\subseteq K$, $H'\cap K\subseteq K'$, $(H\cap H')K=H$ and $(H\cap H')K'=H'$. It follows that minimal Yamabe pairs are unique up to equivalence:
\begin{lem} \label{l:lemma uniqueness minimal lie core} Let $G$ be a locally compact topological group and $(K_1,H_1)$ and $(K_2,H_2)$ two minimal Yamabe pairs with Lie cores $\pi_1\coloneqq \pi_{H_1/K_1}:\ H_1\rightarrow L_1$ and $\pi_2\coloneqq \pi_{H_2/K_2}:\ H_2\rightarrow L_2$:
\begin{enumerate}[label={\rm{(\arabic*)}}, ref={\rm{\arabic*}}, wide]
\item \label{itm:lemma uniqueness minimal lie core:reduce domain} Let $H'\leq H_1$ be an open subgroup. Then, $(K_1\cap H',H')$ is a minimal Yamabe pair of $G$ equivalent to $(K_1,H_1)$ and $[K_1:K_1\cap H']$ is finite.
\item \label{itm:lemma uniqueness minimal lie core:smaller kernel} $K_1\subseteq K_2\Leftrightarrow H_1\subseteq H_2\Leftrightarrow K_2\cap H_1=K_1$. In particular, $K_1=K_2$ if and only if $H_1=H_2$.
\item \label{itm:lemma uniqueness minimal lie core:intersection} $(K_1\cap K_2,H_1\cap H_2)$ is a minimal Yamabe pair with $K_1\cap H_2=K_1\cap K_2=K_2\cap H_1$. In particular, $[K_1:K_1\cap K_2]$ and $[K_2:K_1\cap K_2]$ are finite.
\end{enumerate}
\begin{proof}
\begin{enumerate}[label={\rm{(\arabic*)}}, wide]
\item[\hspace{-1.6em}\setcounter{enumi}{1}\theenumi] By connectedness, $\pi_1[H']=L_1$. Thus, by the Closed Isomorphism Theorem (\cref{o:remarks topological groups}(\ref{itm:remarks topological groups:closed isomorphism theorem})), we conclude that $(K_1\cap H',H')$ is a minimal Yamabe pair of $G$ equivalent to $(K_1,H_1)$. Finally, as $K_1\cap H'$ is an open subset of $K_1$, by compactness, $[K_1:K_1\cap H']$ is finite.
\item We already have $H_1\subseteq H_2\Rightarrow K_2\cap H_1=K_1\Rightarrow K_1\subseteq K_2$ by the previous point (\ref*{itm:lemma uniqueness minimal lie core:reduce domain}). On the other hand, by connectedness, $\pi_1[H_1\cap H_2]=L_1$. If $K_1\subseteq K_2$, then $K_1\leq H_1\cap H_2$, so $H_1=\pi^{-1}_1[\pi_1[H_1\cap H_2]]=H_1\cap H_2\subseteq H_2$.
\item By point (\ref*{itm:lemma uniqueness minimal lie core:reduce domain}), $(K_1\cap H_2,H_1\cap H_2)$ and $(K_2\cap H_1,H_2\cap H_1)$ are minimal Yamabe pairs with $[K_1:K_1\cap H_2]$ and $[K_2:H_1\cap K_2]$ finite. By point (\ref*{itm:lemma uniqueness minimal lie core:smaller kernel}), $K_1\cap K_2=K_1\cap H_2=K_2\cap H_1$. \qedhere
\end{enumerate}
\end{proof}
\end{lem}
As an immediate consequence of the previous proposition we get the following corollary:
\begin{coro} \label{c:uniqueness minimal lie core} Every locally compact topological group has a unique minimal Yamabe pair up to equivalence.
\end{coro}
The previous uniqueness statement implies that the minimal Lie core $L$ is unique up to isomorphism of topological groups, but it is far stronger than that. Indeed, it also says that there is a \emph{global minimal Lie core map} extending all the minimal Yamabe pairs which is unique up to isomorphisms of $L$:
\begin{prop}\label{p:global minimal lie core map}
Let $G$ be a locally compact topological group and $L$ its minimal Lie core. Let $D_L$ be the union of all the domains of minimal Yamabe pairs of $G$. Then, there is a map $\pi_L:\ D_L\rightarrow L$ such that, for any minimal Yamabe pair $(K,H)$ of $G$, $\pi_{L\mid H}$ is a continuous, closed, open and proper surjective homomorphism with kernel $K$. Furthermore, $\pi_L$ is unique up to isomorphisms of $L$.
\begin{proof} Let $\mathcal{Y}$ be the set of all minimal Yamabe pairs of $G$ and fix any minimal Yamabe pair $(K_0,H_0)\in \mathcal{Y}$ and $L=\rfrac{H_0}{K_0}$. Now, for any Yamabe pair $(K,H)$, we define $\pi_{L\mid H}\coloneqq \eta_{(K,H)}\circ \pi_{H/K}:\ H\rightarrow L$ where $\eta_{(K,H)}:\ \rfrac{H}{K}\rightarrow L$ is the canonical isomorphism given by the equivalence between $(K,H)$ and $(K_0,H_0)$. Take any $(K,H),(K',H')\in \mathcal{Y}$. By \cref{l:lemma uniqueness minimal lie core,c:uniqueness minimal lie core}, $(K\cap K',H\cap H')$ is a minimal Yamabe pair and equivalent to $(K,H)$, $(K',H')$ and $(K_0,H_0)$. For any $h\in H\cap H'$, as $(K\cap K',H\cap H')$ is equivalent to $(K_0,H_0)$, there is $h_0\in H\cap H'\cap H_0$ such that $h^{-1}h_0\in K\cap K'\cap K_0$. Thus, $\pi_{L\mid H}(h)=\pi_{L\mid H}(h_0)=\pi_{H_0/K_0}(h_0)=\pi_{L\mid H'}(h_0)=\pi_{L\mid H'}(h)$. Take $D_L=\bigcup_{\mathcal{Y}} H$ and define $\pi_L=\bigcup_{\mathcal{Y}} \pi_{L\mid H}$. By \cref{o:remarks topological groups}(\ref{itm:remarks topological groups:quotient}), we get a global map $\pi_L:\ D_L\rightarrow L$ such that $\pi_{L\mid H}:\ H\rightarrow L$ is a continuous, closed, open and proper surjective homomorphism with kernel $K$ for each minimal Yamabe pair $(K,H)$. 

Suppose $\pi'_L:\ D_L\rightarrow L$ is any other map such that $\pi'_{L\mid H}$ is a continuous, closed, open and proper surjective group homomorphism with kernel $K$ for any minimal Yamabe pair $(K,H)$. Then, by the Closed Isomorphism Theorem (\cref{o:remarks topological groups}(\ref{itm:remarks topological groups:closed isomorphism theorem})), we get an isomorphism $\eta:\ \rfrac{H_0}{K_0}\rightarrow L$ such that $\pi'_{L\mid H_0}=\eta\circ\pi_{L\mid H_0}$. Now, for any $(K,H)\in\mathcal{Y}$ and $g\in H$, there is $g_0\in H\cap H_0$ such that $g\in g_0K$. Then, $\pi'_L(g)=\pi'_L(g_0)=\eta\circ\pi_L(g_0)=\eta\circ\pi_L(g)$, concluding that $\pi'_L=\eta\circ\pi_L$. 
\end{proof}
\end{prop}
Note that $D_L=\mathrm{Dom}(\pi_L)$ is the union of all the domains of minimal Yamabe pairs of $G$ and $\ker(\pi_L)\coloneqq \pi_L^{-1}(1)$ is the union of all the kernels of minimal Yamabe pairs of $G$. Consequently, $D_L$ and $\ker(\pi_L)$ are invariant by any automorphism of $G$ as topological group. In particular, both are normal sets (i.e. conjugate invariant sets).\smallskip

Among all the minimal Yamabe pairs, it could be natural to look for the ones with maximal domain. We have the following criterion.  \begin{prop} \label{p:minimal yamabe pairs with maximal domain} Let $G$ be a locally compact topological group and $(K,H)$ a minimal Yamabe pair of $G$. Let $K'$ be a compact subgroup of $G$ with $K\leq K'$ such that $H$ normalises $K'$ (i.e. $hK'=K'h$ for any $h\in H$). Then, $K=H\cap K'$, $[K':K]$ is finite and $(K',H')$ is a minimal Yamabe pair of $G$ with $H'=K'H$. Furthermore, $H'$ is a finite union of cosets of $H$. 

In particular, $(K,H)$ is a minimal Yamabe pair with maximal domain if and only if there is no compact subgroup $K'\leq G$ normalised by $H$ with $K< K'$.
\begin{proof} Clearly, $K\leq K'\cap H\trianglelefteq H$ is compact. Then, as $K$ is the maximal compact normal subgroup of $H$, we conclude $K'\cap H=K$ and $\rfrac{H}{K'\cap H}=\rfrac{H}{K}=L$. As $K=K'\cap H$ is open in $K'$ compact, $[K':K]$ is finite. Take $\Delta\subseteq K'$ finite such that $K'=\Delta K$. Note that $\Delta H=K'H$ is a clopen subgroup and $K'\trianglelefteq K'H$. Write $H'=\Delta H$. Then, $\pi_{H'/K'\ \mid H}:\ H\rightarrow \rfrac{H'}{K'}$ is a continuous and closed onto homomorphism. Therefore, by the Closed Isomorphism Theorem (\cref{o:remarks topological groups}(\ref{itm:remarks topological groups:closed isomorphism theorem})), $\rfrac{H'}{K'}$ and $\rfrac{H}{K'\cap H}=L$ are isomorphic. That concludes that $(K',H')$ is also a minimal Yamabe pair of $G$. Using also \cref{l:lemma uniqueness minimal lie core}(\ref{itm:lemma uniqueness minimal lie core:reduce domain},\ref{itm:lemma uniqueness minimal lie core:smaller kernel}), we conclude that this gives a necessary and sufficient condition for the maximality of the domain. 
\end{proof}
\end{prop}
Similarly, it is natural to look at minimal Yamabe pairs with minimal kernel. In this case, this question is related to the connected component of $G$.\smallskip

Recall that the \emph{quasicomponent} of a point in a topological space is the intersection of all its clopen neighbourhoods. By definition, quasicomponents are closed sets containing the connected components. In locally connected spaces, connected components are clopen, so quasicomponents and connected components coincide. Similarly, in every compact Hausdorff space, connected components and quasicomponents coincide \cite[Lemma 29.6]{willard1970general}. In general, however, they may be different --- even for locally compact Hausdorff topological spaces. 

In a topological group $G$, as the inversion, the conjugations and the translations are homeomorphisms, the connected component $G^0$ and the quasicomponent $G^{\mathrm{qs}}$ of the identity are both normal closed subgroups of $G$. When $G$ is locally compact, it is a well{\hyp}known fact that $G^0=G^{\mathrm{qs}}$ is the intersection of all the open subgroups of $G$ \cite[Theorem 2.1.4(b)]{dikranjan1990topological}.
Hence, we conclude the following criterion for the existence of a minimal Yamabe pair with minimal kernel. 
\begin{prop} \label{p:minimal yamabe pair with minimal domain}
Let $G$ be a locally compact topological group. Then, there is a minimal Yamabe pair with minimal kernel if and only if $G^0$ is open (i.e. $G$ is locally connected). Furthermore, in that case, for any other minimal Yamabe pair $(K,H)$, $(K\cap G^0,G^0)$ is the minimal Yamabe pair of $G$ with minimal kernel.
\begin{proof} Suppose that $(K,H)$ is a minimal Yamabe pair with minimal kernel. As $H$ is clopen by \cref{o:remarks topological groups}(\ref{itm:remarks topological groups:open subgroup}), we have that $G^0\subseteq H$. On the other hand, for any other open subgroup $H'\leq G$, by \cref{l:lemma uniqueness minimal lie core}(\ref{itm:lemma uniqueness minimal lie core:reduce domain}), we have that $(K\cap H',H\cap H')$ is a minimal Yamabe pair. As $(K,H)$ is the one with minimal kernel, it follows that $H\cap H'=H$, so $H\subseteq H'$. As $G^0$ is the intersection of all the open subgroups, we conclude that $G^0=H$. Conversely, suppose that $G^0$ is open. Then, by \cref{l:lemma uniqueness minimal lie core}(\ref{itm:lemma uniqueness minimal lie core:reduce domain}), for any minimal Yamabe pair $(K,H)$, we have that $(K\cap G^0,G^0)$ is a minimal Yamabe pair. Thus, for any minimal Yamabe pairs $(K,H)$ and $(K',H')$, we have that $(K\cap G^0,G^0)$ and $(K'\cap G^0,G^0)$ are minimal Yamabe pairs, and so $K\cap G^0=K'\cap G^0$ by \cref{l:lemma uniqueness minimal lie core}(\ref{itm:lemma uniqueness minimal lie core:smaller kernel}). Therefore, $(K\cap G^0,G^0)$ is the minimal Yamabe pair with minimal kernel.  
\end{proof}
\end{prop}
In general, even if $G^0$ is not open, a similar conclusion is ``asymptotically'' true:
\begin{prop} \label{p:asymptotic minimal yamabe pair with minimal domain}
Let $G$ be a locally compact topological group and $L$ its minimal Lie core. Then, the restriction to $G^0$ of the global minimal Lie core map $\pi_{L\mid G^0}:\ G^{0}\rightarrow L$ is a continuous, open, closed and proper surjective group homomorphism.
\begin{proof} Take $(K,H)$ minimal Yamabe pair and $\pi_{L\mid H}:\ H\rightarrow L$. By \cref{p:global minimal lie core map}, $\pi_{L\mid H}$ is a continuous, open, closed and proper surjective group homomorphism. By definition, $G^0\leq H$. Thus, consider the restriction $\pi_{L\mid G^0}:\ G^0\rightarrow L$. As $G^0$ is a closed subgroup, $\pi_{L\mid G^0}$ is also a continuous, closed and proper group homomorphism. It remains to show that it is onto and open. Let $b\in L$. We want to show that $\pi^{-1}_{L\mid H}(b)\cap G^0\neq \emptyset$. Let $H'\leq G$ be a clopen subgroup such that $H'\leq H$. Then, $\pi_{L\mid H}[H']$ is clopen in $L$. As $L$ is connected, we get that $\pi^{-1}_{L\mid H}(b)\cap H'\neq \emptyset$. Since $\pi^{-1}_{L\mid H}(b)$ is compact and $H'$ is arbitrary, we conclude that $\pi^{-1}_{L\mid H}(b)\cap G^0\neq \emptyset$. Therefore, $\pi_{L\mid G^0}:\ G^0\rightarrow L$ is onto. We conclude that it is also open by the Closed Isomorphism Theorem (\cref{o:remarks topological groups}(\ref{itm:remarks topological groups:closed isomorphism theorem})). 
\end{proof}
\end{prop}
\subsection{Local compactness and generic pieces}
A \emph{piecewise $A${\hyp}hyperdefinable group} is a group whose universe is piecewise $A${\hyp}hyperdefinable and whose operations are piecewise bounded $\bigwedge_A${\hyp}definable. 
\begin{ex}
Let $G$ be a definable group and $X\subseteq G$ a symmetric definable subset. Then, the subgroup $H\leq G$ generated by $X$ is a countably piecewise definable group. If $K\trianglelefteq H$ is a piecewise $\bigwedge${\hyp}definable normal subgroup, the quotient $\rfrac{H}{K}=\underrightarrow{\lim} \rfrac{X^n}{K}$ is a countably piecewise hyperdefinable group too. If $K\subseteq X^n$ for some $n$, then $\rfrac{H}{K}$ is also locally hyperdefinable. This corresponds to the case studied in \cite{hrushovski2011stable}.   
\end{ex}
\begin{rmk} Piecewise $\bigwedge${\hyp}definable subgroups of piecewise hyperdefinable groups are piecewise hyperdefinable groups. The quotient of a piecewise hyperdefinable group by a normal piecewise $\bigwedge${\hyp}definable subgroup is a piecewise hyperdefinable group.
\end{rmk}
Note that the group operations are continuous between the logic topologies by \cref{p:functions logic topology}(\ref{itm:functions logic topology:continuous}). However, the product topology and the logic topology may differ, so piecewise hyperdefinable groups with the logic topologies do not need to be topological groups. 
\begin{prop} \label{p:small piecewise hyperdefinable groups and translations} Let $G$ be a piecewise hyperdefinable group with a global logic topology. Then, every translation is a homeomorphism in its global logic topology.
\begin{proof} Trivial by \cref{p:functions logic topology}(\ref{itm:functions logic topology:continuous}) and \cref{p:uniqueness of hausdorff logic topologies}. 
\end{proof}
\end{prop}
\begin{rmk} Groups with a $\mathrm{T}_1$ topology such that every translation is continuous are called \emph{semitopological groups} --- see \cite{husain2018introduction} for an introduction to semitopological groups.
\end{rmk} 
\begin{theo}\label{t:countably piecewise hyperdefinable groups and topological groups}
Let $G$ be a countably piecewise hyperdefinable group with a global logic topology. Then, $G$ is a topological group with the global logic topology. 
\begin{proof} Clear by \cref{p:functions logic topology}(\ref{itm:functions logic topology:continuous}) and \cref{p:product of countably piecewise hyperdefinable sets}. 
\end{proof}
\end{theo}
\begin{ex} \label{e:example no topological group}
We show now an example of a piecewise hyperdefinable group with a global logic topology that is not a topological group. We simply adapt the fundamental example given in \cite[Example 1.2]{tatsuuma1998group} to the case of piecewise hyperdefinable groups. 

First, recall that $\Q^n$ with the usual topology is piecewise hyperdefinable with a global logic topology in the theory of real closed fields by \cref{e:example 4}. Now, the inclusion $\psi_{n,m}:\ \Q^m\rightarrow \Q^n$ given by $\psi(x)=(x,0,\ldots,0)$ for $n>m$ is a piecewise bounded $\bigwedge${\hyp}definable $1${\hyp}to{\hyp}$1$ map. Also, the set of pairwise sums of two compact countable subsets of $\Q^n$ is a compact countable subset of $\Q^n$, so $+$ is a piecewise bounded $\bigwedge${\hyp}definable map. Then, $\bigoplus_{\N}\Q=\underrightarrow{\lim}\, \Q^n$ is a piecewise hyperdefinable group. Now, $\bigoplus_\N \Q$ is not a topological group. Consider the set $U=\{x\sth |x_j|<|\cos(jx_0)|\mathrm{\ for\ }j\in\N_{>0}\}$. As $x_0\in\Q$ for any $x\in \bigoplus_{\N} \Q$, we have $\cos(jx_0)\neq 0$, so $U$ is an open neighbourhood of $0$. However, there is no open neighbourhood $V$ of $0$ such that $V+V\subseteq U$, concluding that $\bigoplus_{\N}\Q$ is not a topological group. Aiming a contradiction, suppose otherwise; take $V$ an open neighbourhood of $0$ such that $V+V\subseteq U$. As $V$ is an open neighbourhood of $0$, there is $\varepsilon_0\in\R_{>0}$ such that $\{x\sth |x_0|<\varepsilon_0 \mathrm{\ and\ }x_i=0\mathrm{\ for\ }i\in\N_{>0}\}\subseteq V$. Take $n\in\N_{>0}$ such that $2n\varepsilon_0>\pi$. There is then $\varepsilon_1\in\R_{>0}$ such that $\{x\sth |x_n|<\varepsilon_1\mathrm{\ and\ }x_i=0\mathrm{\ for\ }i\neq n\}\subseteq V$. Hence, $\{x\sth |x_0|<\varepsilon_0,\ |x_n|<\varepsilon_1\mathrm{\ and\ }x_i=0\mathrm{\ for\ }i\in\N\setminus \{0,n\}\}\subseteq V+V\subseteq U$. In particular, $(-\varepsilon_0,\varepsilon_0)_{\Q}\times(-\varepsilon_1,\varepsilon_1)_{\Q}\subseteq \{(x_0,x_1)\in \Q\times\Q\sth |x_1|<|\cos(nx_0)|\}$. However, this is impossible when $2n\varepsilon_0>\pi$, getting a contradiction.
\end{ex}
\begin{theo} \label{t:piecewise hyperdefinable and local compactness} 
Let $G$ be a locally hyperdefinable group with a global logic topology. Then, $G$ is a locally compact topological group with this topology.

Furthermore, a countably piecewise hyperdefinable group $G$ is a locally compact topological group with some logic topology if and only if this logic topology is the global logic topology and $G$ is locally hyperdefinable.
\begin{proof} We know that the  global logic topology of $G$ is locally compact by \cref{p:locally compact locally hyperdefinable} and Hausdorff by \cref{p:hausdorff locally hyperdefinable}. By \cref{p:product of locally hyperdefinable sets} and \cref{p:functions logic topology}(\ref{itm:functions logic topology:continuous}), we conclude that $G$ is a locally compact topological group.

By definition, a logic topology is $\mathrm{T}_1$ if and only if it is the  global logic topology. On the other hand, assuming $G=\underrightarrow{\lim}_{n\in\N}G_n$, by Baire's Category Theorem \cite[Theorem 48.2]{munkres1999topology}, if $G$ is locally compact Hausdorff, there are $h$ and $n\in \N$ such that $G_n$ is a neighbourhood of $h$. Thus, for any $g\in G$, $gh^{-1}G_n$ is an $\bigwedge${\hyp}definable neighbourhood of $g$. 
\end{proof}
\end{theo}
\begin{ex} \label{e:example no locally compact topological group} We give an example of a countably piecewise hyperdefinable group with a global logic topology which is not locally hyperdefinable; this is the infinite countable direct sum of circles with the inductive topology. Denote the unit circle by $\mathbb{S}\coloneqq\{ x\in \mathbb{C}\sth |x|=1\}$, which is hyperdefinable in the theory of real closed fields as quotient of the common definable circle by the infinitesimals. For $n>m$, we take the map $\psi_{n,m}:\ \mathbb{S}^m\rightarrow \mathbb{S}^n$ by $\psi_{n,m}(x)=(x,1,\ldots,1)$. Then, $\bigoplus_{\N} \mathbb{S}\coloneqq \underrightarrow{\lim}\, \mathbb{S}^n$ is a countably piecewise hyperdefinable group with a global logic topology that is not locally hyperdefinable. 

A local base of open neighbourhoods of the identity in the global logic topology of $\bigoplus_{\N} \mathbb{S}$ is the family of subsets $U_\varepsilon\coloneqq\{x\sth  \mathrm{d}_{\mathbb{S}}(x_i,1)<\varepsilon_i\mathrm{\ for\ }i\in\N\}$ for sequences $\varepsilon=(\varepsilon_i)_{i\in\N}$ with $\varepsilon_i\in (0,1]$, where $\mathrm{d}_{\mathbb{S}}$ is the normalised usual distance in the unit circle. Indeed, suppose $U$ is an open neighbourhood of $1$ in $\bigoplus_{\N} \mathbb{S}$. By using compactness in $\mathbb{S}^n$ for each $n\in\N$, we recursively find a sequence $\varepsilon=(\varepsilon_i)_{i\in\N}$ with $\varepsilon_i\in (0,1]$ such that $\{x\sth \mathrm{d}_{\mathbb{S}}(x_i,1)\leq \varepsilon_i\mathrm{\ for\ }i\leq n\}\subseteq U\cap \mathbb{S}^n$.
\end{ex}
Unfortunately, proving that a particular piecewise hyperdefinable set is locally hyperdefinable may be truly hard, as it requires to check a property about a topological space that we understand only vaguely. Until now, the only method available to show that a piecewise hyperdefinable set is locally hyperdefinable relies on \cref{p:mapping locally hyperdefinable sets} and the fact that piecewise definable sets are trivially locally hyperdefinable (i.e. piecewise definable and locally definable are the same). Sometimes that is not enough. To solve this problem we introduce generic sets.

Let $G$ be a piecewise hyperdefinable group. A \emph{generic} subset is an $\bigwedge${\hyp}definable subset $V$ such that, for any other $\bigwedge${\hyp}definable subset $W$, $[W:V]\coloneqq\min\{|\Delta|\sth W\subseteq \Delta V\}$ is finite, i.e. there is a finite $\Delta\subseteq G$ with $W\subseteq \Delta V$. Obviously, if $V$ is a generic set and $V\subseteq W$ for $\bigwedge${\hyp}definable $W$, then $W$ is also a generic set. In other words, if there are generic sets, the generic sets form an upper set of the family of $\bigwedge${\hyp}definable subsets. Hence, if there is a generic set, then there is in particular a generic piece, i.e. a piece which is a generic set. 

The following theorem is a generalisation of an unpublished example due to Hrushovski.
\begin{theo}[Generic Set Lemma]  \label{t:generic set lemma} 
Let $G$ be a piecewise hyperdefinable group and $V$ a symmetric generic set. Then, for each $n\in\N$, $V^{n+2}$ is a neighbourhood of $V^n$ in some logic topology. In particular, if $G$ has a generic piece, then $G$ is locally hyperdefinable. Furthermore, when $G$ is small, $G$ is locally hyperdefinable if and only if it has a generic piece.
\begin{proof}Using that $V$ is generic, find a well{\hyp}ordered sequence $(a_\xi)_{\xi\in \alpha}$ in $G$ such that, for every $\bigwedge${\hyp}definable subset $W\subseteq G$, there is $\Delta_W\subseteq \alpha$ finite with $W\subseteq \bigcup_{\xi\in \Delta_W}a_\xi V$. Let $A$ be a set of parameters containing $\{a_\xi\}_{\xi\in\alpha}$ and such that $G=\underrightarrow{\lim}\ G_i$ is a piecewise $A${\hyp}hyperdefinable group and $V$ is $\bigwedge_A${\hyp}definable. From now on, we work on the $A${\hyp}logic topology. We want to show that $V^{n+2}$ is a neighbourhood of $V^n$. Let $\Delta\subseteq \alpha$ be finite and minimal such that $V^{n+4}\subseteq V^{n+2}\cup \bigcup_{\xi\in \Delta}a_\xi V$. Let $U=V^{n+4}\setminus \bigcup_{\xi\in \Delta}a_\xi V$. Obviously, $U\subseteq V^{n+2}$ and $U$ is open in $V^{n+4}$. Note also that $V^n\subseteq U$; otherwise, taking $a\in V^n\setminus U$, there is $\xi\in\Delta$ such that $a\in a_{\xi}V$, so $a_{\xi}\in V^{n+1}$ and $a_{\xi}V\subseteq V^{n+2}$, contradicting minimality of $\Delta$. Similarly, for each piece $G_i$ such that $V^{n+4}\subseteq G_i$, pick a finite and minimal subset $\Delta_i\subseteq \alpha$ such that $G_i\subseteq V^{n+2}\cup \bigcup_{\xi\in \Delta_i}a_\xi V$ and $\Delta\subseteq \Delta_i$. Define $U_i=G_i\setminus \bigcup_{\xi\in\Delta_i}a_\xi V$. Again, it is clear that $U_i\subseteq V^{n+2}\subseteq V^{n+4}$ and $U_i$ is open in $G_i$. Also, by minimality of $\Delta_i$, it follows that $V^n\subseteq U_i$.

We claim that $U_i=U$ for any $i\in I$ such that $V^{n+4}\subseteq G_i$. It is clear by definition that $U_i\subseteq U$. On the other hand, take $a\in V^{n+2}\setminus U_i$. As $a\notin U_i$ and $a\in V^{n+2}\subseteq V^{n+4}\subseteq G_i$, there is $\xi\in\Delta_i$  such that $a\in a_\xi V$. As $a\in V^{n+2}$, $a_\xi\in a\cdot V\subseteq V^{n+3}$, concluding $a_\xi V\subseteq V^{n+4}$. Then, by minimality of $\Delta_i$, it follows that $\xi\in \Delta$, so $a\notin U$. That shows that $U$ is open in $G$, so $V^{n+2}$ is a neighbourhood of $V^n$. 

In particular, $V^3$ is a neighbourhood of $V$. As $G=\bigcup_{\xi\in\alpha} a_{\xi}V$, for every $a\in G$ there is $\xi\in\alpha$ such that $a\in a_\xi V\subseteq a_\xi U\subseteq a_\xi V^3$ with $a_\xi U$ open in the $A${\hyp}logic topology. Thus, $a_\xi V^3$ is an $\bigwedge_A${\hyp}definable neighbourhood of $a$ in the $A${\hyp}logic topology, concluding that $G$ is locally $A${\hyp}hyperdefinable. 

On the other hand, if $G$ is locally hyperdefinable and small, it is a locally compact topological group by \cref{t:piecewise hyperdefinable and local compactness}. Therefore, the identity is in the interior of some piece of $G$. Every $\bigwedge${\hyp}definable subset $W$ of $G$ is compact, and so covered by finitely many translates of this piece. As $W$ is arbitrary, this piece is generic. 
\end{proof}
\end{theo}
If $G$ has a generic subset, so has $\rfrac{G}{K}$ for $K\trianglelefteq G$ piecewise $\bigwedge${\hyp}definable. Therefore, by the Generic Set Lemma (\cref{t:generic set lemma}), when $G$ has a generic subset, $\rfrac{G}{K}$ is a locally hyperdefinable group for any piecewise $\bigwedge${\hyp}definable normal subgroup $K\trianglelefteq G$.\medskip

A \emph{$T${\hyp}rough $k${\hyp}approximate subgroup} of a group $G$ is a subset $X\subseteq G$ such that $X^2\subseteq \Delta XT$ with $|\Delta|\leq k\in\N_{>0}$ and $1\in T\subseteq G$. In particular, a \emph{$k${\hyp}approximate subgroup} is a $1${\hyp}rough $k${\hyp}approximate subgroup. 

It is clear from the definitions that every symmetric generic set is in particular an approximate subgroup. Conversely, any $\bigwedge${\hyp}definable approximate subgroup is a generic set of the piecewise hyperdefinable group that it generates. Thus, we can understand generic sets as a strengthening of $\bigwedge${\hyp}definable approximate subgroups.

The next corollary follows from \cref{t:piecewise hyperdefinable and local compactness} and the Generic Set Lemma (\cref{t:generic set lemma}) as a particular case:
\begin{coro} \label{c:corollary locally hyperdefinable and generic piece} Let $G=\underrightarrow{\lim}\, X^n$ be a piecewise hyperdefinable group generated by an $\bigwedge${\hyp}definable symmetric set $X$ and $T\trianglelefteq G$ be a normal piecewise $\bigwedge${\hyp}definable subgroup of small index. Then, $\rfrac{G}{T}$ is a locally compact topological group if and only if $X^n$ is a $T${\hyp}rough approximate subgroup for some $n$. In particular, if $T$ is $\bigwedge${\hyp}definable, $\rfrac{G}{T}$ is a locally compact topological group if and only if $X^n$ is an approximate subgroup for some $n$.
\end{coro}

An \emph{isomorphism} of piecewise $A${\hyp}hyperdefinable groups is an isomorphism of groups which is also an isomorphism of piecewise $A${\hyp}hyperdefinable sets.

\begin{theo}[Isomorphism Theorem] \label{t:isomorphism theorem} 
Let $f:\ G\rightarrow H$ be an onto piecewise bounded and proper $\bigwedge_A${\hyp}definable homomorphism of piecewise $A${\hyp}hyperdefinable groups. Then, for each $\bigwedge_A${\hyp}definable subgroup $S\leq K\coloneqq \ker(f)$ with $S\trianglelefteq G$, there is a unique map $\widetilde{f}_S:\ \rfrac{G}{S}\rightarrow H$ such that $f=\widetilde{f}_S\circ \pi_{G/S}$. This map $\widetilde{f}_S$ is a piecewise bounded and proper $\bigwedge_A${\hyp}definable homomorphism of groups with kernel $\rfrac{K}{S}$. In particular, $\widetilde{f}_K:\ \rfrac{G}{K}\rightarrow H$ is a piecewise $\bigwedge_A${\hyp}definable isomorphism. Furthermore, if $G$ has a global logic topology, then $f$ is an open map between the  global logic topologies. 
\begin{proof} Existence and uniqueness are given by the usual Isomorphism Theorem. Say $G=\underrightarrow{\lim}\, G_i$ and $H=\underrightarrow{\lim}\, H_j$. We get that $\widetilde{f}_S$ is obviously piecewise $\bigwedge_A${\hyp}definable as, for any pieces $G_i$ and $H_j$, $\quot^{-1}_{(\rfrac{G_i}{S})\times H_j}[\widetilde{f}_S]=\quot^{-1}_{G_i\times H_j}[f]$. It is trivially piecewise bounded and proper as $f$ and $\pi_{G/S}$ are so. Assuming that $G$ has a global logic topology, one{\hyp}side translations are continuous. Since, by \cref{p:functions logic topology}(\ref{itm:functions logic topology:homeomorphism}), $\widetilde{f}_K$ is a piecewise $\bigwedge_A${\hyp}definable homeomorphism, it is in particular an open map. Since $\pi^{-1}_{G/K} [\pi_{G/K}[U]]=KU=\bigcup_{x\in K} xU$, we see that $\pi_{G/K}$ is open. Therefore, $f=\widetilde{f}_K\circ\pi_{G/K}$ is open too.
\end{proof}
\end{theo}

\subsubsection*{Digression on Machado's Closed Approximate Subgroup Theorem:} 
The Generic Set Lemma (\cref{t:generic set lemma}) can be rewritten in purely topological terms. Written in this way, it is likely that this result was already (partially) known in the theory of topological groups. Recall that a semitopological group is a group with a $\mathrm{T}_1$ topology such that left and right translations are continuous.
\begin{theo} \label{t:generics}
Let $G$ be a semitopological group with a coherent covering $\mathcal{C}$ by closed symmetric subsets such that, for any $A,B\in\mathcal{C}$, there is $C\in\mathcal{C}$ with $AB\subseteq C$. Let $V$ be a \emph{generic} piece, i.e. an element $V\in\mathcal{C}$ such that $[W:V]$ is finite for every $W\in \mathcal{C}$. Then, $V$ has non{\hyp}empty interior. In particular, if $V$ is compact Hausdorff with the subspace topology, then $G$ is a locally compact topological group. Furthermore, if $\mathcal{C}$ is a countable covering by compact Hausdorff subsets, then $G$ is locally compact Hausdorff if and only if $\mathcal{C}$ has a generic piece. 
\begin{proof} Mimicking the proof of \cref{t:generic set lemma}, we show that $V^2$ is a neighbourhood of the identity, and so $V$ has non{\hyp}empty interior as finitely many translates of it cover $V^2$. If $V$ is compact Hausdorff, as left translations are homeomorphisms, we get that every point has a compact Hausdorff neighbourhood. Therefore, $G$ is locally compact Hausdorff, so a topological group by Ellis's Theorem \cite[Theorem 2]{ellis1957locally}. The ``furthermore'' part is an immediate consequence of Baire's Category Theorem \cite[Theorem 48.2]{munkres1999topology}. 
\end{proof}
\end{theo}
As a corollary we get the following notable Closed Approximate Subgroups Theorem, which was first proved by Machado in \cite[Theorem 1.4]{machado2023closed}. Recall that the \emph{commensurator} in $G$ of an approximate subgroup $\Lambda\subseteq G$ is the subset $\mathrm{Comm}(\Lambda)=\{g\in G\sth g^{-1}\Lambda g\mathrm{\ and\ }\Lambda$ are commensurable$\}$. The commensurator was first introduced by Hrushovski in \cite{hrushovski2022beyond}. The following is one of the most fundamental results about the commensurator:
\begin{lem}{\textnormal{{\cite[Lemma 5.1]{hrushovski2022beyond}}}} \label{l:commensurator} Let $G$ be a group and $\Lambda\subseteq G$ an approximate subgroup. Then, $\mathrm{Comm}(\Lambda)=\bigcup \{\Lambda'\subseteq G$ approximate $\mathrm{subgroup}\sth \Lambda'$ and $\Lambda$ are commensurable$\}$ and $\mathrm{Comm}(G)\leq G$.
\end{lem}
\begin{coro}[Machado's Closed Approximate Subgroups Theorem] \label{c:machado}
Let $G$ be a locally compact topological group and $X$ a closed approximate subgroup. Take $\Lambda=\overline{X^2\cap K^2}$ where $K$ is a compact symmetric neighbourhood of the identity. Then, $\mathrm{Comm}(\Lambda)$, with the direct limit topology given by $\Omega=\{\Lambda'\sth \Lambda'$ compact approximate subgroup commensurable to $\Lambda\}$, is a locally compact topological group such that $\inc:\ \mathrm{Comm}(\Lambda)\rightarrow G$ is a continuous $1${\hyp}to{\hyp}$1$ group homomorphism, $X\subseteq \mathrm{Comm}(\Lambda)$, $\inc_{\mid X}$ is a homeomorphism and $X$ has non{\hyp}empty interior. Furthermore, if $G$ is a Lie group, then $\mathrm{Comm}(\Lambda)$ is a Lie group too.
\begin{proof} As $\Lambda$ is compact, for any $B$ commensurable to $\Lambda$, we have that $\overline{B}$ is compact and commensurable to $\Lambda$. Thus, $\Omega$ is a covering of $\mathrm{Comm}(\Lambda)$ by \cref{l:commensurator}. By construction, $\Lambda$ is generic in $\Omega$. Also, by an easy computation, if $\Lambda_1$ and $\Lambda_2$ are commensurable approximate subgroups, then $\Lambda_1\Lambda_2\cup \Lambda_2\Lambda_1$ is an approximate subgroup commensurable to $\Lambda_1$ and to $\Lambda_2$. Thus, by \cref{t:generics}, $\mathrm{Comm}(\Lambda)$ with the direct limit topology is a locally compact topological group such that $\inc:\ \mathrm{Comm}(\Lambda)\rightarrow G$ is a continuous proper $1${\hyp}to{\hyp}$1$ group homomorphism and $\Lambda$ has non{\hyp}empty interior. As any two compact neighbourhoods of the identity are commensurable, by \cite[Lemma 2.2, Lemma 2.3]{machado2023closed}, we get that $X^2\cap (K')^2\subseteq \mathrm{Comm}(\Lambda)$ for any compact neighbourhood of the identity $K'$, concluding that $X\subseteq X^2\subseteq \mathrm{Comm}(\Lambda)$. Since $\Lambda\subseteq X^2$ has non{\hyp}empty interior, we get that $X^2$ has non{\hyp}empty interior and so $X$ has non{\hyp}empty interior. If $Y\subseteq X$ is closed in $G$, by continuity of $\inc$, it is closed in $\mathrm{Comm}(\Lambda)$. On the other hand, if $Y\subseteq X$ is closed in $\mathrm{Comm}(\Lambda)$, then $Y\cap (K')^2=Y\cap \overline{X^2\cap (K')^2}$ is closed in $\overline{X^2\cap (K')^2}$ (and so in $G$) for any compact symmetric neighbourhood of the identity $K'$ of $G$. By local compactness of $G$, the family of compact symmetric neighbourhoods of the identity forms a local covering of $G$. Thus, by the Local Covering \cref{l:local coherent}, $Y$ is closed in $G$. Since $X$ is closed, we conclude that the subspace topologies of $X$ in $G$ and in $\mathrm{Comm}(\Lambda)$ coincide. Finally, if $G$ is a Lie group, $\mathrm{Comm}(\Lambda)$ is a Lie group by \cite[Chapter III, \S 8.2, Corollary 1]{bourbaki1975lie}.   
\end{proof}
\end{coro}
\begin{rmk}{\textnormal{{\cite[Theorem 4.1]{machado2023closed}}}} If $X$ is compact, we do not need to assume that $G$ is locally compact; in that case, we can simply take $\Lambda=X$.
\end{rmk}

\subsection{Model{\hyp}theoretic components}
We define now some model{\hyp}theoretic components for piecewise hyperdefinable groups. Let $G$ be a piecewise $A${\hyp}hyperdefinable group.\smallskip\\ 
The \emph{invariant component of $G$ over $A$} is
\[G^{000}_A\coloneqq\bigcap\left\{H\leq G\sth H\mbox{ is }A\mbox{{\hyp}invariant with }[G:H]\mbox{ small}\right\}.\]
The \emph{infinitesimal component of $G$ over $A$} is
\[G^{00}_A\coloneqq\bigcap\left\{H\leq G\sth H\mbox{ is p/w. }{\textstyle{\bigwedge}}\mbox{{\hyp}def. with } G^{000}_A\leq H\right\};\]
The \emph{connected component of $G$ over $A$} is 
\[G^0_A\coloneqq\bigcap \left\{ H\leq G\sth H\ \mathrm{and}\ G\setminus H \mbox{ are p/w. }{\textstyle{\bigwedge}}\mbox{{\hyp}def. with }G^{000}_A\leq H\right\}.\]
Obviously, $G^{000}_A\leq G^{00}_A\leq G^0_A\leq G$. 
\begin{lem} \label{l:normal subgroup of small index} 
Let $G$ be a piecewise $A${\hyp}hyperdefinable group and $T\leq G$ be an $A${\hyp}invariant subgroup of small index. Then, there is a unique maximal normal subgroup $\widetilde{T}\trianglelefteq G$ contained in $T$. Furthermore, $\widetilde{T}$ is $A${\hyp}invariant and has small index. Moreover, $\widetilde{T}$ is piecewise $\bigwedge${\hyp}definable when $T$ is so.
\begin{proof} Take $\widetilde{T}=\bigcap_{i\in [G:T]} T^{g_i}$ with $\{g_i\}_{i\in [G:T]}$ set of representatives. 
\end{proof}
\end{lem}
Let $B$ be a set of parameters with $A\subseteq B$. The \emph{$B${\hyp}logic topology} in $\rfrac{G}{G^{000}_A}$ is the one such that $V\subseteq \rfrac{G}{G^{000}_A}$ is closed if and only if $\pi^{-1}_{G/G^{000}_A}[V]$ is piecewise $\bigwedge_B${\hyp}definable. The global logic topology in $\rfrac{G}{G^{000}_A}$ is the one such that $V\subseteq \rfrac{G}{G^{000}_A}$ is closed if and only if $\pi^{-1}_{G/G^{000}_A}[V]$ is piecewise $\bigwedge${\hyp}definable.
\begin{theo}
Let $G$ be a piecewise $A${\hyp}hyperdefinable group. Then:
\begin{enumerate}[label={\rm{(\arabic*)}}, wide]
\item $G^{000}_A$ is an $A${\hyp}invariant normal subgroup of $G$ and $[G:G^{000}_A]$ is small. In fact, $[V:G^{000}_A]\leq 2^{|A|+|\lang|}$ for any $\bigwedge${\hyp}definable subset $V\subseteq G$.
\item Let $B$ be a small set of parameters with $A\subseteq B$. The inversion map is continuous on $\rfrac{G}{G^{000}_A}$ with the $B${\hyp}logic topology.
\item The global logic topology on $\rfrac{G}{G^{000}_A}$ coincides with the $B${\hyp}logic topology for every small $B$ containing $A$ and a set of representatives of $\rfrac{G}{G^{000}_A}$. Every translation map is continuous on $\rfrac{G}{G^{000}_A}$ with the global logic topology.
\item $G^{00}_A$ is a piecewise $\bigwedge_A${\hyp}definable normal subgroup of $G$. Furthermore, $\rfrac{G^{00}_A}{G^{000}_A}$ is the closure of the identity in the global logic topology of $\rfrac{G}{G^{000}_A}$.
\item Let $\pi:\ \rfrac{G}{G^{000}_A}\rightarrow \rfrac{G}{G^{00}_A}$ be the natural quotient map given by $\pi_{G/G^{00}_A}=\pi\circ\pi_{G/G^{000}_A}$. Then, $\pi$ is a Kolmogorov map between the $B${\hyp}logic topologies for any $B$ small with $A\subseteq B$. In particular, it is the Kolmogorov quotient between the global logic topologies. If $\rfrac{G}{G^{00}_A}$ is a topological group, then the group operations in $\rfrac{G}{G^{000}_A}$ are continuous.
\item $G^0_A$ is a piecewise $\bigwedge_A${\hyp}definable normal subgroup of $G$.
\item $\rfrac{G}{G^{00}_A}$ is a locally compact topological group whenever $G$ has a generic piece modulo $G^{00}_A$. In that case, $\rfrac{G^0_A}{G^{00}_A}$ is the connected component of $\rfrac{G}{G^{00}_A}$ in the global logic topology and $\rfrac{G^0_A}{G^{000}_A}$ is the connected component of $\rfrac{G}{G^{000}_A}$ in the global logic topology.
\end{enumerate}
\begin{proof} 
\begin{enumerate}[label={\rm{(\arabic*)}}, wide]
\item[\hspace{-1.4em}\setcounter{enumi}{1}\theenumi] Trivially, $G^{000}_A$ is $A${\hyp}invariant and $[G:G^{000}_A]$ is small. Using \cref{l:normal subgroup of small index}, it is obvious that $G^{000}_A$ is a normal subgroup. Finally, recall that, for real elements, an $A^*${\hyp}invariant equivalence relation on a definable set with a small amount of equivalence classes has at most $2^{|A^*|+|\lang|}$ equivalence classes --- indeed, such a relation is coarser than having the same type over an elementary substructure containing $A^*$. The same holds for $A${\hyp}invariant equivalence relations on hyperdefinable sets. Consequently, we actually have $[V:G^{000}_A]\leq 2^{|A|+|\lang|}$ for any $\bigwedge${\hyp}definable $V$.
\item Trivial.
\item Take $B$ small containing $A$ and a set of representatives of $\rfrac{G}{G^{000}_A}$. We have that, for any $V\subseteq \rfrac{G}{G^{000}_A}$, $\pi^{-1}_{G/G^{000}_A}[V]$ is $B${\hyp}invariant. 
In particular, the global logic topology on $\rfrac{G}{G^{000}_A}$ is the $B${\hyp}logic topology, so it is a well{\hyp}defined topology. Finally, as $\pi_{G/G^{000}_A}$ is a homomorphism and every translation map is piecewise bounded $\bigwedge${\hyp}definable, we conclude, by \cref{p:functions logic topology}(\ref{itm:functions logic topology:continuous}), that every translation map is continuous in $\rfrac{G}{G^{000}_A}$ with the global logic topology.
\item Let $\widehat{V}$ be the closure of the identity on the global logic topology of $\rfrac{G}{G^{000}_A}$ and $V=\pi^{-1}_{G/G^{000}_A}[\widehat{V}]$. Obviously, $V\subseteq G^{00}_A$. On the other hand, as translations are continuous by point {\rm{(3)}}, it follows that $\widehat{V}^{-1}\widehat{V}\subseteq \widehat{V}$ and $\bar{g}\widehat{V}\bar{g}^{-1}\subseteq \widehat{V}$ for any $\bar{g}\in\rfrac{G}{G^{000}_A}$. Thus, $V$ is a normal subgroup, concluding $G^{00}_A=V$. In particular, $G^{00}_A$ is a piecewise $\bigwedge_A${\hyp}definable normal subgroup.
\item As translations are continuous, for any $\bar{g}\in\rfrac{G}{G^{000}_A}$, the closure of $\bar{g}$ in the global logic topology is $\bar{g}\cdot\rfrac{G^{00}_A}{G^{000}_A}$. Thus, $\pi$ is the Kolmogorov quotient between the global logic topologies. 
By \cref{l:kolmogorov maps}, $\Ima\, \pi$ is a lattice isomorphism between the global logic topologies with inverse $\Ima^{-1}\pi$. Now, as $\pi$ is $A${\hyp}invariant, it follows that $\Ima\, \pi$ is a lattice isomorphism between the $B${\hyp}logic topologies with inverse $\Ima^{-1}\pi$. Therefore, by \cref{l:kolmogorov maps}, we conclude that $\pi$ is a Kolmogorov map between the $B${\hyp}logic topologies.   

Now, suppose $\rfrac{G}{G^{00}_A}$ is a topological group. Take $U\subseteq \rfrac{G}{G^{000}_A}$ open and $\bar{g}_1\bar{g}_2\in U$. Then, $\pi(\bar{g}_1)\pi(\bar{g}_2)\in \pi[U]$. As $\pi[U]$ is open, there are $\widehat{U}_1$ and $\widehat{U}_2$ open in $\rfrac{G}{G^{000}_A}$ such that $\pi(\bar{g}_1)\in \widehat{U}_1$ and $\pi(\bar{g}_2)\in \widehat{U}_2$ and $\widehat{U}_1\widehat{U}_2\subseteq \pi[U]$. Write $U_1=\pi^{-1}[\widehat{U}_1]$ and $U_2=\pi^{-1}[\widehat{U}_2]$. By \cref{l:kolmogorov maps}, $\pi^{-1}[\pi[U]]=U$. Then, $U_1U_2\subseteq U$ with $\bar{g}_1\in U_1$, $\bar{g}_2\in U_2$ and $U_1$ and $U_2$ open. As $U$, $\bar{g}_1$ and $\bar{g}_2$ are arbitrary, we conclude that the product operation is continuous.
\item For any $H\leq G$ with $G^{000}_A\leq H$ such that $H$ and $G\setminus H$ are piecewise $\bigwedge${\hyp}definable, we have that $\rfrac{H}{G^{00}_A}$ is a clopen subgroup in the global logic topology. Hence, $\rfrac{G^0_A}{G^{00}_A}$ is the intersection of all the clopen subgroups of $\rfrac{G}{G^{00}_A}$ in the global logic topology. It follows that $\rfrac{G^0_A}{G^{00}_A}$ is closed in the global logic topology, so $G^0_A$ is piecewise $\bigwedge${\hyp}definable. As it is $A${\hyp}invariant, we conclude that $G^0_A$ is piecewise $\bigwedge_A${\hyp}definable. Finally, as every translation is a homeomorphism in $\rfrac{G}{G^{00}_A}$, we conclude that $\rfrac{G^0_A}{G^{00}_A}$ is a normal subgroup of $\rfrac{G}{G^{00}_A}$. Thus, $G^0_A\trianglelefteq G$.
\item By the Generic Set Lemma (\cref{t:generic set lemma}) and \cref{t:piecewise hyperdefinable and local compactness}, $\rfrac{G}{G^{00}_A}$ with the global logic topology is a locally compact topological group if it has a generic piece. In that case, by \cite[Theorem 2.1.4(b)]{dikranjan1990topological}, $\rfrac{G^0_A}{G^{00}_A}$ is the connected component of $\rfrac{G}{G^{00}_A}$. By point {\rm{(5)}}, $\rfrac{G^0_A}{G^{000}_A}$ is the connected component of $\rfrac{G}{G^{000}_A}$. \qedhere
\end{enumerate}
\end{proof}
\end{theo}
We define now a final special model{\hyp}theoretic component which has no analogues in the definable or hyperdefinable case. I want to thank Hrushovski for all his help via private conversations in relation with this result. Recall that an aperiodic topological group is a topological group that has no non{\hyp}trivial compact normal subgroups.

Let $G$ be a piecewise hyperdefinable group with a generic piece. The \emph{aperiodic component} $G^{\mathrm{ap}}$ of $G$ is the smallest piecewise $\bigwedge${\hyp}definable normal subgroup of $G$ with small index such that $\rfrac{G}{G^{\mathrm{ap}}}$ is an aperiodic locally compact topological group with its global logic topology.
\begin{lem} \label{l:lemma aperiodic} Let $G$ be an aperiodic locally compact topological group which is the union of $\lambda$ compact subsets with $\lambda\geq \aleph_0$. Then, $|G|\leq (\lambda+2^{\aleph_0})^{\lambda}$. 
\begin{proof} By Gleason{\hyp}Yamabe \cref{t:gleason yamabe,c:minimal yamabe pairs}, there is an open subgroup $H\leq G$ and a compact subgroup $K\trianglelefteq H$ such that $\rfrac{H}{K}$ is an aperiodic connected Lie group. Thus, $[H:K]\leq 2^{\aleph_0}$. On the other hand, as $G$ is a union of $\lambda$ compact subsets, $[G:H]\leq \lambda$, so $[G:K]\leq \lambda+2^{\aleph_0}$. Now, take $\{a_i\}_{i\in \lambda}$ a set of representatives of $\rfrac{G}{H}$ and write $K_i\coloneqq K^{a_i}=a_iKa_i^{-1}$ for each $i$. As $K\trianglelefteq H$, for any $g\in G$, we have $K^g=K_i$ for some $i\in \lambda$. Note that $[G:K_i]\leq \lambda+ 2^{\aleph_0}$, so $[G:\bigcap_{i\in \lambda} K_i]\leq (\lambda+ 2^{\aleph_0})^{\lambda}$. As $G$ is aperiodic, $\bigcap_{i\in \lambda} K_i=1$, so we conclude $|G|\leq (\lambda+ 2^{\aleph_0})^{\lambda}$. 
\end{proof}
\end{lem}
\begin{theo}  \label{t:the aperiodic component} Let $G$ be a piecewise $0${\hyp}hyperdefinable group with a generic piece. Then, $G^{\mathrm{ap}}$ exists, is piecewise $\bigwedge_0${\hyp}definable and does not change by expansions of the language.  
\begin{proof} Let $\lambda=\cf(G)$ and $\tau=(\lambda+2^{\aleph_0})^{\lambda}$. By assumption $\kappa>\lambda+2^{\aleph_0}$ is a strong limit cardinal, and so $\kappa>\tau^+$. For any small subset of parameters $A$, let $\mathcal{H}_A$ be the family of normal subgroups $H\trianglelefteq G$ of small index which are piecewise $\bigwedge_{\leq \tau}${\hyp}definable with parameters from $A$ (i.e. piecewise $\bigwedge_B${\hyp}definable for some $B\subseteq A$ with $|B|\leq \tau$) such that $\rfrac{G}{H}$, with its global logic topology, is an aperiodic locally compact topological group.

\begin{claim} For any $A$ small, $\mathcal{H}_A$ is closed under arbitrary intersections. Furthermore, for any $\mathcal{F}\subseteq \mathcal{H}_A$ there is $\mathcal{F}_0\subseteq \mathcal{F}$ with $|\mathcal{F}_0|\leq \tau$ and $\bigcap \mathcal{F}_0=\bigcap\mathcal{F}$.
\begin{proof} Take $\mathcal{F}\subseteq \mathcal{H}_A$. Then, $H\coloneqq \bigcap \mathcal{F}$ is a piecewise $\bigwedge_A${\hyp}definable normal subgroup of $G$ of small index. As $G$ has a generic piece, $\rfrac{G}{H}$ has a generic piece as well, so $\rfrac{G}{H}$ is a locally compact topological group with its global logic topology by the Generic Set Lemma (\cref{t:generic set lemma}) and \cref{t:piecewise hyperdefinable and local compactness}. Take $K\leq G$ such that $\rfrac{K}{H}\trianglelefteq \rfrac{G}{H}$ is a compact normal subgroup in its global logic topology. Take $F\in\mathcal{F}$ arbitrary. We know that $\pi:\ \rfrac{G}{H}\rightarrow\rfrac{G}{F}$ is a continuous homomorphism between the global logic topologies, so $\rfrac{K}{F}$ is a compact normal subgroup of $\rfrac{G}{F}$. As $\rfrac{G}{F}$ is aperiodic, we conclude that $K\leq F$. As $F\in\mathcal{F}$ is arbitrary, we conclude that $K\subseteq\bigcap \mathcal{F}=H$, so $\rfrac{G}{H}$ is aperiodic. Since $\rfrac{G}{H}$ is an aperiodic locally compact topological group, by \cref{l:lemma aperiodic}, $[G:H]\leq \tau$.

If $|\mathcal{F}|\leq \tau$, then $H$ is piecewise $\bigwedge_{\leq \tau}${\hyp}definable, so $H\in\mathcal{H}$. Now, we claim that there is $\mathcal{F}_0\subseteq \mathcal{F}$ with $|\mathcal{F}_0|\leq \tau$ such that $\bigcap \mathcal{F}_0=H$. Indeed, suppose otherwise, then we can find recursively a sequence $(F_i)_{i\in\tau^+}$ of elements in $\mathcal{F}$ such that $H<\bigcap_{i\in \alpha} F_i\subset \bigcap_{i\in \beta}F_i$ for $\beta<\alpha<\tau^+$. Thus, we can find a sequence $(g_i)_{i\in \tau^+}$ such that $g_i\in \rfrac{F_j}{H}$ for $j\leq i$ and $g_i\notin \rfrac{F_{i+1}}{H}$; contradicting that $[G:H]\leq \tau$. \claimqed
\end{proof}
\end{claim}

Take a $\tau^+${\hyp}saturated elementary substructure $\mathfrak{N}\preceq \mathfrak{M}$. Set $G^{\mathrm{ap}}=\bigcap \mathcal{H}_N$, so $G^{\mathrm{ap}}\in \mathcal{H}_N$. By $\tau^+${\hyp}saturation of $\mathfrak{N}$, we have that $G^{\mathrm{ap}}$ is $0${\hyp}invariant, so it is piecewise $\bigwedge_0${\hyp}definable. As $G^{\mathrm{ap}}$ is $0${\hyp}invariant and $\mathfrak{N}$ is $\tau^+${\hyp}saturated, we conclude that $G^{\mathrm{ap}}=\bigcap\mathcal{H}_A$ for any $A$. Therefore, $G^{\mathrm{ap}}$ is the smallest piecewise $\bigwedge_{\leq \tau}${\hyp}definable normal subgroup of $G$ with small index such that $\rfrac{G}{G^{\mathrm{ap}}}$ is an aperiodic locally compact topological group with its global logic topology.

Now, we show that $G^{\mathrm{ap}}$ does not change by expansions of the language. In any $\kappa${\hyp}saturated and strongly $\kappa${\hyp}homogeneous $\lang'/\lang${\hyp}expansion $\mathfrak{M}'$ of an elementary extension of $\mathfrak{M}$, we find a piecewise $\bigwedge_0${\hyp}definable normal subgroup $G^{\mathrm{ap}}_{\lang'}\trianglelefteq G$ which is the smallest $\bigwedge_{\leq \tau}${\hyp}definable normal subgroup of $G$ of small index such that $\rfrac{G}{G^{\mathrm{ap}}_{\lang'}}$ is an aperiodic locally compact topological group with its global logic topology. Since $[G:G^{\mathrm{ap}}_{\lang'}]\leq\tau$ by \cref{l:lemma aperiodic}, there is a large enough $\kappa${\hyp}saturated and strongly $\kappa${\hyp}homogeneous $\lang_0/\lang${\hyp}expansion $\mathfrak{M}_0$ of an elementary extension of $\mathfrak{M}$ such that, for any further $\kappa${\hyp}saturated and strongly $\kappa${\hyp}homogeneous $\lang_1/\lang_0${\hyp}expansion $\mathfrak{M}_1$ of an elementary extension of $\mathfrak{M}_0$, we have that $G^{\mathrm{ap}}_0(\mathfrak{M}_1)=G^{\mathrm{ap}}_1$ --- where $G^{\mathrm{ap}}_1$ is computed in $\mathfrak{M}_1$, $G^{\mathrm{ap}}_0$ is computed in $\mathfrak{M}_0$ and $G^{\mathrm{ap}}_0(\mathfrak{M}_1)$ is the realisation of $G^{\mathrm{ap}}_0$ in ${\mathfrak{M}_1}_{\mid \lang_0}$.

Take this base expansion $\mathfrak{M}_0$ of an elementary extension of $\mathfrak{M}$, with language $\lang_0$. By replacing $\mathfrak{M}_0$ by an elementary extension, we can assume that its $\lang${\hyp}reduct ${\mathfrak{M}_0}_{\mid \lang}$ is a $\kappa${\hyp}saturated and strongly $\kappa${\hyp}homogeneous elementary extension of $\mathfrak{M}$. Let $G^{\mathrm{ap}}_0$ be computed in $\mathfrak{M}_0$. Take $\alpha\in\Aut({\mathfrak{M}_0}_{\mid\lang})$. Consider the language $\lang_1$ expanding $\lang_0$ by adding a new symbol $\alpha x$ of the same sort that $x$ for each symbol $x$ of $\lang_0$ which is not in $\lang$. Take the $\lang_1${\hyp}expansion $\mathfrak{M}'_1$ of $\mathfrak{M}_0$ given by interpreting $(\alpha x)^{\mathfrak{M}'_1}=\alpha( x^{\mathfrak{M}_0})$. Let $\mathfrak{M}_1$ be a $\kappa${\hyp}saturated and strongly $\kappa${\hyp}homogeneous elementary extension of $\mathfrak{M}'_1$. Let $G^{\mathrm{ap}}_0(\mathfrak{M}_1)$ be the realisation of $G^{\mathrm{ap}}_0$ in ${\mathfrak{M}_1}_{\mid\lang_0}$, $G^{\mathrm{ap}}_1$ be computed in $\mathfrak{M}_1$ and $G^{\mathrm{ap}}_\alpha$ be computed in ${\mathfrak{M}_1}_{\mid \lang_1\setminus (\lang_0\setminus \lang)}$. By the choice of $\mathfrak{M}_0$, we know that $G^{\mathrm{ap}}_0(\mathfrak{M}_1)=G^{\mathrm{ap}}_1\subseteq G^{\mathrm{ap}}_0\cap G^{\mathrm{ap}}_\alpha$. Since $G^{\mathrm{ap}}_\alpha(\mathfrak{M}_0)$ is precisely $\alpha(G^{\mathrm{ap}}_0)$, we get that $G^{\mathrm{ap}}_0\subseteq \alpha(G^{\mathrm{ap}}_0)$. Therefore, as $\alpha$ is arbitrary, we conclude that $G^{\mathrm{ap}}_0$ is $\Aut({\mathfrak{M}_0}_{\mid\lang})${\hyp}invariant. By Beth's Definability Theorem \cite[Exercise 6.1.4]{tent2012course}, we conclude that $G^{\mathrm{ap}}_0$ is already piecewise $\bigwedge_0${\hyp}definable in $\lang$. Therefore, $G^{\mathrm{ap}}_0=G^{\mathrm{ap}}$ and $G^{\mathrm{ap}}$ does not change by expansions of the language.

Finally, by invariance under arbitrary expansions of the language, we have that $G^{\mathrm{ap}}$ is the smallest $\bigwedge${\hyp}definable normal subgroup of $G$ of small index such that $\rfrac{G}{G^{\mathrm{ap}}}$ is an aperiodic locally compact topological group with its global logic topology. 
\end{proof}
\end{theo}

\subsection{Lie cores}
Now, we adapt the results about Lie cores to piecewise hyperdefinable groups. \smallskip

Let $G$ be a piecewise $A${\hyp}hyperdefinable group. An \emph{$A${\hyp}Yamabe pair} of $G$ is a pair $(K,H)$ of subgroups $K\trianglelefteq H\leq G$ which is a Yamabe pair modulo $G^{00}_A$ for the global logic topology of $\rfrac{G}{G^{00}_A}$ and $\rfrac{K}{G^{00}_A}$ is $\bigwedge${\hyp}definable. In other words, it is a pair satisfying the following three properties:
\begin{enumerate}[label={\rm{(\roman*)}}, wide]
\item $H\leq G$ is a piecewise $\bigwedge${\hyp}definable subgroup whose complement is also piecewise $\bigwedge${\hyp}definable.
\item $K\trianglelefteq H$ is a piecewise $\bigwedge${\hyp}definable normal subgroup of $H$ such that $G^{00}_A\leq K$ and $\rfrac{K}{G^{00}_A}$ is $\bigwedge${\hyp}definable.
\item $L\coloneqq \rfrac{H}{K}$ with the respective global logic topology is a finite dimensional Lie group. 
\end{enumerate}

We say that $H$ is the \emph{domain}, $K$ is the \emph{kernel} and $L$ is the \emph{Lie core}. Write $\pi\coloneqq\pi_{H/K}:\ H\rightarrow L$ and $\widetilde{\pi}\coloneqq\widetilde{\pi}_{H/K}:\ \rfrac{H}{G^{00}_A}\rightarrow L$ for the quotient maps with $\pi=\widetilde{\pi}\circ\pi_{G/G^{00}_A}$ where $\pi_{G/G^{00}_A}:\ G\rightarrow \rfrac{G}{G^{00}_A}$. We say that a Lie group is an \emph{$A${\hyp}Lie core} of $G$ if it is isomorphic, as Lie group, to the Lie core of some $A${\hyp}Yamabe pair of $G$. 
\begin{rmk} \label{o:remarks lie core} Let $G$ be a piecewise $A${\hyp}hyperdefinable group. Let $(K,H)$ be an $A${\hyp}Yamabe pair of $G$ and $\pi\coloneqq \pi_{H/K}:\ H\rightarrow L$ its Lie core and $\widetilde{\pi}\coloneqq\widetilde{\pi}_{H/K}:\ \rfrac{H}{G^{00}_A}\rightarrow L$ such that $\pi=\widetilde{\pi}\circ\pi_{G/G^{00}_A}$.
\begin{enumerate}[label={\rm{(\arabic*)}}, ref={\rm{\arabic*}}, wide]
\item \label{itm:remarks lie core:modulo G^00} The fact that $G^{00}\leq K$ means that $[G:H]<\kappa$. Indeed, as $\rfrac{H}{K}$ is Hausdorff, we already have that $[H:K]<\kappa$. Therefore, if $[G:H]<\kappa$, we conclude that $[G:K]<\kappa$, so $G^{00}\leq K$ (for some set of parameters). Thus, the condition of working modulo $G^{00}$ is saying that $H$ is large in some sense. 

On the other hand, the condition that $\rfrac{K}{G^{00}_A}$ is $\bigwedge${\hyp}definable may seem superfluous. Indeed, as $(\rfrac{K}{G^{00}_A},\rfrac{H}{G^{00}_A})$ is a Yamabe pair, $\rfrac{K}{G^{00}_A}$ is compact. If $\rfrac{G}{G^{00}_A}$ has a generic piece (as we will assume in the rest of the section), this is enough to conclude that $\rfrac{K}{G^{00}_A}$ is $\bigwedge${\hyp}definable. However, in general, $\rfrac{G}{G^{00}_A}$ may have compact sets that are not $\bigwedge${\hyp}definable, so we need to add this condition. 
\item \label{itm:remarks lie core:quotient map} Note that $\widetilde{\pi}$ is a piecewise bounded and proper $\bigwedge${\hyp}definable map. Thus, by \cref{p:functions logic topology}, $\widetilde{\pi}$ is continuous and closed between the global logic topologies. By the Isomorphism \cref{t:isomorphism theorem}, since $\widetilde{\pi}$ is an onto homomorphism, we have that $\rfrac{(H/G^{00}_A)}{(K/G^{00}_A)}$ and $\rfrac{H}{K}$ are isomorphic as piecewise hyperdefinable groups and $\widetilde{\pi}$ is also an open map between the global logic topologies. Finally, as it has compact fibres, $\widetilde{\pi}$ is also proper \cite[Theorem 3.7.2]{engelking1989general}. In sum, $\widetilde{\pi}$ deeply connects the logic topology of $\rfrac{G}{G^{00}_A}$ and the geometry of $L$.

When $G^{00}_A$ is $\bigwedge_A${\hyp}definable, $\pi$ is also a piecewise bounded and proper $\bigwedge${\hyp}definable surjective homomorphism. Consequently, $\pi$ is a continuous, closed and proper map between the logic topologies (with enough parameters), concluding that the relation between $\rfrac{G}{G^{00}_A}$ and $L$ can be mostly lifted into a relation between $G$ and $L$. In that special case, we say that $\pi:\ H\rightarrow L$ is a \emph{Lie model} of $G$. In general, however, $G^{00}_A$ is only piecewise $\bigwedge_A${\hyp}definable and $\pi$ is only piecewise bounded --- we could even have $G=G^{00}_A$ (e.g. take $G=\bigcup_{n\in\N} [-a_n,a_n]$ with $a_n\in o(a_{n+1})$ for each $n\in\N$ in the theory of real closed fields). In Section 3, we give sufficient conditions to show that $G^{00}_A$ is $\bigwedge_A${\hyp}definable.
\item \label{itm:remarks lie core:minimal lie core} $L$ is a minimal Lie core if and only if it is an aperiodic connected Lie group. In that case, $\rfrac{K}{G^{00}_A}$ is the maximal compact normal subgroup of $\rfrac{H}{G^{00}_A}$, so it is its maximal $\bigwedge${\hyp}definable normal subgroup. In particular, every automorphism over $A$ leaving $H$ invariant leaves $K$ invariant.
\item \label{itm:remarks lie core:generic piece} As $L$ is locally compact and $\widetilde{\pi}$ is proper and continuous, we conclude that $\rfrac{H}{G^{00}_A}$ is locally compact. Since $\rfrac{H}{G^{00}_A}$ is open, we conclude that $\rfrac{G}{G^{00}_A}$ is locally compact with the global logic topology. Thus, if $G$ is a countably piecewise hyperdefinable group, by \cref{t:countably piecewise hyperdefinable groups and topological groups,t:piecewise hyperdefinable and local compactness} and the Generic Set Lemma (\cref{t:generic set lemma}), we conclude that $G$ has a generic piece modulo $G^{00}_A$. For that reason, we will assume from now on that $G$ has a generic piece modulo $G^{00}_A$.
\item \label{itm:remarks lie core:lie model} Note that our definitions of Lie core and Lie model extend the notion of Lie model used in \cite{hrushovski2011stable}. Suppose that $G$ is a piecewise $A${\hyp}definable group and $G^{00}_A$ is $\bigwedge_A${\hyp}definable. Then, $\pi$ is proper and continuous so, for any $\Gamma\subseteq U\subseteq L$ with $\Gamma$ compact and $U$ open, $\pi^{-1}[\Gamma]$ is $\bigwedge${\hyp}definable and $\pi^{-1}[U]$ is $\bigvee${\hyp}definable. Hence, there is a definable subset $D\subseteq G$ such that
\[\pi^{-1}[\Gamma]\subseteq D\subseteq \pi^{-1}[U].\]
As $L$ is first countable, it is metrisable by Birkhoff{\hyp}Kakutani Theorem \cite[Theorem 1.5.2]{tao2014hilbert}. Thus, every closed set in $L$ is $\mathrm{G_{\delta}}$ \cite[Example 2, page 249]{munkres1999topology}. In particular, it follows that the preimage of any compact set is $\bigwedge_{\omega}${\hyp}definable. When $L$ is connected, it is also second countable, so the preimage of every open set is $\bigvee_{\omega}${\hyp}definable.
\end{enumerate}
\end{rmk}
\begin{theo} \label{t:logic existence of lie core} Let $G$ be a piecewise $A${\hyp}hyperdefinable group with a generic piece modulo $G^{00}_A$. Then, $G$ has an $A${\hyp}Lie core. If $G$ is countably piecewise hyperdefinable, this condition is also necessary. 

Furthermore, for any $U$ such that $\rfrac{U}{G^{00}_A}$ is an open neighbourhood of the identity in the global logic topology of $\rfrac{G}{G^{00}_A}$, there is an $A${\hyp}Yamabe pair $(K,H)$ of $G$ with $K\subseteq UG^{00}_A$.
\begin{proof} ($\Rightarrow$) From Gleason{\hyp}Yamabe \cref{t:gleason yamabe}, using the Generic Set Lemma (\cref{t:generic set lemma}), and \cref{t:piecewise hyperdefinable and local compactness,p:compact locally hyperdefinable}. ($\Leftarrow$) We have explained the necessity of the generic piece assumption in \cref{o:remarks lie core}(\ref{itm:remarks lie core:generic piece}). 
\end{proof}
\end{theo}

\begin{theo} \label{t:logic existence of minimal lie core}
Let $G$ be a piecewise $A${\hyp}hyperdefinable group with a generic piece modulo $G^{00}_A$ and $(K_1,H_1)$ an $A${\hyp}Yamabe pair of $G$. Then, $G$ has a minimal $A${\hyp}Yamabe pair $(K,H)$ smaller than or equal to $(K_1,H_1)$. Furthermore, $H\subseteq U^2$ for any $U$ containing $K_1$ such that $\rfrac{U}{G^{00}_A}$ is a clopen neighbourhood of $\rfrac{K_1}{G^{00}_A}$ in the global logic topology. 
\begin{proof} Clear from \cref{c:minimal yamabe pairs}, using the Generic Set Lemma (\cref{t:generic set lemma}), and \cref{t:piecewise hyperdefinable and local compactness,p:compact locally hyperdefinable}. 
\end{proof}
\end{theo}
Applying also \cref{c:uniqueness minimal lie core}, we conclude the following result.
\begin{theo} \label{t:logic uniqueness of minimal lie core}
Let $G$ be a piecewise hyperdefinable group with a generic piece modulo $G^{00}_A$. Then, $G$ has a unique minimal $A${\hyp}Yamabe pair up to equivalence.
\end{theo}
By \cref{p:global minimal lie core map}, we get a \emph{global minimal Lie core map} $\widetilde{\pi}_L:\ \rfrac{D_L}{G^{00}_A}\rightarrow L$ extending all the minimal $A${\hyp}Yamabe pairs, which is unique up to isomorphisms of $L$. Let $\pi_L:\ D_L\rightarrow L$ be the map given by $\pi_L=\widetilde{\pi}_L\circ \pi_{G/G^{00}_A\ \mid D_L}$. Here, $D_L$ is the union of all the domains of minimal $A${\hyp}Yamabe pairs and $\ker(\pi_L)\coloneqq \pi_L^{-1}(1)$ is the union of all the kernels of minimal $A${\hyp}Yamabe pairs. Consequently, $D_L$ and $\ker(\pi_L)$ are $A${\hyp}invariant.
\begin{rmk} As noted in \cite{hrushovski2011stable}, the uniqueness of the minimal Lie core is achieved at a price. Indeed, while Gleason{\hyp}Yamabe \cref{t:gleason yamabe} gives us Yamabe pairs of arbitrarily small kernel, we have lost the control over the kernel in \cref{t:logic existence of minimal lie core}. If we do not care about uniqueness (as in \cite{breuillard2012structure} or \cite{massicot2015approximate}), it could be better just to apply \cref{t:logic existence of lie core} to find a Yamabe pair $(K,H)$ with arbitrarily small kernel. Also, it may be natural to apply Gleason{\hyp}Yamabe{\hyp}Kreitlon \cref{t:gleason yamabe carolino} rather than Gleason{\hyp}Yamabe \cref{t:gleason yamabe} to get some extra control on some parameters.
\end{rmk}
In \cref{p:minimal yamabe pairs with maximal domain}, we gave a criterion to find minimal Yamabe pairs with maximal domain in topological groups. Applying it modulo $G^{00}_A$, this result can be easily adapted to piecewise hyperdefinable groups.
Similarly, we can easily adapt \cref{p:minimal yamabe pair with minimal domain,p:asymptotic minimal yamabe pair with minimal domain} to the context of piecewise hyperdefinable groups by applying them modulo $G^{00}_A$. In particular, it follows that the minimal $A${\hyp}Lie core is piecewise $A${\hyp}hyperdefinable:
\begin{prop} \label{p:connected component and minimal lie core}
Let $G$ be a piecewise $A${\hyp}hyperdefinable group with a generic piece modulo $G^{00}_A$ and let $L$ be the minimal $A${\hyp}Lie core of $G$. Then, the restriction to $G^0_A$ of the global minimal Lie core map $\pi_{L\mid G^0_A}:\ G^0_A\rightarrow L$ is a piecewise bounded $\bigwedge${\hyp}definable surjective group homomorphism. Furthermore, we conclude that $L\cong \rfrac{G^0_A}{\ker(\pi_{L\mid G^0_A})}$. 
\end{prop}
The previous \cref{p:connected component and minimal lie core} gives us a canonical presentation of the minimal $A${\hyp}Lie core of $G$ as the piecewise $A${\hyp}hyperdefinable group $\rfrac{G^0_A}{K}$ with $K\coloneqq\ker(\pi_{L\mid G^0_A})$. Similarly, we get a canonical $A${\hyp}invariant presentation of the global minimal Lie core map $\pi_L:\ D_L\rightarrow L$ by taking $\pi_{L\mid G^0_A}=\pi_{G^0_A/K}$. We now give a more precise description of $K\coloneqq\ker(\pi_{L\mid G^0_A})$ using $G^{\mathrm{ap}}$.
\begin{lem} \label{l:minimal yamabe pair and aperiodic component} Let $G$ be a piecewise $0${\hyp}hyperdefinable group with a generic piece and $G^{\mathrm{ap}}$ the aperiodic component of $G$. Let $(K,H)$ be a minimal $A${\hyp}Yamabe pair of $G$ with Lie core $\pi_{H/K}:\ H\rightarrow L$. Then, $G^{\mathrm{ap}}\cap H\leq K$.  
\begin{proof}  Let $B$ be the small set of parameters with $A\subseteq B$ such that the $B${\hyp}logic topology of $\rfrac{G}{G^{00}_A}$ is its global logic topology. Denote by $\mathrm{cl}_B$ the closure in the $B${\hyp}logic topology of $G$. We define by recursion the sequence $(J_\alpha)_{\alpha\in\Ord}$ of piecewise $\bigwedge_B${\hyp}definable normal subgroups of $G$ given by $J_0=G^{00}_A$, $J_{\gamma}=\mathrm{cl}_B(\bigcup_{i\in \gamma} J_i)$ for $\gamma$ limit and $J_{\alpha+1}=\mathrm{cl}_B(\bigcup \mathcal{K}_\alpha)$ where $\mathcal{K}_\alpha$ is the family of piecewise $\bigwedge${\hyp}definable normal subgroups which are $\bigwedge${\hyp}definable modulo $J_\alpha$.

\begin{claim} For $\alpha_0\in\Ord$ large enough, $J_\alpha=G^{\mathrm{ap}}$ for any $\alpha\geq \alpha_0$.
\begin{proof} We prove inductively that $J_\alpha\subseteq G^{\mathrm{ap}}$ for any $\alpha\in\Ord$. Obviously, $J_0\subseteq G^{\mathrm{ap}}$. For $\gamma$ limit, assuming that $J_i\subseteq G^{\mathrm{ap}}$ for each $i\in \gamma$, as $G^{\mathrm{ap}}$ is piecewise $\bigwedge_B${\hyp}definable, we get that $J_\gamma\subseteq G^{\mathrm{ap}}$. Finally, assuming that $J_\alpha\subseteq G^{\mathrm{ap}}$, we have a piecewise bounded $\bigwedge${\hyp}definable onto homomorphism $\pi:\ \rfrac{G}{J_{\alpha}}\rightarrow \rfrac{G}{G^{\mathrm{ap}}}$. As $\rfrac{G}{G^{\mathrm{ap}}}$ is aperiodic, we have that every piecewise $\bigwedge${\hyp}definable normal subgroup of $G$ which is $\bigwedge${\hyp}definable modulo $J_{\alpha}$ must be contained in $G^{\mathrm{ap}}$. Therefore, $\bigcup \mathcal{K}_\alpha\subseteq G^{\mathrm{ap}}$, concluding that $J_{\alpha+1}\subseteq G^{\mathrm{ap}}$. 

On the other hand, as there is a small amount of piecewise $\bigwedge_B${\hyp}definable normal subgroups of $G$, for large $\alpha_0$, we must have $J_{\alpha_0+1}=J_{\alpha_0}$. In particular, that means that $\rfrac{G}{J_{\alpha_0}}$ is an aperiodic locally compact topological group with its global logic topology, so $G^{\mathrm{ap}}\leq J_{\alpha_0}$. Thus,  $G^{\mathrm{ap}}=J_\beta$ for $\beta\geq \alpha_0$. \claimqed
\end{proof}
\end{claim}

Now, we prove by induction that, for any $\alpha\in \Ord$, we have $J_\alpha\cap H\leq K$. Obviously, $J_0\cap H\leq K$. For $\gamma$ limit, assuming that $J_i\cap H\leq K$ for any $i\in \gamma$, we have $H\cap \bigcup_{i\in\gamma} J_i\leq K$. Therefore, $H\cap J_{\gamma}=\mathrm{cl}_B(H\cap \bigcup_{i\in\gamma}J_i)\subseteq K$, as $H$ is clopen and $K$ is closed in the $B${\hyp}logic topology. Finally, assuming $J_\alpha\cap H\leq K$, we have a piecewise bounded $\bigwedge${\hyp}definable onto homomorphism $\pi:\ \rfrac{H}{J_\alpha}\rightarrow L=\rfrac{H}{K}$. As $L$ is aperiodic, we have that $H\cap \bigcup\mathcal{K}_\alpha\subseteq K$. Therefore, $H\cap J_{\alpha+1}=\mathrm{cl}_B(H\cap \bigcup \mathcal{K}_\alpha)\subseteq K$, as $H$ is clopen and $K$ is closed in the $B${\hyp}logic topology. In particular, by the claim, we conclude that $G^{\mathrm{ap}}\cap H\leq K$.  
\end{proof}
\end{lem}

\begin{theo}  \label{t:canonical representation of the minimal lie core} Let $G$ be a piecewise $A${\hyp}hyperdefinable group with a generic piece. Let $L$ be the minimal $A${\hyp}Lie core of $G$ and $\pi_{L\mid G^0_A}:\ G^0_A\rightarrow L$ the restriction to $G^0_A$ of the global minimal Lie core map of $L$. Then, $G^0_A\cap G^{\mathrm{ap}}=\ker(\pi_{L\mid G^0_A})$. In particular, $L\cong\rfrac{G^0_A}{G^{\mathrm{ap}}}$. 
\begin{proof} By \cref{l:minimal yamabe pair and aperiodic component}, we know that $G^{\mathrm{ap}}\cap G^0_A\leq \ker(\pi_{L\mid G^0_A})$. On the other hand, $\ker(\pi_{L\mid G^0_A})=G^0_A\cap \pi_L^{-1}(1)\trianglelefteq G$ is $\bigwedge${\hyp}definable modulo $G^{00}_A\leq G^{\mathrm{ap}}$. Therefore, $\rfrac{\ker(\pi_{L\mid G^0_A})}{G^{\mathrm{ap}}}$ is a compact normal subgroup of $\rfrac{G}{G^{\mathrm{ap}}}$. As it is aperiodic, we conclude that $\ker(\pi_{L\mid G^0_A})\leq G^{\mathrm{ap}}$, so $\ker(\pi_{L\mid G^0_A})=G^{\mathrm{ap}}\cap G^0_A$. 
\end{proof}
\end{theo}
In sum, for any piecewise $A${\hyp}hyperdefinable group $G$ with a generic piece, we have the following structure in terms of the components $G^0_A$, $G^{00}_A$ and $G^{\mathrm{ap}}$:
\[\begin{array}{ccl}1.& \rfrac{G^{\mathrm{ap}}\cap G^0_A}{G^{00}_A}& \mathrm{is\ a\ compact\ topological\ group.}\\
2.&\rfrac{G^0_A}{G^{\mathrm{ap}}}&\mathrm{is\ an\ aperiodic\ connected\ Lie\ group.}\\
3.&\rfrac{G}{G^0_A}&\mathrm{is\ a\ totally\ disconnected\ locally\ compact}\\ & & \mathrm{topological\ group,\ i.e.\ a}\ locally\ profinite\ group.
\end{array}\]
Unfortunately, if $G^{00}_A=G^{\mathrm{ap}}=G^0_A=G$, all the previous results say nothing about $G$. In the following section, we extend the Stabilizer Theorem to the context of piecewise hyperdefinable groups. This theorem gives sufficient conditions to conclude that, with enough parameters, $G^{00}$ is $\bigwedge${\hyp}definable. As we pointed out at the beginning of the section, in this particular situation, the minimal Lie core gives very precise information about $G$ since the quotient homomorphism is piecewise bounded and proper.

\

We now note that the minimal Lie core is independent of the parameters. Furthermore, as $G^{\mathrm{ap}}$ is independent of expansions of the language, so is the minimal Lie core. 
\begin{coro} \label{c:independence of expansions of the minimal lie core} Let $G$ be a piecewise $0${\hyp}hyperdefinable group with a generic piece. Then, the minimal $A${\hyp}Lie core of $G$ is isomorphic to the minimal $0${\hyp}Lie core of $G$. Furthermore, the minimal Lie core of $G$ does not change by expansions of the language.
\begin{proof} By \cref{l:minimal yamabe pair and aperiodic component}, we have that the minimal $A${\hyp}Lie core of $G$ is isomorphic to the minimal Lie core of $\rfrac{G}{G^{\mathrm{ap}}}$. As the latter does not depend on parameters or expansions of the language, we conclude. 
%
\end{proof}
\end{coro}

The \emph{$A${\hyp}Lie rank} $\mathrm{Lrank}_A(G)$ of a piecewise hyperdefinable group $G$ with a generic piece modulo $G^{00}_A$ is the dimension of its minimal $A${\hyp}Lie core. As a consequence of \cref{t:logic uniqueness of minimal lie core}, $\mathrm{Lrank}_A(G)$ is a well{\hyp}defined invariant. Note that, by \cref{c:independence of expansions of the minimal lie core}, we have that $\mathrm{Lrank}(G)\coloneqq\mathrm{Lrank}_A(G)$ does not depend on the parameters when $G$ has a generic piece.
\begin{prop} Let $G$ be a piecewise $A${\hyp}hyperdefinable group with a generic piece modulo $G^{00}_A$ and $N\trianglelefteq G$ a piecewise $\bigwedge_A${\hyp}definable normal subgroup of small index. Then, 
\[\mathrm{Lrank}_A(G)\geq \mathrm{Lrank}_A(G/N)+\mathrm{Lrank}_A(N).\]
More precisely, let $L$ and $L_{G/N}$ be the minimal $A${\hyp}Lie cores of $G$ and $\rfrac{G}{N}$ respectively, presented with the canonical piecewise $A${\hyp}hyperdefinable structures given by \cref{p:connected component and minimal lie core}. Let $\pi_L$ and $\pi_{L_{G/N}}$ be their canonical global minimal Lie core maps. Let $(K,H)$ be a minimal $A${\hyp}Yamabe pair of $G$. Then:
\begin{enumerate}[label={\rm{(\arabic*)}}, wide]
\item $(K\cap N,H\cap N)$ is an $A${\hyp}Yamabe pair of $N$ with Lie core $\rfrac{H\cap N}{K\cap N}\cong L_N\coloneqq\pi_{L\mid H}[N]$. The connected component of $L_N$ is aperiodic so, in particular, $\mathrm{Lrank}_A(N)=\mathrm{dim}(L_N)$.
\item There is an $\bigwedge${\hyp}definable normal subgroup $T\trianglelefteq \rfrac{H}{N}$ with $\rfrac{K}{N}\trianglelefteq T$ such that $(T,\rfrac{H}{N})$ is a minimal $A${\hyp}Yamabe pair of $\rfrac{G}{N}$.
\item There is a piecewise bounded $\bigwedge_A${\hyp}definable surjective group homomorphism $\psi:\ L\rightarrow L_{G/N}$ such that $\psi\circ \pi_L=\pi_{L_{G/N}}\circ \pi_{G/N}$ (on the domain of $\pi_L$).
\end{enumerate}
\begin{proof} \begin{enumerate}[label={\rm{(\arabic*)}}, wide]
\item[\hspace{-1.4em}\setcounter{enumi}{1}\theenumi] Firstly, note that $G^{00}_A\leq N$, so we get $N^{00}_A=G^{00}_A$. We have that $\rfrac{H\cap N}{K\cap N}\cong \rfrac{H\cap N}{K}\cong L_N$. Note that $\rfrac{H\cap N}{G^{00}_A}$ is a closed normal subgroup of $\rfrac{H}{G^{00}_A}$. As $\widetilde{\pi}_{H/K}:\ \rfrac{H}{G^{00}_A}\rightarrow \rfrac{H}{K}$ is closed between the global logic topologies, we get that $L_N\cong \rfrac{H\cap N}{K}\cong \rfrac{((H\cap N)/G^{00}_A)}{(K/G^{00}_A)}$ is closed in $L\cong \rfrac{H}{K}\cong \rfrac{(H/G^{00}_A)}{(K/G^{00}_A)}$. Therefore, $L_N$ is a closed normal subgroup of $L$, concluding that $L_N$ is a Lie group. Then, $(K\cap N,H\cap N)$ is an $A${\hyp}Yamabe pair of $N$ with Lie core (isomorphic to) $L_N$. Let $L_N^0$ be the connected component of $L_N$. By \cref{l:maximal compact normal subgroup}, $L^0_N$ has a unique maximal compact normal subgroup $K_{N}$. Thus, by uniqueness, $K_N$ is characteristic in $L^0_N$, which is characteristic in $L_N$. Therefore, $K_N$ is a compact normal subgroup of $L$, so it is trivial by aperiodicity of $L$. In other words, $L^0_N$ is aperiodic. Thus, $L^0_N$ is the minimal $A${\hyp}Lie core of $N$. In particular, $\mathrm{Lrank}_A(N)=\mathrm{dim}(L_N)$. 
\item By point {\rm{(1)}}, we know that $L_N\trianglelefteq L$ is a closed subgroup, so $L_0\coloneqq\rfrac{L}{L_N}$ is a Lie group too. As $L$ is connected, $L_0$ is connected too. By \cref{l:maximal compact normal subgroup}, there is a maximal compact normal subgroup $T_0\trianglelefteq L_0$. Then, $\rfrac{L_0}{T_0}$ is a connected aperiodic Lie group. 

Since $\rfrac{G}{G^{00}_A}$ is a topological group, we know that the quotient homomorphism \[\pi_{(G/G^{00}_A)/(N/G^{00}_A)}:\ \rfrac{G}{G^{00}_A}\rightarrow  \rfrac{(G/G^{00}_A)}{(N/G^{00}_A)}\cong \rfrac{G}{N}\] is an open map. Therefore, $\rfrac{H}{N}$ is an open subgroup of $\rfrac{G}{N}$. Now, $\rfrac{G}{N}$ is small and contains a generic set, so $\rfrac{G}{N}$ is a locally compact topological group by the Generic Set Lemma (\cref{t:generic set lemma}) and \cref{t:piecewise hyperdefinable and local compactness}. Thus, $\rfrac{H}{N}$ is clopen. Consider, $\phi_0:\ \rfrac{H}{N}\rightarrow L_0$ given by $\phi_0\circ\pi_{(G/N)\mid H}=\pi_{L/L_N}\circ\pi_{L\mid H}$. It is clear that $\phi_0$ is a piecewise bounded $\bigwedge${\hyp}definable surjective group homomorphism with kernel $\rfrac{K}{N}$, which is $\bigwedge${\hyp}definable. Therefore, it is also piecewise proper. By the Isomorphism \cref{t:isomorphism theorem}, it follows that $(\rfrac{K}{N},\rfrac{H}{N})$ is an $A${\hyp}Yamabe pair of $\rfrac{G}{N}$ with Lie core (isomorphic to) $L_0$. 

As $L_0$ is locally hyperdefinable, $T_0$ is $\bigwedge${\hyp}definable. Thus, $\pi_{L_0/T_0}:\ L_0\rightarrow \rfrac{L_0}{T_0}$ is piecewise bounded and proper $\bigwedge${\hyp}definable surjective group homomorphism. Take $\phi=\pi_{L_0/T_0}\circ\phi_0:\ \rfrac{H}{N}\rightarrow \rfrac{L_0}{T_0}$. Then, $\phi$ is a piecewise bounded and proper $\bigwedge${\hyp}definable surjective group homomorphism. By the Isomorphism \cref{t:isomorphism theorem}, we conclude that $\rfrac{L_0}{T_0}\cong \rfrac{(H/N)}{T}$ where $T\coloneqq \ker(\phi)$ is an $\bigwedge${\hyp}definable normal subgroup of $\rfrac{H}{N}$ with $\rfrac{K}{N}\trianglelefteq T$. Consequently, $(T,\rfrac{H}{N})$ is a minimal $A${\hyp}Yamabe pair of $\rfrac{G}{N}$ with Lie core (isomorphic to) $\rfrac{L_0}{T_0}\cong L_{G/N}$.
\item By the Isomorphism \cref{t:isomorphism theorem}, take $\eta:\ \rfrac{L_0}{T_0}\rightarrow L_{G/N}$ isomorphism such that $\pi_{L_{G/N}\mid (\rfrac{H}{N})}=\eta\circ\phi$. Consider $\psi=\eta \circ \pi_{L_0/T_0}\circ\pi_{L/L_N}:\ L\rightarrow L_{G/N}$ --- note that, a priori, the definition of $\psi$ depends on $(K,H)$. Obviously, $\psi$ is a piecewise bounded $\bigwedge${\hyp}definable onto group homomorphism. Also, $\psi\circ \pi_{L\mid H}=\pi_{L_{G/N}\mid (\rfrac{H}{N})}\circ \pi_{(G/N)\mid H}$. 

Let $(K',H')$ be any other minimal $A${\hyp}Yamabe pair of $G$. For $h'\in H'$, there is $h\in H\cap H'$ such that $\pi_L(h')=\pi_L(h)$, i.e. $h^{-1}h'\in K'$. By point {\rm{(2)}}, it follows that $\pi_{L_{G/N}}(\pi_{G/N}(h))=\pi_{L_{G/N}}(\pi_{G/N}(h'))$. Thus, $\psi(\pi_L(h'))=\psi(\pi_L(h))=\pi_{L_{G/N}}(\pi_{G/N}(h))=\pi_{L_{G/N}}(\pi_{G/N}(h'))$. Therefore, $\psi\circ\pi_L=\pi_{L_{G/N}}\circ\pi_{G/N}$ on the domain of $\pi_L$.  

It remains to show that $\psi$ is $A${\hyp}invariant. Take $\sigma\in\Aut(\mathfrak{M}/A)$ and $x\in L$. Take $h$ such that $\pi_L(h)=x$. Using that $\pi_L$, $\pi_{L_{G/N}}$ and $\pi_{G/N}$ are $A${\hyp}invariant, we get that \[\psi(\sigma(x))=\psi\circ \pi_L(\sigma(h))=\pi_{L_{G/N}}\circ\pi_{G/N}(\sigma(x))=\sigma(\pi_{L_{G/N}}\circ\pi_{G/N}(x))=\sigma(\psi(h)).\]
\end{enumerate}
Finally, putting everything together, we conclude that 
\[\mathrm{Lrank}_A(G/N)=\mathrm{dim}(L_{G/N})\leq \mathrm{dim}(L_0)=\mathrm{Lrank}_A(G)-\mathrm{Lrank}_A(N).\]
\end{proof}
\end{prop}
The following proposition is a direct consequence of the results of \cite{jing2023nonabelian}. I am very grateful to Chieu{\hyp}Minh for telling me about it. 

\begin{prop}[An{\hyp}Jing{\hyp}Tran{\hyp}Zhang bound] Let $G$ be a piecewise hyperdefinable group and $X$ a symmetric generic set of $G$. Then, $\mathrm{Lrank}(G)\leq 12 \log_2(k)^2$ where $k=[X^2:X]$.    
\begin{proof} 
By working modulo $G^{00}$, we may assume that $G$ is small. Let $(K,H)$ be a minimal Yamabe pair of $G$ and $\pi\coloneqq \pi_{H/K}: H \rightarrow L$ the minimal Lie core. By the Generic Set Lemma (\cref{t:generic set lemma}), we have that $X^2$ is a symmetric compact neighbourhood of the identity in the global logic topology. Thus, as $H$ is open, $Y=H\cap X^2$ is also a symmetric compact neighbourhood of the identity in the global logic topology. By \cite[Lemma 2.3]{machado2023closed}, $Y$ is a $k^3${\hyp}approximate subgroup. As $\pi$ is a continuous and open homomorphism, $\pi[Y]$ is a compact neighbourhood of the identity and a $k^3${\hyp}approximate subgroup. As it is a neighbourhood of the identity, it has positive Haar measure, so the general Brunn{\hyp}Minkowski Inequality \cite[Theorem 1.1]{jing2023nonabelian} applies and we get that $2 \leq k^{\sfrac{3}{\alpha}}$ with $\alpha\coloneqq d-m-h$, where $d=\dim(L)$, $m=\max\{\dim(\Gamma)\sth \Gamma\leq L \mathrm{\ compact}\}$ and $h$ is the helix dimension of $L$. On the other hand, by \cite[Corollary 2.15]{jing2023nonabelian}, $h\leq \sfrac{n}{3}$ with $n\coloneqq d-m$, so $2n\leq 9\log_2(k)$. Finally, as $L$ has no compact normal subgroups, by \cite[Fact 3.6, Lemma 3.9]{an2021small}, we conclude that $d\leq \frac{n(n+1)}{2}$, so \[d\leq \frac{1}{2}\left\lfloor \frac{9}{2}\log_2(k)\right\rfloor^2+\frac{1}{2}\left\lfloor \frac{9}{2}\log_2(k)\right\rfloor\leq 12\log_2(k)^2.\]
\end{proof}
\end{prop}

\subsubsection{Digression on Breuillard{\hyp}Green{\hyp}Tao Classification} Let $\mathcal{A}$ be a particular class of piecewise hyperdefinable groups and $L$ a finite dimensional Lie group. It is then natural to study the \emph{direct problem} for $\mathcal{A}$ or (dually) the \emph{inverse problem} for $L$:  
\begin{center}
What are the Lie cores of elements of $\mathcal{A}$? 

Is $L$ a Lie core of an element of $\mathcal{A}$?
\end{center}

The classification theorem by Breuillard, Green and Tao for finite approximate subgroups \cite{breuillard2012structure} can be interpreted as an answer to the direct problem for the class of piecewise definable subgroups generated by pseudo{\hyp}finite definable approximate subgroups. Indeed, it can be restated to say that Lie cores of piecewise definable subgroups generated by pseudo{\hyp}finite definable approximate subgroups are nilpotent \cite[Proposition 9.6]{breuillard2012structure}. We consider it interesting to discuss the direct and inverse problems in general. In fact, one may expect that solutions of these questions would yield classification results similar to the one of \cite{breuillard2012structure}.\smallskip

The following easy examples show that the inverse problem is trivially solved for some basic classes of piecewise hyperdefinable groups. In particular, these examples show that no general classification result for the Lie cores analogous to the one of \cite{breuillard2012structure} could be found in those cases. 

\begin{ex} Any Lie group $L$ is a Lie core of some piecewise hyperdefinable group in the theory of real closed fields. Indeed, $L$ is in particular a second countable manifold, so it is a locally  hyperdefinable subset of countable cofinality. By \cref{p:functions logic topology 2}, the group operations are piecewise bounded $\bigwedge${\hyp}definable, so $L$ is a locally hyperdefinable group. Clearly, $L$ is its own Lie core. \smallskip

For a slightly more explicit construction, note that any connected Lie group $L$ is the Lie core of some piecewise definable group generated by a definable approximate subgroup. Indeed, connected Lie groups are metrisable  by Birkhoff{\hyp}Kakutani Theorem \cite[Theorem 1.5.2]{tao2014hilbert}, so take a left invariant metric $d$ for $L$. As $L$ is locally compact, we may assume that the closed unit ball $\overline{\mathbb{D}}$ is a compact symmetric neighbourhood of the identity. Consider the structure of $L$ with the language of groups, a sort for $\R$ with the language of ordered rings and a function symbol for the metric $d$. Let $L'$ be an $|L|${\hyp}saturated elementary extension of it. Consider the subgroup $H\leq L'$ generated by the closed unit ball and $E=\{(a,b)\sth d(a,b)<\sfrac{1}{n}\mathrm{\ for\ any\ }n\in \N\}$. Clearly, $\overline{\mathbb{D}}$ is a definable approximate subgroup and $H$ is the piecewise definable group generated by it. Then, $L$ is a Lie core of $H$, as, in fact, we have $\rfrac{H}{E}\cong L$. \smallskip

In the case of a linear connected Lie group $L$, we can combine both examples. Indeed, in that case, in the real numbers, $L$ is piecewise definable and its metric and group operations are definable, so, after saturation, we just need to take the piecewise definable group generated by the closed unit ball and quotient out by the infinitesimals as in the previous example.
\end{ex}

\section{Stabilizer Theorem}
In this section, we aim to extend the Stabilizer Theorem \cite[Theorem 3.5]{hrushovski2011stable} to piecewise hyperdefinable groups. The main point of this theorem is that it provides sufficient conditions to conclude that $G^{00}$ is $\bigwedge${\hyp}definable. As we already noted in the previous section, this result has significant consequences for the projection $\pi_L$ of the minimal Lie core.

To prove the Stabilizer Theorem, we need first to extend the model{\hyp}theoretic notions of dividing, forking and stable relation to piecewise hyperdefinable sets. Once these model{\hyp}theoretic notions have been properly defined for piecewise hyperdefinable groups, adapting the original proof of the Stabilizer Theorem is straightforward.

Forking and dividing for hyperimaginaries have already been well studied by many authors (e.g. \cite{hart2000coordinatisation,wagner2010simple,kim2014simplicity}). Here, in the first subsection, we rewrite the definitions of dividing and forking for hyperdefinable sets in a slightly different way which we find more natural from the point of view of this paper. After that, it is trivial to extend the definition to the context of piecewise hyperdefinable sets.
\subsection{Dividing and forking}
Let $(P_t)_{t\in T}$ be $A${\hyp}hyperdefinable sets. Write $P\coloneqq \prod P_t$. For infinite tuples $\overline{a},\overline{b}\in P$ of hyperimaginaries, we write $\tp(\overline{a}/A)=\tp(\overline{b}/A)$ to mean that $\tp(a_{\mid T_0}/A)=\tp(b_{\mid T_0}/A)$ for any $T_0\subseteq T$ finite.

Let $(I,<)$ be a linear order. An \emph{$A${\hyp}indiscernible} sequence in $P$ indexed by $(I,<)$ is a sequence $(\overline{a}_i)_{i\in I}$ of hyperimaginary tuples $\overline{a}_i\in P$ such that, for any $n\in \N$ and any $i_1<\cdots<i_n$ and $j_1<\cdots<j_n$,
\[\tp(\overline{a}_{i_1},\ldots,\overline{a}_{i_n}/A)=\tp(\overline{a}_{j_1},\ldots,\overline{a}_{j_n}/A).\]

\begin{lem}[Standard Lemma] \label{l:standard lemma for hyperimaginaries} Let $P=\prod_{t\in T}P_t$ be a product of $A${\hyp}hyperdefinable sets and $(I,<)$ and $(J,<)$ two infinite linear orders, with $|T|^+ +|J|\leq\kappa$. Then, given a sequence $(\overline{a}_i)_{i\in I}$ of hyperimaginary tuples in $P$, for any set of representatives $A^*$, there is a sequence $(\overline{b}_j)_{j\in J}$ of hyperimaginary tuples in $P$ with an $A^*${\hyp}indiscernible sequence of representatives $(\overline{b}^*_j)_{j\in J}$ such that
\[({\overline{b}_{j_1}}_{\mid T_0},\ldots,{\overline{b}_{j_n}}_{\mid T_0})\in W\]
for any $n\in \N$, any $j_1<\ldots<j_n$, any $T_0\subseteq T$ finite and any $\bigwedge_{A}${\hyp}definable set $W\subseteq \left(\prod_{t\in T_0} P_t\right)^n$ such that $({\overline{a}_{i_1}}_{\mid T_0},\ldots,{\overline{a}_{i_n}}_{\mid T_0})\in W$ for every $i_1<\cdots<i_n$.
\begin{proof} Trivial from the classic Ehrenfeucht{\hyp}Mostowski Standard Lemma proved using Ramsey's Theorem \cite[Theorem 5.1.5]{tent2012course}. 
\end{proof}
\end{lem}
\begin{coro} \label{c:indiscernible hyperimaginaries} 
Let $P=\prod_{t\in T}P_t$ be a product of $A${\hyp}hyperdefinable sets and $(I,<)$ an infinite linear order, with $|T|^+ +|I|\leq\kappa$. Then, a sequence $(\overline{b}_i)_{i\in I}$ of hyperimaginary tuples in $P$ is $A${\hyp}indiscernible if and only if, for some set of representatives $A^*$ of $A$, there is an $A^*${\hyp}indiscernible sequence of representatives $(\overline{b}^*_i)_{i\in I}$.

Furthermore, if $(\overline{b}_i)_{i\in I}$ is $A${\hyp}indiscernible, for any set of representatives $A^*$ there is an $A^{**}${\hyp}indiscernible sequence of representatives $(\overline{b}^*_i)_{i\in I}$ of $(\overline{b}_i)_{i\in I}$ where $A^{**}$ is another set of representatives of $A$ with $\tp(A^*)=\tp(A^{**})$. 
\begin{proof} By the Standard \cref{l:standard lemma for hyperimaginaries,c:types over hyperimaginaries}. 
\end{proof}
\end{coro}
Let $P=\rfrac{X}{E}$ be an $A${\hyp}hyperdefinable set and $\quot_P:\ X\rightarrow P$ its quotient map. An $\bigwedge${\hyp}definable subset $V\subseteq P$ \emph{divides over $A$} if and only if $\quot^{-1}_P[V]$ divides over $A^*$ for some set of representatives of $A$. 
\begin{lem} \label{l:finite character of dividing} 
Let $P$ be an $A${\hyp}hyperdefinable set and $V\subseteq P$ an $\bigwedge_{A,B}${\hyp}definable subset. Then, $V$ divides over $A$ if and only if there are a finite tuple $b$ from $B$ and an $\bigwedge_{A,b}${\hyp}definable subset $W\subseteq P$ such that $V\subseteq W$ and $W$ divides over $A$.
\begin{proof} One direction is trivial. Let us check the other. Take a uniform definition $\underline{V}$ of $V$, which exists by \cref{l:uniform definition}. Take representatives $A^*$ of $A$ such that $\underline{V}(x,A^*,B^*)$ divides over $A^*$. There is then $b$ finite such that $\underline{W}(x,A^*,b^*)\coloneqq \underline{V}(x,A^*,B^*)\cap \For^x(\lang(A^*,b^*))$ divides over $A^*$. Now, as $\underline{V}$ is a uniform definition, $W$ is $\bigwedge_{A,b}${\hyp}definable. Obviously, $V\subseteq W$ and $W$ divides over $A$. 
\end{proof}
\end{lem} 
\begin{lem} \label{l:dividing} 
Let $P$ and $Q$ be $A${\hyp}hyperdefinable sets. Let $V\subseteq Q\times P$ be an $\bigwedge_{A}${\hyp}definable set and $b\in Q$. Then, $V(b)$ divides over $A$ if and only if there is an $A${\hyp}indiscernible sequence $(b_i)_{i\in \omega}$ of hyperimaginaries from $Q$ such that $\tp(b_0/A)=\tp(b/A)$ and $\bigcap_{i\in \omega} V(b_i)=\emptyset$.

Furthermore, $V(b)$ divides over $A$ if and only if for any set of representatives $A^*$ of $A$ there is another set of representatives $A^{**}$ such that $\tp(A^*)=\tp(A^{**})$ and $V(b)$ divides over $A^{**}$.
\begin{proof} Assume $V(b)$ divides over $A$. Take a uniform definition $\underline{V}$ of $V$, which exists by \cref{l:uniform definition}. Take representatives $A^*$ of $A$ such that $\underline{V}(x,A^*,b^*)$ divides over $A^*$. There is then an $A^*${\hyp}indiscernible sequence $(b^*_i)_{i\in \omega}$ such that $\tp(b^*_0/A^*)=\tp(b^*/A^*)$ and $\bigcup_{i\in\omega}\underline{V}(x,A^*,b^*_i)$ is not finitely satisfiable. Then, $(b_i)_{i\in \omega}$ given by $b_i=\quot_Q(b^*_i)$ is $A${\hyp}indiscernible and $\tp(b_0/A)=\tp(b/A)$. Also, $\quot^{-1}_P\left[\bigcap_{i\in \omega} V(b_i)\right]=\bigcap_{i\in \omega} \quot^{-1}_P[V(b_i)]=\bigcap_{i\in\omega}\underline{V}(\mathfrak{M},A^*,b^*_i)=\emptyset$, so $\bigcap_{i\in\omega} V(b_i)=\emptyset$ by surjectivity of $\quot_P$.

On the other hand, assume there is an $A${\hyp}indiscernible sequence $(b_i)_{i\in \omega}$ in $Q$ such that $\tp(b_0/A)=\tp(b/A)$ and $\bigcap_{i\in \omega} V(b_i)=\emptyset$. Take  a uniform definition $\underline{V}$ of $V$, which exists by \cref{l:uniform definition}. By \cref{c:indiscernible hyperimaginaries}, there is an $A^{*}${\hyp}indiscernible sequence $(b^*_i)_{i\in \omega}$ of representatives of $(b_i)_{i\in \omega}$ with $A^{*}$ representatives of $A$. Now, as $\tp(b_0/A)=\tp(b/A)$, by \cref{c:types over hyperimaginaries}, there is a representative $b^{**}\in b$ and a set of representatives $A^{**}$ of $A$ such that $\tp(b^*_0, A^{*})=\tp(b^{**},A^{**})$. Take $\sigma\in\Aut(\mathfrak{M})$ mapping $(b_0^*,A^*)$ to $(b^{**},A^{**})$, and write $b'_i\coloneqq \sigma(b_i)$ and ${b'}^*_i\coloneqq \sigma(b^*_i)$ for $i\in\omega$. Then, $({b'}^*_i)_{i\in\omega}$ is an $A^{**}${\hyp}indiscernible sequence with ${b'}^*_0=b^{**}$. For this sequence, it follows that
$\emptyset=\quot^{-1}_P\left[\bigcap_{i\in \omega}V(b'_i)\right]=\bigcap_{i\in\omega}\quot^{-1}_P[V(b'_i)]=\bigcap_{i\in \omega}\underline{V}(\mathfrak{M},A^{**},{b'}^{*}_i),$ 
concluding that $\underline{V}$ divides over $A^{**}$. In particular, $V$ divides over $A$.

For the ``furthermore part'', note that in the previous paragraph $\tp(A^*)$ can be chosen arbitrarily by \cref{c:indiscernible hyperimaginaries}. 
\end{proof}
\end{lem}
Combining both propositions we conclude that the definition given here is equivalent to the one studied previously by other authors (e.g. \cite{hart2000coordinatisation,wagner2010simple,kim2014simplicity}). It is now straightforward to prove the following basic proposition.
\begin{lem} \label{l:dividing and functions} 
Let $P$ and $Q$ be $A${\hyp}hyperdefinable sets. Let $V\subseteq P$ be $\bigwedge${\hyp}definable and $f:\ P\rightarrow Q$ a $1${\hyp}to{\hyp}$1$ $\bigwedge_A${\hyp}definable function. Then, $f[V]$ divides over $A$ provided that $V$ divides over $A$.
\end{lem}
Let $P$ be an $A${\hyp}hyperdefinable set. The \emph{forking ideal} $\mathfrak{f}_P(A)$ of $P$ over $A$ is the ideal of $\bigwedge${\hyp}definable subsets of $P$ generated by the ones dividing over $A$. An $\bigwedge${\hyp}definable subset of $P$ \emph{forks over $A$} if it is in that forking ideal. 
\begin{rmk} Trivially, we have that $V$ forks over $A$ implies that $\underline{V}$ forks over $A$. In the case of simple theories, the converse also holds as forking and dividing are the same \cite[Proposition 3.2.7]{wagner2010simple}. In general, however, the converse is not true. For instance, let $\mathfrak{M}$ be the monster model of the theory of dense circular orders in the usual language, let $M$ be the whole $1${\hyp}ary universe of $\mathfrak{M}$ and $E$ the trivial equivalence relation given by $xEy$ for every $x,y\in M$. Obviously, $\rfrac{M}{E}$ is a singleton, so it does not fork over $\emptyset$. However, $M$ forks over $\emptyset$.
\end{rmk}
Let $P$ be a piecewise $A${\hyp}hyperdefinable set. An $\bigwedge${\hyp}definable set $V\subseteq P$ \emph{divides} over $A$ if it divides over $A$ as subset of some/any piece $P_i$ containing $V$. Note that, by \cref{l:dividing and functions}, this is well{\hyp}defined. An $\bigwedge${\hyp}definable subset $V$ of $P$ forks over $A$ if and only if $V$ forks over $A$ as subset of some/any piece $P_i$ containing $V$. The \emph{forking ideal} $\mathfrak{f}_{P}(A)$ of $P$ over $A$ is the family of $\bigwedge${\hyp}definable subsets of $P$ forking over $A$. Clearly, $\mathfrak{f}_P(A)$ is the ideal of $\bigwedge${\hyp}definable subsets of $P$ generated by the ones dividing over $A$.


\subsection{Ideals}
From now on, $\lambda$ is a cardinal with $\kappa\geq\lambda>|\lang|+|A|$.\medskip

Let $P$ be a piecewise $A${\hyp}hyperdefinable set and $\mu$ an ideal of $\bigwedge_{<\lambda}${\hyp}definable subsets in $P$. We say that $\mu$ is \emph{$A${\hyp}invariant} if it is invariant under $\Aut(\mathfrak{M}/A)$, i.e. $W(\overline{b})\in\mu$ implies $W(\overline{b}')\in\mu$ for any $\overline{b}'$ with $\tp(\overline{b}/A)=\tp(\overline{b}'/A)$ and any $\bigwedge_{\bar{b}}${\hyp}definable subset $W(\overline{b})\subseteq P$ with $|\overline{b}|<\lambda$. We say that an $\bigwedge_{<\lambda}${\hyp}definable subset is \emph{$\mu${\hyp}negligible} if it is in $\mu$ and that it is \emph{$\mu${\hyp}wide} if it is not in $\mu$. We say that $\mu$ is \emph{locally atomic} if, for any wide $\bigwedge_B${\hyp}definable subset $V$ with $|B|<\lambda$, there is $a\in V$ such that $\tp(a/B)$ is wide. 
\begin{rmk} For $X\subseteq P$, write $\mu_{\mid X}\coloneqq\{W\in \mu\sth W\subseteq X\}$. Clearly, $\mu_{\mid X}$ is an ideal of $\bigwedge_{<\lambda}${\hyp}definable subsets of $X$. Say $P=\underrightarrow{\lim}\, P_i$, then write $\mu_i\coloneqq \mu_{\mid P_i}$ for each piece. We have from the definitions that $\mu$ is $A${\hyp}invariant if and only if each $\mu_i$ is so. Similarly, it is locally atomic if and only if each $\mu_i$ is so.
\end{rmk}

We say that an $\bigwedge_A${\hyp}definable subset $V$ is \emph{$A${\hyp}medium} if for any $\bigwedge_{\bar{b}}${\hyp}definable subset $W(\overline{b})\subseteq V$ with $|\bar{b}|<\lambda$, we have 
\[W(\overline{b})\in \mu\Leftrightarrow W(\overline{b}_0)\cap W(\overline{b}_1)\in\mu\]
for any $A${\hyp}indiscernible sequence $(\overline{b}_i)_{i\in\omega}$ realising $\tp(\overline{b}/A)$. We say that an $\bigwedge_{<\lambda}${\hyp}definable subset $V$ is $A${\hyp}medium if there is an $A${\hyp}medium $\bigwedge_A${\hyp}definable subset $V_0$ with $V\subseteq V_0$. We say that $V$ is \emph{strictly $A${\hyp}medium} if it is wide and $A${\hyp}medium. We say that $\mu$ is \emph{$A${\hyp}medium\footnote{We use the terminology of \cite{montenegro2018stabilizers}. In \cite{hrushovski2011stable}, it is said that the ideal has the $S_1$ property.}} if every $\bigwedge_{<\lambda}${\hyp}definable subset is $A${\hyp}medium for $\mu$. 
\begin{rmk} Note that, by definition, if $\mu$ is $A${\hyp}medium, $\mu$ is in particular $A${\hyp}invariant.
Note that the family of $A${\hyp}medium sets of $\mu$ is an $A${\hyp}invariant ideal of $\bigwedge_{<\lambda}${\hyp}definable subsets.
\end{rmk}
\begin{ex} 
\begin{enumerate}[label={\rm{(\arabic*)}}, wide] 
\item[\hspace{-1.4em}\setcounter{enumi}{1}\theenumi] The forking ideal over $A$ is an $A${\hyp}invariant ideal of $\bigwedge${\hyp}definable subsets. In simple theories, the forking ideal over $A$ is locally atomic and $A${\hyp}medium. Indeed, it is locally atomic by the extension property of forking \cite[Theorem 3.2.8]{wagner2010simple}. On the other hand, suppose $V(\overline{b})$ is $\bigwedge_{\overline{b}}${\hyp}definable and $V(\overline{b}_0)\cap V(\overline{b}_1)$ forks over $A$ with $(\overline{b}_i)_{i\in\omega}$ an $A${\hyp}indiscernible sequence realising $\tp(\overline{b}/A)$. Then, by simplicity, $V(\overline{b}_0)\cap V(\overline{b}_1)$ divides over $A$, so $\bigcap^n_{i=0} V(\overline{b}_{2i})\cap V(\overline{b}_{2i+1})=\emptyset$ for some $n\in\N$. Thus, $V(\overline{b})$ divides over $A$, so it forks over $A$. As $V$ is arbitrary, we conclude that the forking ideal over $A$ is $A${\hyp}medium.
\item If we have an $A${\hyp}invariant measure on the lattice of $\bigwedge_{<\lambda}${\hyp}definable sets, the ideal of zero measure $\bigwedge_{<\lambda}${\hyp}definable subsets is an $A${\hyp}invariant ideal. In this case, every $\bigwedge_A${\hyp}definable subset of finite measure is $A${\hyp}medium.
\end{enumerate} 
\end{ex}
\begin{lem}{\textnormal{{\cite[Lemma 2.9]{hrushovski2011stable}}}} \label{l:s1 and forking} 
Let $P$ be a piecewise $A${\hyp}hyperdefinable set and $\mu$ an ideal of $\bigwedge_{<\lambda}${\hyp}definable subsets on $P$. Let $V\subseteq P$ be a strictly $A${\hyp}medium $\bigwedge_{<\lambda}${\hyp}definable subset. Then, $V$ does not fork over $A$.
\begin{proof} Take an $A${\hyp}medium $\bigwedge_A${\hyp}definable set $V_0$ such that $V\subseteq V_0$. Suppose that $V$ forks over $A$. Then, there are $W'_1,\ldots,W'_n$ $\bigwedge${\hyp}definable subsets dividing over $A$ such that $V\subseteq\bigcup_i W'_i$. Applying \cref{l:finite character of dividing}, we may assume that each $W'_i$ is $\bigwedge_{<\lambda}${\hyp}definable. Now, $V\subseteq \bigcup (W'_i\cap V_0)$ and $W'_i\cap V_0$ divides over $A$. Write $W_i(\overline{b})=W'_i\cap V_0$. By \cref{l:dividing}, there is $(\overline{b}_j)_{j\in\omega}$ $A${\hyp}indiscernible such that $\bigcap_{j\in\omega}W_i(\overline{b}_j)=\emptyset$. Hence, $\bigcap^k_{j=0}W_i(\overline{b}_j)=\emptyset\in \mu$ for some $k$, concluding that $W_i(\overline{b})\in \mu$ by $A${\hyp}mediumness of $V_0$. This concludes that $V\in\mu$, contradicting that $V$ is wide. 
\end{proof}
\end{lem}

\begin{lem}\label{l:translation of medium}
Let $\mathfrak{N}\preceq \mathfrak{M}$ with $|N|<\lambda$ and $G=\underrightarrow{\lim}\, G_k$ be a piecewise $N${\hyp}hyperdefinable group. Let $\mu$ be an $N${\hyp}invariant ideal of $\bigwedge_{<\lambda}${\hyp}definable subsets of $G$. Let $W$ and $U$ be non{\hyp}empty $\bigwedge_N${\hyp}definable subsets of $G$. If $\mu$ is invariant under left translations and $U\cdot W$ is $N${\hyp}medium, then $W$ is $N${\hyp}medium too. Similarly, if $\mu$ is invariant under right translations and $W\cdot U$ is medium over $A$, then $W$ is medium over $A$ too.  
\begin{proof} Let $X(\overline{b})\subseteq W$ be $\bigwedge_{\bar{b}}${\hyp}definable with $|\overline{b}|<\lambda$. Take an $N${\hyp}indiscernible sequence $(\overline{b}^*_i)_{i\in\omega}$ of representatives of elements in $\tp(\overline{b}/N)$ and write $X_i\coloneqq X(\overline{b}_i)$. Fix $k$ such that $U\subseteq G_k$ and take $\underline{U}$ defining $U$. As $\mathfrak{N}\prec\mathfrak{M}$, $\underline{U}$ is finitely satisfiable in $N$, so $\underline{U}$ does not fork over $N$ \cite[Lemma 7.1.10]{tent2012course}. Thus, $\underline{U}$ has a non{\hyp}forking extension to a complete type $p$ over $N,\overline{b}^*_0$. As $p$ does not fork over $N$, it does not divide over $N$. By \cite[Lemma 7.1.5]{tent2012course}, there is $a^*$ realising $p$ such that $(\overline{b}^*_i)_{i\in \omega}$ is $N,a^*${\hyp}indiscernible. As $a^*$ realises $p$, in particular, $a\in U$. 

Suppose $\mu$ is invariant under left translations and $U\cdot W$ is $N${\hyp}medium. Then, we get that $a\cdot X_0\cap a\cdot X_1\in \mu$ if and only if $a\cdot X_0\in\mu$. Therefore, provided left translational invariance, 
\[X_0\cap X_1\in \mu\Leftrightarrow a\cdot X_0\cap a\cdot X_1\in \mu\Leftrightarrow a\cdot X_0\in\mu\Leftrightarrow X_0\in\mu.\]

We similarly prove the case with right translations.
\end{proof}
\end{lem}
\subsection{Stable relations}
Let $P$ and $Q$ be piecewise $A${\hyp}hyperdefinable sets and $V,W\subseteq P\times Q$ disjoint $A${\hyp}invariant subsets. We say that $V$ and $W$ are \emph{unstably separated over $A$} if there is an infinite $A${\hyp}indiscernible sequence $(a_i,b_i)_{i\in \omega}$ such that $(a_0,b_1)\in V$ and $(a_1,b_0)\in W$. We say that they are \emph{stably separated over $A$} if they are not unstably separated. We say that an $A${\hyp}invariant binary relation $R\subseteq P\times Q$ is \emph{stable over $A$} if $R$ and $R^c$ are stably separated over $A$. 
\begin{rmk} 
\begin{enumerate}[label={\rm{(\arabic*)}}, wide] 
\item[\hspace{-1.4em}\setcounter{enumi}{1}\theenumi] Clearly, being stably separated over $A$ is invariant under piecewise $\bigwedge_A${\hyp}definable isomorphisms.
\item Note that $V$ and $W$ are stably separated over $A$ if and only if $V\cap (P_i\times Q_j)$ and $W\cap (P_i\times Q_j)$ are stably separated over $A$ for any pieces $P_i$ and $Q_j$.
\item Let $P$ and $Q$ be $A${\hyp}hyperdefinable. By \cref{c:indiscernible hyperimaginaries}, $V$ and $W$ are stably separated over $A$ if and only if $\quot^{-1}_{P\times Q}[V]$ and $\quot^{-1}_{P\times Q}[W]$ are stably separated over $A^*$ for any set of representatives $A^*$ of $A$. In particular, if $V$ and $W$ are $\bigwedge_A${\hyp}definable, we have that $V$ and $W$ are stably separated over $A$ if and only if $\underline{V}$ and $\underline{W}$ are stably separated over $A^*$ for any set of representatives $A^*$ of $A$ and any partial types $\underline{V}$ and $\underline{W}$ defining $V$ and $W$.
\item Note that, by the symmetry of stability, $R$ is stable over $A$ if and only if $R^c$ is stable over $A$. 
\end{enumerate}
\end{rmk}
We need the following two lemmas from \cite{hrushovski2011stable} for the Stabilizer Theorem. Their original proofs can be easily adapted. 

Recall that we say that a partial type $\Sigma(x,y)$ over $A^*$ divides over $A^*$ \emph{with respect to} a global $A^*${\hyp}invariant type $\widehat{q}(y)$ if there is a sequence $(b^*_i)_{i\in \omega}$ such that $\tp(b^*_i/A^*,b^*_0,\ldots,b^*_{i-1})\models \widehat{q}$ and $\bigwedge \Sigma(x,b^*_i)$ is inconsistent.  
\begin{lem}{\textnormal{{\cite[Lemma 2.2]{hrushovski2011stable}}}} 
Let $P=\rfrac{X}{E}$ and $Q=\rfrac{Y}{F}$ be $A^*${\hyp}hyperdefinable sets and $\widehat{q}$ a global $A^*${\hyp}invariant type with $\underline{Y}\in\widehat{q}$. Let $W_1,W_2\subseteq P\times Q$ be $\bigwedge_{A^*}${\hyp}definable sets stably separated over $A^*$. Assume that $a\in P$ is such that $\underline{W_2(a)}\subseteq \widehat{q}$. Then, $V=W_1\cap \left(\tp(a/A^*)\times Q\right)$ divides over $A^*$. Furthermore, $\underline{V}$ divides over $A^*$ with respect to $\widehat{q}$.
\end{lem}
For sets $V$ and $W$, write 
\[V\times_{\mathrm{nf}(A)}W\coloneqq\{(a,b)\in V\times W\sth \tp(b/aA)\mbox{ does not fork over }A\}.\] 
Similarly, $V\times_{\mathrm{ndiv}(A)}W$, $V\hbox{\hss $\prescript{}{\mathrm{nf}(A)}{\times}$ \hss} W$ and $V\hbox{\hss $\prescript{}{\mathrm{ndiv}(A)}{\times}$ \hss} W$.
\begin{lem}{\textnormal{{\cite[Lemma 2.3]{hrushovski2011stable}}}} \label{l:fundamental lemma for the stabilizer theorem}
Let $P$ and $Q$ be piecewise $A^*${\hyp}hyperdefinable sets and $R\subseteq P\times Q$ a stable binary relation over $A^*$. Let $p\subseteq P$ and $q\subseteq Q$ be types over $A^*$. Assume that there is an $A^*${\hyp}invariant global type $\widehat{q}$ extending a partial type $\underline{q}$ over $A^*$ defining $q$.
\begin{enumerate}[label={\rm{(\arabic*)}}, ref={\rm{\arabic*}}, wide]
\item \label{itm:fundamental lemma for the stabilizer theorem 1} Take $a\in p$ and $b\in q$ such that $(a,b)\in R$. Suppose that there are representatives $a^*\in a$ and $b^*\in b$ such that $b^*\models\widehat{q}_{\mid A^*,a^*}$. Then, we have that $(a',b)\in R$ for any $a'\in p$ such that $\tp(a'/A^*,b)$ does not divide over $A^*$.
\item \label{itm:fundamental lemma for the stabilizer theorem 2} Take $a',a\in p$ and $b\in q$. Suppose that $\tp(a/A^*,b)$ and $\tp(a'/A^*,b)$ do not divide over $A^*$. Then, $(a,b)\in R$ if and only if $(a',b)\in R$.
\item \label{itm:fundamental lemma for the stabilizer theorem 3} Assume that there is also an $A^*${\hyp}invariant global type $\widehat{p}$ extending a partial type $\underline{p}$ over $A^*$ defining $p$. Then, the following conditions are equivalent:
\[\begin{array}{llll}
\mathrm{i.}& p\times_{\mathrm{ndiv}(A^*)}q\cap R\neq \emptyset & \mathrm{v.}& p\hbox{\hss $\prescript{}{\mathrm{ndiv}(A^*)}{\times}$ \hss} q\cap R\neq \emptyset\\
\mathrm{ii.}& p\times_{\mathrm{nf}(A^*)}q\cap R\neq \emptyset & \mathrm{vi.}& p\hbox{\hss $\prescript{}{\mathrm{nf}(A^*)}{\times}$ \hss} q\cap R\neq \emptyset\\
\mathrm{iii.}& p\times_{\mathrm{nf}(A^*)}q\subseteq R &\mathrm{vii.} & p\hbox{\hss $\prescript{}{\mathrm{nf}(A^*)}{\times}$ \hss} q\subseteq R\\
\mathrm{iv.}& p\times_{\mathrm{ndiv}(A^*)}q\subseteq R &\mathrm{viii.} & p\hbox{\hss $\prescript{}{\mathrm{ndiv}(A^*)}{\times}$ \hss} q\subseteq R\end{array}\]
\end{enumerate}
\end{lem}
\subsection{Stabilizer Theorem}
From now on, $\lambda$ is a cardinal with $\kappa\geq\lambda>|\lang|+|A|$.\medskip

Let $G$ be a piecewise $A${\hyp}hyperdefinable group and $\mu$ an ideal of $\bigwedge_{<\lambda}${\hyp}definable subsets of $G$. Let $V,W\subseteq G$ be two $\bigwedge_{<\lambda}${\hyp}definable subsets. The \emph{(left) $\mu${\hyp}stabilizer of $V$ with respect to $W$} is the set $\St_\mu(V,W)=\{a\in G\sth a^{-1}\cdot V\cap W\notin \mu\}$, and the \emph{(left) $\mu${\hyp}stabilizer group of $V$ with respect to $W$} is the subgroup $\Stab_\mu(V,W)$ of $G$ generated by $\St_\mu(V,W)$. Write $\St_\mu(V)\coloneqq \St_\mu(V,V)$ and $\Stab_\mu(V)\coloneqq \St_\mu(V,V)$. We omit the subscript $\mu$ if there is no confusion. 
\begin{lem} \label{l:basic remarks} 
Let $G$ be a piecewise $A${\hyp}hyperdefinable group, $\mu$ an $A${\hyp}invariant ideal of $\bigwedge_{<\lambda}${\hyp}definable subsets of $G$ and $V$ and $W$ two $\bigwedge_A${\hyp}definable wide subsets. Then:
\begin{enumerate}[label={\rm{(\arabic*)}}, ref={\rm{\arabic*}}, wide]
\item \label{itm:basic remarks:locally atomic} Take any $B\supseteq A$ with $|B|<\lambda$. If $\mu$ is locally atomic, $g\in \St(V,W)$ if and only if there are $a\in V$ and $b\in W$ such that $g=ab^{-1}$ with $\tp(b/B,g)\notin \mu$. 
\item \label{itm:basic remarks:invariant under translations} If $\mu$ is invariant under left translations, then $\St(V,W)^{-1}=\St(W,V)$. In particular, $\St(V)$ is a symmetric subset.
\end{enumerate}
\end{lem}
We say that a subset $X\subseteq G$ is \emph{stable over $A$} if the relation $\{(a,b)\sth a^{-1}b\in X\}$ is stable over $A$. 
\begin{ex}{\textnormal{{\cite[Lemma 2.10]{hrushovski2011stable}\&\cite[Lemma 2.8]{montenegro2018stabilizers}}}} \label{e:mos 2.8} Let $\mu$ be an $A${\hyp}invariant ideal of $\bigwedge_{<\lambda}${\hyp}definable subsets of a piecewise $A${\hyp}hyperdefinable group. If $\mu$ is invariant under left translations and $V$ and $W$ are $A${\hyp}medium $\bigwedge_A${\hyp}definable subsets, then $\St(V,W)$ is stable over $A$. 
\end{ex}
\begin{rmk} Since ${\scriptstyle \bullet}^{-1}:\ G\rightarrow G$ is a piecewise $\bigwedge_A${\hyp}definable isomorphism, $X$ is stable over $A$ if and only if $X^{-1}$ is stable over $A$. Also, $X$ is stable over $A$ if and only if $\{(a,b)\sth ab\in X\}$ is stable over $A$, if and only if $\{(a,b)\sth ab^{-1}\in X\}$ is stable over $A$.
\end{rmk} 
For subsets $V$ and $W$ of a piecewise hyperdefinable group $G$, write 
\[V\cdot_{\mathrm{nf}(A)}W\coloneqq \{a\cdot b\sth (a,b)\in V\times_{\mathrm{nf}(A)}W\}.\] 
Similarly, $V\cdot_{\mathrm{ndiv}(A)}W$, $V\hbox{\hss $\prescript{}{\mathrm{nf}(A)}{\cdot}$ \hss} W$ and $V\hbox{\hss $\prescript{}{\mathrm{ndiv}(A)}{\cdot}$ \hss} W$.
\begin{lem} \label{l:lemma stabilizer} 
Let $G$ be a piecewise $A^*${\hyp}hyperdefinable group and $\mu$ a locally atomic $A^*${\hyp}invariant ideal of $\bigwedge_{<\lambda}${\hyp}definable subsets of $G$ which is also invariant under left translations. Let $p$ be an $A^*${\hyp}medium $A^*${\hyp}type and $X\subseteq G$ a stable subset over $A^*$ containing $p$. Suppose that there is an $A^*${\hyp}invariant global type $\widehat{p}$ extending a partial type $\underline{p}$ over $A^*$ defining $p$. Then, $\Stab(p)\subseteq X\cdot X^{-1}$. Furthermore, $\Stab(p)\cdot_{\mathrm{nf}(A^*)} p\subseteq X$.
\begin{proof} We may assume that $p$ is strictly $A${\hyp}medium; otherwise, $\St(p)=\emptyset$ and $\Stab(p)=\{1\}$, so the lemma holds trivially. Since $\mu$ is invariant under left translations, $\St(p)^{-1}=\St(p)$ by \cref{l:basic remarks}(\ref{itm:basic remarks:invariant under translations}). Thus, by definition, $\Stab(p)=\bigcup^{\infty}_{n=0} \St(p)^n$. We prove by induction on $n$ that $b\cdot c\in X$ for $b\in\St(p)^n$ and $c\in p$ such that $\tp(c/A^*,b)$ does not fork over $A^*$.

By hypothesis, $p\subseteq X$, so we are done for $n=0$. Assume it is true for $n-1$ with $n\in\N_{>0}$. Take $b=b_1\cdot \cdots\cdot b_n$ with $b_i\in \St(p)$ for each $i\in\{1,\ldots,n\}$ and $c\in p$ such that $\tp(c/A^*,b)$ does not fork over $A^*$. We want to prove that $b\cdot c\in X$. 
As $b_n\in\St(p)$, by \cref{l:basic remarks}(\ref{itm:basic remarks:locally atomic}), there is $c'\in p$ such that $b_n\cdot c'\in p$ and $\tp(c'/A^*,b_1,\ldots,b_n)\notin \mu$. In particular, as $p$ is $A^*${\hyp}medium, $\tp(c'/A^*,b_1,\ldots,b_n)$ does not fork over $A^*$ by \cref{l:s1 and forking}. Hence, $\tp(c'/A^*,b)$ and $\tp(b_n\cdot c'/A^*,b_1,\ldots,b_{n-1})$ do not fork over $A^*$. By induction hypothesis, $b\cdot c'=(b_1\cdots b_{n-1})\cdot b_n\cdot c'\in X$. Since $(b,c')\in \tp(b/A^*)\times_{\mathrm{nf}(A^*)}p$, by \cref{l:fundamental lemma for the stabilizer theorem}(\ref{itm:fundamental lemma for the stabilizer theorem 2}), we conclude that $b\cdot c\in X$ too.

Now, for any $b\in \Stab(p)$ we can find $c\in p$ such that $\tp(c/A^*,b)$ does not fork over $A^*$ --- namely, choose $c^*\models \widehat{p}_{\mid A^*,b^*}$. Therefore, we conclude that $\Stab(p)\subseteq X\cdot p^{-1}\subseteq X\cdot X^{-1}$.
\end{proof}
\end{lem} 
As an immediate corollary we get the following result:
\begin{coro} Let $\mathfrak{N}\prec\mathfrak{M}$ with $|N|<\lambda$. Let $G$ be a piecewise $N${\hyp}hyperdefinable group and $\mu$ a locally atomic $N${\hyp}invariant ideal of $\bigwedge_{<\lambda}${\hyp}definable subsets invariant under left translations. Let $p$ be an $N${\hyp}medium $N${\hyp}type with $p\subseteq \St(p)$. Then, $\Stab(p)=\St(p)^2=(p\cdot p^{-1})^2$.
\begin{proof} If $p\subseteq \St(p)$, by \cref{l:basic remarks}, we get that $\St(p)\subseteq p\cdot p^{-1}\subseteq \St(p)^2$. By \cref{e:mos 2.8}, $\St(p)$ is a stable subset. Hence, by \cref{l:lemma stabilizer}, taking $\widehat{p}$ a coheir of $\underline{p}$, we conclude that $\Stab(p)\subseteq \St(p)^2\subseteq (p\cdot p^{-1})^2\subseteq \Stab(p)$.
\end{proof}  
\end{coro}
\begin{lem}{\textnormal{{\cite[Lemma 2.11]{montenegro2018stabilizers}}}} \label{l:mos 2.11}
Let $G$ be a piecewise $A^*${\hyp}hyperdefinable group and $\mu$ an $A^*${\hyp}invariant ideal of $\bigwedge_{<\lambda}${\hyp}definable subsets of $G$ which is also invariant under left translations. Let $p$ be an $A^*${\hyp}type and $W$ a strictly $A^*${\hyp}medium $\bigwedge_{A^*}${\hyp}definable subset such that $p^{-1}W$ is $A^*${\hyp}medium too. Suppose that there is an $A^*${\hyp}invariant global type $\widehat{p}$ extending a partial type $\underline{p}$ over $A^*$ defining $p$. Then, $p\cdot_{\mathrm{nf}(A^*)} p^{-1}\subseteq \St(W)$.
\begin{proof} Take $(a^*_i)_{i\in \omega}$ such that $a^*_i\models \widehat{p}_{\mid A^*,a^*_0,\ldots,a^*_{i-1}}$. Then,   $(a_i)_{i\in \omega}$ is an $A^*${\hyp}indiscernible sequence of realisations of $p$ with $\tp(a_1/A^*,a_0)$ non{\hyp}forking. As $\mu$ is invariant under left translations, $a^{-1}_0W$ is wide. Since $p^{-1}W$ is $A^*${\hyp}medium, we conclude that $a^{-1}_0W\cap a^{-1}_1W\not\in\mu$. Thus, by invariance under left translations, $a_0a^{-1}_1\in \St(W)$. By \cref{l:fundamental lemma for the stabilizer theorem}(\ref{itm:fundamental lemma for the stabilizer theorem 2}) and \cref{e:mos 2.8}, we conclude that $p\cdot_{\mathrm{nf}(A^*)}p^{-1}\subseteq \St(W)$.  
\end{proof}
\end{lem}
For two subsets $X$ and $Y$ of a group $G$, recall that $[X:Y]\coloneqq \min\{|\Delta|\sth X\subseteq \Delta Y\}$. Note that $[X^{-1}:Y^{-1}]=\min\{|\Delta|\sth X\subseteq Y\Delta\}$.
\begin{lem} \label{l:index lemma 1}
Let $G$ be a piecewise $A^*${\hyp}hyperdefinable group and $S\leq G$ an $\bigwedge_{A^*}${\hyp}definable subgroup. Let $p\subseteq G$ be an $A^*${\hyp}type such that there is an $A^*${\hyp}invariant global type $\widehat{p}$ extending a partial type $\underline{p}$ defining $p$ over $A^*$.
\begin{enumerate}[label={\rm{(\arabic*)}}, ref={\rm{\arabic*}}, wide]
\item \label{itm:index lemma 1:left} If $[p:S]$ is small, then $[p:S]=1$ and $p^{-1}\cdot p\subseteq S$.
\item \label{itm:index lemma 1:right} If $[p^{-1}:S]$ is small, then $[p^{-1}:S]=1$ and $p\cdot p^{-1}\subseteq S$. 
\end{enumerate}
\begin{proof} Both are very similar, so we just prove {(\ref*{itm:index lemma 1:left})}. Let $(a_i)_{i\in \alpha}$ be a sequence in $p$ such that $\alpha=[p:S]$ and $a_jS\cap a_iS=\emptyset$ for $i\neq j$, and $p\subseteq \bigcup_{i\in\alpha}a_iS$. Suppose that there is no $i\in\alpha$ such that $\underline{(a_iS)\cap p}$ is contained in $\widehat{p}$. Then, there are formulas $\varphi_i\in \underline{(a_iS)\cap p}$ for each $i\in\alpha$ such that $\{\neg\varphi_i\}_{i\in\alpha}\subseteq \widehat{p}$. In particular, $\{\neg\varphi_i\}_{i\in\alpha}\cup\underline{p}$ is finitely satisfiable. Thus, by saturation, there is $c\in p$ such that $c\notin \bigcup_{i\in \alpha} a_iS$, getting a contradiction. That concludes that there is $i\in \alpha$ such that $\underline{(a_iS)\cap p}\subseteq \widehat{p}$. Since $a_iS\cap a_jS=\emptyset$ for any $i\neq j$, we conclude that there is just one such $i\in\alpha$. Since $\widehat{p}$ and $S$ are $A^*${\hyp}invariant, we get that $a_iS$ is $A^*${\hyp}invariant too. Indeed, take $\sigma\in\Aut(\mathfrak{M}/A^*)$ arbitrary. Then, $\sigma[\underline{(a_iS)\cap p}]\subseteq \widehat{p}$ defines $(\sigma(a_i)S)\cap p$, concluding that $\sigma(a_i)S=a_iS$. As $\sigma$ is arbitrary, we conclude that $a_iS$ is $A^*${\hyp}invariant. Thus, $a_iS$ is $\bigwedge_{A^*}${\hyp}definable by \cref{c:parameters of definition}. As $(a_iS)\cap p\neq \emptyset$, by minimality of $p$, $p\subseteq a_iS$. So, $[p:S]=1$ and $p^{-1}\cdot p\subseteq S$.

For (\ref*{itm:index lemma 1:right}), one argues analogously to point (\ref*{itm:index lemma 1:left}) but working with right cosets --- as $[p^{-1}:S]=\min\{|\Delta|\sth p\subseteq S\Delta\}$. 
\end{proof}  
\end{lem}
\begin{lem} \label{l:index lemma 2}
Let $G$ be a piecewise $A${\hyp}hyperdefinable group, $S\leq G$ an $\bigwedge_A${\hyp}definable subgroup and $V\subseteq G$ an $\bigwedge_A${\hyp}definable set. Suppose $[V:S]$ is not small. Then, there is an $A${\hyp}indiscernible sequence $(a_i)_{i\in\omega}$ in $V$ such that $a_i\cdot S\cap a_j\cdot S=\emptyset$ for $i\neq j$. 
\begin{proof}
Let $\mathbf{G}$, $\mathbf{S}$ and $\mathbf{V}$ be the interpretations in the monster model $\boldsymbol{\mathfrak{C}}$. Since $\boldsymbol{\mathfrak{C}}$ is the monster model, $[\mathbf{V}: \mathbf{S}]\notin\Ord$. Recall that there is $\tau\in\Ord$ depending only on $|A^*|$ such that, for any sequence $(a^*_i)_{i\in\tau}$ of elements in $\boldsymbol{\mathfrak{C}}$, there is an $A^*${\hyp}indiscernible sequence $(b^*_j)_{j\in\omega}$ of elements in $\boldsymbol{\mathfrak{C}}$ with the property that, for any $j_1<\cdots<j_n$, $\tp(a^*_{i_1},\ldots,a^*_{i_n}/A^*)=\tp(b^*_{j_1},\ldots,b^*_{j_n}/A^*)$ for some $i_1<\cdots<i_n$ --- see \cite[Lemma 7.2.12]{tent2012course}. In particular, take $(a_i)_{i\in \tau}$ a sequence of hyperimaginaries of $\mathbf{V}$ such that $a_i\cdot \mathbf{S}\cap a_j\cdot \mathbf{S}=\emptyset$ for each $i\neq j$. Let $(a^*_i)_{i\in\tau}$ be representatives of $(a_i)_{i\in \tau}$. Then, there is a sequence $(\widetilde{b}^*_j)_{j\in\omega}$ of $A^*${\hyp}indiscernible elements such that $\tp(\widetilde{b}^*_0,\widetilde{b}^*_1/A^*)=\tp(a^*_{j'},a^*_{j''}/A^*)$  for some $j'<j''$. Now, by $\kappa${\hyp}saturation of $\mathfrak{M}$, we can find $(b^*_j)_{j\in\omega}$ elements in $\mathfrak{M}$ realising the same type that $(\widetilde{b}^*_j)_{j\in\omega}$ over $A^*$. So, the projections $b_j\in V$ form an $A${\hyp}indiscernible sequence $(b_j)_{j\in\omega}$ in $V$ such that $b_i\cdot S\cap b_j\cdot S=\emptyset$ for $i\neq j$. 
\end{proof}
\end{lem}
We now prove the Stabilizer Theorem for piecewise hyperdefinable groups. Below, \cref{c:stabilizer theorem mos b2} corresponds to \cite[Theorem 2.12 (B2)]{montenegro2018stabilizers}; \cref{t:stabilizer theorem 2} corresponds to \cite[Theorem 3.5]{hrushovski2011stable} and \cite[Theorem 2.12 (B1)]{montenegro2018stabilizers}; and \cref{t:stabilizer theorem 2}(\ref{itm:stabilizer theorem 2:c}) corresponds to \cite[Proposition 2.14]{montenegro2018stabilizers}.
\begin{theo} \label{t:stabilizer theorem 1} Let $\mathfrak{N}\prec\mathfrak{M}$ with $|N|<\lambda$. Let $G$ be a piecewise $N${\hyp}hyperdefinable group and $\mu$ a locally atomic $N${\hyp}invariant ideal of $\bigwedge_{<\lambda}${\hyp}definable subsets invariant under left translations. Let $p\subseteq G$ be an $N${\hyp}medium $N${\hyp}type. Assume that there is a wide $N${\hyp}type $q\subseteq \St(p)$ such that $p^{-1}\cdot q$ is $N${\hyp}medium. Then, $\Stab(p)=\Stab(q)=\St(p)^2=\St(q)^4=(p\cdot p^{-1})^2$ is a wide $\bigwedge_N${\hyp}definable subgroup of $G$ without proper $\bigwedge_N${\hyp}definable subgroups of small index. Furthermore:
\begin{enumerate}[label={\rm{(\alph*)}}, ref={\rm{\alph*}}, wide]
\item \label{itm:stabilizer theorem 1:a} $\Stab(q)\cdot_{\mathrm{nf}(N)}q\subseteq \St(p)$.
\item \label{itm:stabilizer theorem 1:b} Every strictly $N${\hyp}medium $N${\hyp}type of $\Stab(p)$ lies in $\St(p)$.  
\end{enumerate}
\begin{proof} We take $\widehat{p}$ and $\widehat{q}$ coheirs of $\underline{p}$ and $\underline{q}$. By \cref{l:translation of medium}, $q$ is $N${\hyp}medium. By \cref{l:mos 2.11}, $p\cdot_{\mathrm{nf}(N)}p^{-1}\subseteq \St(q)$. Then, $\St(p)\subseteq p\cdot p^{-1}\subseteq \St(q)^2$. Indeed, for any $a,b\in p$ we can find $c\in p$ such that $\tp(c/N,a,b)$ does not fork over $N$ --- simply, take $c^*$ realising $\widehat{p}_{\mid N,a^*,b^*}$. Thus, $ac^{-1},bc^{-1}\in \St(q)$, concluding that $ab^{-1}\in\St(q)^2$ by \cref{l:basic remarks}(\ref{itm:basic remarks:invariant under translations}). 

Since $q\subseteq \St(p)$, $\Stab(q)\cdot_{\mathrm{nf}(N)} q\subseteq \St(p)$ by \cref{l:lemma stabilizer,e:mos 2.8}. In particular, $\Stab(q)\subseteq \St(p)^2$ by \cref{l:basic remarks}(\ref{itm:basic remarks:invariant under translations}). Thus, $\Stab(p)=\Stab(q)=\St(p)^2=(p\cdot p^{-1})^2=\St(q)^4$ is an $\bigwedge_N${\hyp}definable subgroup. Since $q\subseteq \St(p)\subseteq \Stab(p)$, we have that $\Stab(p)$ is wide. 

Take an $\bigwedge_N${\hyp}definable subgroup $T\leq \Stab(p)$ such that $[\Stab(p):T]$ is small. Since $p\cdot p^{-1}\subseteq \Stab(p)$, $[p^{-1}:T]$ is also small. By \cref{l:index lemma 1}(\ref{itm:index lemma 1:right}), $p\cdot p^{-1}\subseteq T$. Therefore, $\Stab(p)=(p\cdot p^{-1})^2\subseteq T$. In other words, $\Stab(p)$ does not have proper $\bigwedge_N${\hyp}definable subgroups of small index.

Finally, we prove property (\ref*{itm:stabilizer theorem 1:b}). Take $r\subset \Stab(p)$ a strictly $N${\hyp}medium $N${\hyp}type. Set $c\in q$. Since $\mu$ is locally atomic, there is $b\in r$ such that $\tp(b/N,c)\notin \mu$. By \cref{l:s1 and forking}, $\tp(b/N,c)$ does not fork over $N$. Then, by \cref{l:dividing and functions,l:types and infinite definable functions}, $\tp(b^{-1}c^{-1}/N,c)$ does not fork over $N$. Since $c,b\in \Stab(q)$, we have $b^{-1}c^{-1}\in\Stab(q)$. Write $r'=\tp(b^{-1}c^{-1}/N)\subseteq \Stab(q)$. By {\rm{(a)}}, $r'\cdot_{\mathrm{nf}(N)}q\subseteq \St(p)$. By \cref{l:fundamental lemma for the stabilizer theorem,e:mos 2.8}, we conclude that $b^{-1}=b^{-1}\cdot c^{-1}\cdot c\in \St(p)$, so $b\in \St(p)$. Since $r$ is a type over $N$, we conclude $r\subseteq \St(p)$. 
\end{proof}
\end{theo} 
\begin{coro} \label{c:stabilizer theorem mos b2}
Let $\mathfrak{N}\prec\mathfrak{M}$ with $|N|<\lambda$ and $G$ a piecewise $N${\hyp}hyperdefinable group. Let $\mu$ be an $N${\hyp}invariant locally atomic ideal of $\bigwedge_{<\lambda}${\hyp}definable subsets of $G$ invariant under left and right translations. Let $p$ be a wide type over $N$ and assume that $p^{-1}\cdot p\cdot p^{-1}$ is $N${\hyp}medium. Then, $\Stab(p)=\St(p)^2=(p\cdot p^{-1})^2$ is a wide $\bigwedge_N${\hyp}definable subgroup of $G$ without proper $\bigwedge_N${\hyp}definable subgroups of small index such that every strictly $N${\hyp}medium $N${\hyp}type of $\Stab(p)$ lies in $\St(p)$.
\begin{proof} Take $a\in p$ arbitrary. By the local atomic property, there is $b\in p$ such that $\tp(b/N,a)$ is wide. By \cref{l:translation of medium}, $p\cdot p^{-1}$ is $N${\hyp}medium. Again by \cref{l:translation of medium} but now using invariance under right translations, we get that $p$ is $N${\hyp}medium. By \cref{l:s1 and forking}, it follows that $\tp(b/N,a)$ does not fork over $N$. Consider the $N${\hyp}type $q=\tp(ba^{-1}/N)\subseteq p\prescript{}{\mathrm{nf}(N)}{\cdot} p^{-1}$. 

By invariance under right translations, we know that $q$ is a wide $\bigwedge_N${\hyp}definable set. As $q\subseteq p\cdot p^{-1}$ with $p\cdot p^{-1}$ $N${\hyp}medium, we have that $q$ is $N${\hyp}medium too. Also, note that $p^{-1}\cdot q\subseteq p^{-1}\cdot p\cdot p^{-1}$, so $p^{-1}\cdot q$ is $N${\hyp}medium too. On the other hand, by \cref{l:translation of medium} using invariance under right translations, $p^{-1}\cdot p$ is $N${\hyp}medium, so $p\cdot_{\mathrm{nf}(N)}p^{-1}\subseteq \St(p)$ by \cref{l:mos 2.11}. By \cref{e:mos 2.8} and \cref{l:fundamental lemma for the stabilizer theorem}(\ref{itm:fundamental lemma for the stabilizer theorem 3}), $p\prescript{}{\mathrm{nf}(N)}{\cdot} p^{-1}\subseteq \St(p)$ too, so $q\subseteq \St(p)$. By \cref{t:stabilizer theorem 1}, we conclude that $\Stab(p)=\St(p)^2=(p\cdot p^{-1})^2$ is a wide $\bigwedge_N${\hyp}definable subgroup of $G$ without proper $\bigwedge_N${\hyp}definable subgroups of small index such that every strictly $N${\hyp}medium $N${\hyp}type of $\Stab(p)$ lies in $\St(p)$. 
\end{proof}
\end{coro}
\begin{theo}[Stabilizer Theorem]  \label{t:stabilizer theorem 2}
Let $\mathfrak{N}\prec\mathfrak{M}$ with $|N|<\lambda$ and $G=\underrightarrow{\lim}\, X^n$ a piecewise $N${\hyp}hyperdefinable group generated by a symmetric $\bigwedge_N${\hyp}definable subset $X$. Let $\mu$ be an $N${\hyp}invariant locally atomic ideal of $\bigwedge_{<\lambda}${\hyp}definable subsets of $G$ invariant under left translations. Let $p\subseteq X$ be a wide type over $N$ and assume that $X^3$ is $N${\hyp}medium. Then, $\Stab(p)=\St(p)^2=(p\cdot p^{-1})^2$ is a wide and $N${\hyp}medium $\bigwedge_N${\hyp}definable normal subgroup of small index of $G$ without proper $\bigwedge_N${\hyp}definable subgroups of small index. Furthermore: 
\begin{enumerate}[label={\rm{(\alph*)}}, ref={\rm{\alph*}}, wide]
\item \label{itm:stabilizer theorem 2:a} $p\cdot p\cdot p^{-1}=p\cdot \Stab(p)$ is a left coset of $\Stab(p)$.  
\item \label{itm:stabilizer theorem 2:b} Every wide $N${\hyp}type of $\Stab(p)$ is contained in $\St(p)$.
\item \label{itm:stabilizer theorem 2:c} Assume $\mu$ is also invariant under right translations. Then, $p\cdot p^{-1}\cdot p=\Stab(p)\cdot p$ is a right coset of $\Stab(p)$.
\end{enumerate}
\begin{proof} Take $a\in p$ arbitrary. By the local atomic property, we find $b\in p$ such that $\tp(b/N,a)$ is wide. By \cref{l:translation of medium}, $X$ is $N${\hyp}medium. As $\tp(b/N,a)\subseteq p\subseteq X$, we conclude that $\tp(b/N,a)$ is $N${\hyp}medium, so $\tp(b/N,a)$ does not fork over $N$ by \cref{l:s1 and forking}. Also, $\tp(a^{-1}b/N,a)=a^{-1}\cdot \tp(b/N,a)$ is wide by left invariance of $\mu$, so $q=\tp(a^{-1}b/N)\subseteq p^{-1}\cdot_{\mathrm{nf}(N)}p$ is wide. By \cref{l:translation of medium}, $X^2$ is $N${\hyp}medium, so $p^2$ is $N${\hyp}medium too. Using any coheir extending $\underline{p^{-1}}$, we get $p^{-1}\cdot_{\mathrm{nf}(N)}p\subseteq \St(p)$ by \cref{l:mos 2.11}. In particular, $q\subseteq \St(p)$. As $X$ is symmetric, $p^{-1}\cdot q\subseteq p^{-1}\cdot p^{-1}\cdot p\subseteq X^3$, so $p^{-1}\cdot q$ is $N${\hyp}medium. Thus, by \cref{t:stabilizer theorem 1}, $\Stab(p)=\Stab(q)=\St(p)^2=\St(q)^4=(p\cdot p^{-1})^2$ is a wide $\bigwedge_N${\hyp}definable subgroup without proper $\bigwedge_N${\hyp}definable subgroups of small index. By \cref{t:stabilizer theorem 1}(\ref{itm:stabilizer theorem 1:b}), every strictly $N${\hyp}medium type of $\Stab(p)$ is contained in $\St(p)$. Also, $\Stab(q)\cdot_{\mathrm{nf}(N)}q\subseteq \St(p)$ by \cref{t:stabilizer theorem 1}(\ref{itm:stabilizer theorem 1:a}).

We show that it is a normal subgroup of small index. First of all, note that $[X^2:\Stab(p)]$ is small. Otherwise, by \cref{l:index lemma 2}, there is an $N${\hyp}indiscernible sequence $(b_j)_{j\in\omega}$ in $X^2$ such that $b_i\cdot \Stab(p)\cap b_j\cdot \Stab(p)=\emptyset$ for $i\neq j$. Take $d\in p^{-1}$ arbitrary. As $p\cdot p^{-1}\subseteq \Stab(p)$, we get $(b_i\cdot p\cdot d)\cap (b_j\cdot p\cdot d)=\emptyset$ for $i\neq j$, so $b_i\cdot p\cap b_j\cdot p=\emptyset$ for $i\neq j$. As $X^3$ is $N${\hyp}medium, that implies $b_0\cdot p\in\mu$, so $p\in\mu$ by invariance under left translations, contradicting our hypotheses. For any $c\in X$, we have $[\tp(c/N):\Stab(p)]=1$ by \cref{l:index lemma 1}(\ref{itm:index lemma 1:left}), so $\tp(c/N)\cdot \Stab(p)\cdot \tp(c/N)^{-1}=c\cdot \Stab(p)\cdot c^{-1}$ is $\bigwedge_N${\hyp}definable. As $X$ is symmetric, we also get that $[p^{-1}:c\cdot \Stab(p)\cdot c^{-1}]=[p^{-1}\cdot c:\Stab(p)]\leq [X^2:\Stab(p)]$ is small, so $p\cdot p^{-1}\subseteq c\cdot \Stab(p)\cdot c^{-1}$ by \cref{l:index lemma 1}(\ref{itm:index lemma 1:right}). Therefore, $\Stab(p)\subseteq c\cdot \Stab(p)\cdot c^{-1}$. As $c\in X$ is arbitrary and $X$ is symmetric, $\Stab(p)=c\cdot \Stab(p)\cdot c^{-1}$. Thus, $X\subseteq \mathrm{N}_G(\Stab(p))$, concluding $\Stab(p)\trianglelefteq G$. Since we have proved that $[X:\Stab(p)]$ is small, we conclude that $\Stab(p)$ has small index by normality.

We show now property (\ref*{itm:stabilizer theorem 2:a}), i.e. $p\cdot p\cdot p^{-1}=p\cdot \Stab(p)=y\cdot\Stab(p)$ for any $y\in p$. As $(p\cdot p^{-1})^2=\Stab(p)$, we have $p\cdot p\cdot p^{-1}\subseteq p\cdot \Stab(p)$. As $[p:\Stab(p)]\leq [X:\Stab(p)]$ is small, we get $[p:\Stab(p)]=1$ by \cref{l:index lemma 1}(\ref{itm:index lemma 1:left}). Thus, $p\cdot \Stab(p)=y\cdot \Stab(p)$. On the other hand, take $x\in y\cdot \Stab(p)$ arbitrary, so $x^{-1}=c\cdot y^{-1}$ with $c\in\Stab(p)=\Stab(q)$. By definition of $q$, we know that $q=\tp(a^{-1}b/N)$ with $\tp(b/N,a)$ wide and $a,b\in p$. Take $z_0$ such that $\tp(z_0,y/N)=\tp(b,a/N)$. By $N${\hyp}invariance of $\mu$, $\tp(z_0/N,y)$ is wide. By the local atomic property, we can find $z\in \tp(z_0/N,y)$ such that $\tp(z/N,y,c)$ is wide. Thus, $\tp(y^{-1}z/N,c)$ is wide by invariance under left translations. By \cref{l:s1 and forking}, $\tp(y^{-1}z/N,c)$ does not fork over $N$. Thus, $x^{-1}\cdot z=c\cdot y^{-1}z\in \Stab(q)\cdot_{\mathrm{nf}(N)}q\subseteq \St(p)$. In particular, we conclude $x^{-1}\in \St(p)\cdot p^{-1}\subseteq p\cdot p^{-1}\cdot p^{-1}$, so $x\in p\cdot p\cdot p^{-1}$.

As $X$ is symmetric, $p\cdot \Stab(p)\subseteq X^3$, so $p\cdot \Stab(p)$ is $N${\hyp}medium. By \cref{l:translation of medium}, we conclude that $\Stab(p)$ is $N${\hyp}medium. Thus, we get property (\ref*{itm:stabilizer theorem 2:b}). 

Finally, it remains to prove property (\ref*{itm:stabilizer theorem 2:c}). In other words, assuming that $\mu$ is also invariant under right translations, we want to show that $p\cdot p^{-1}\cdot p=\Stab(p)\cdot p=\Stab(p)\cdot y$ for any $y\in p$. Take $a\in p$ arbitrary and, using the local atomic property, find $b\in p$ such that $\tp(b/N,a)$ is wide. As $p$ is $N${\hyp}medium, by \cref{l:s1 and forking}, $\tp(b/N,a)$ does not fork over $N$. Consider $q_2=\tp(ba^{-1}/N)\subseteq p\prescript{}{\mathrm{nf}(N)}{\cdot} p^{-1}$. Then, by invariance under right translations, $q_2$ is a wide $\bigwedge_N${\hyp}definable set. As $q_2\subseteq X^2$ and $p^{-1}\cdot q_2\subseteq X^3$, we conclude that $q_2$ and $p^{-1}\cdot q_2$ are $N${\hyp}medium. As $p\cdot_{\mathrm{nf}(N)} p^{-1}\subseteq \St(p)$ by \cref{l:lemma stabilizer}, we get using \cref{e:mos 2.8} and \cref{l:fundamental lemma for the stabilizer theorem}(\ref{itm:fundamental lemma for the stabilizer theorem 3}) that $p\prescript{}{\mathrm{nf}(N)}{\cdot} p^{-1}\subseteq \St(p)$ too, so $q_2\subseteq \St(p)$. By \cref{t:stabilizer theorem 1}, we conclude that $\Stab(p)=\Stab(q_2)$ and $\Stab(q_2)\cdot_{\mathrm{nf}(N)}q_2\subseteq \St(p)$. Take $x\in \Stab(p)\cdot y$ arbitrary, so $x=cy$ with $c\in\Stab(p)=\Stab(q_2)$ and $y\in p$. Find $z_0$ such that $\tp(z_0,y/N)=\tp(a,b/N)$, so $\tp(y/N,z_0)$ is wide. Using the local atomic property, we may find $y_1\in \tp(y/N,z_0)$ such that $\tp(y_1/N,z_0,c)$ is wide. Take $z$ such that $\tp(z,y/N,c)=\tp(z_0,y_1/N,c)$. Thus, $yz^{-1}\in q_2$ and $\tp(y/N,z,c)$ is wide. In particular, as $X^2$ is $N${\hyp}medium, by \cref{l:s1 and forking}, $\tp(yz^{-1}/N,c)$ does not fork over $N$. Thus, $x\cdot z^{-1}=c\cdot yz^{-1}\in \Stab(q_2)\cdot_{\mathrm{nf}(N)}q_2\subseteq \St(p)$, concluding $x\in p\cdot p^{-1}\cdot p$. As $x$ is arbitrary, $\Stab(p)\cdot y\subseteq p\cdot p^{-1}\cdot p$. As $p\cdot p^{-1}\cdot p\subseteq \Stab(p)\cdot y$, we get $p\cdot p^{-1}\cdot p=\Stab(p)\cdot p=\Stab(p)\cdot y$ for any $y\in p$. 
\end{proof}
\end{theo}
\begin{rmk} \label{o:remarks stabilizer}
\begin{enumerate}[label={\rm{(\arabic*)}}, ref={\rm{\arabic*}}, wide]
\item[\hspace{-1.4em}\setcounter{enumi}{1}(\theenumi)] \label{itm:remarks stabilizer:no right} The main improvement of \cref{t:stabilizer theorem 2} with respect to the original formulations studied in \cite{hrushovski2011stable,montenegro2018stabilizers} is that we do not require invariance of $\mu$ under right translations. In both papers, while they preferably considered ideals invariant under two{\hyp}sided translations, the authors already anticipated that it should be possible to weaken the invariance under right translations hypothesis, perhaps losing a few properties of $\Stab(p)$. We have shown that, in fact, it can be completely eliminated without significant consequences. Indeed, we only use invariance under right translations for (\ref{itm:stabilizer theorem 2:c}), but this is essentially replaced by (\ref{itm:stabilizer theorem 2:a}). 
\item \label{itm:remarks stabilizer:X^3} Although in \cite[Theorem 2.12]{montenegro2018stabilizers} it was assumed that $X^4$ is $N${\hyp}medium, we only need $N${\hyp}mediumness of $X^3$.
\item \label{itm:remarks stabilizer:locally atomic} Note that the locally atomic property is mostly used only for subsets of $p$. The only time when we apply the locally atomic property for other subsets is in the proof of \cref{t:stabilizer theorem 1}(\ref{itm:stabilizer theorem 1:b}). Thus, when $\mu$ only satisfies the locally atomic property for subsets of $p$, we only miss (\ref{itm:stabilizer theorem 2:b}).  

For property (\ref{itm:stabilizer theorem 2:b}), we have used the locally atomic property for subsets of $\Stab(p)=(p\cdot p^{-1})^2$. However, in fact, it suffices to assume it only for subsets of $p\cdot p\cdot p^{-1}$. Indeed, by \cref{t:stabilizer theorem 2}(\ref{itm:stabilizer theorem 2:a}) and invariance under left translations, if $V\subseteq \Stab(p)$ is $\bigwedge_B${\hyp}definable and wide, then $y\cdot V\subseteq p\cdot p\cdot p^{-1}$ is $\bigwedge_{B,y}${\hyp}definable and wide for any $y\in p$. By the locally atomic property on $p\cdot p\cdot p^{-1}$, there is $b\in V$ such that $\tp(yb/B,y)$ is wide, concluding that $\tp(b/B)$ is wide by invariance under left translations.
\item \label{itm:remarks stabilizer:local} It suffices to have $\mu$ defined only on $X^3$, as, in that case, we may extend the ideal by taking the one generated by the finite unions of left translates of elements of $\mu$. As $\mu$ is invariant under left{\hyp}translations within $X^3$ (i.e. $V\in\mu$ if and only if $gV\in\mu$ for any $V\subseteq X^3$ and $g\in G$ such that $gV\subseteq X^3$), this extension coincides with $\mu$ inside of $X^3$, so we can apply the Stabilizer \cref{t:stabilizer theorem 2} by \cref{o:remarks stabilizer}(\ref{itm:remarks stabilizer:locally atomic}). 
\end{enumerate}
\end{rmk} 
\begin{prop}{\textnormal{{\cite[Corollary 3.11]{hrushovski2011stable}\&\cite[Proposition 2.13]{montenegro2018stabilizers}}}} Let $\mathfrak{N}\prec\mathfrak{M}$ with $|N|<\lambda$ and $G=\underrightarrow{\lim}\, X^n$ a piecewise $N${\hyp}hyperdefinable group generated by a symmetric $\bigwedge_N${\hyp}definable subset $X$. Let $\mu$ be an $N${\hyp}invariant locally atomic ideal of $\bigwedge_{<\lambda}${\hyp}definable subsets of $G$ invariant under left translations such that $X$ is wide and $X^3$ is $N${\hyp}medium. Then, $\mu$ is $N${\hyp}medium on $G$.  
\begin{proof} By the locally atomic property, there is a wide $N${\hyp}type $p\subseteq X$. By the Stabilizer \cref{t:stabilizer theorem 2}, $S\coloneqq \Stab(p)=(p\cdot p^{-1})^2$ is an $N${\hyp}medium normal subgroup of small index. Let $Y(\overline{b})\subseteq G$ be $\bigwedge_{N,\bar{b}}${\hyp}definable with $|\overline{b}|<\lambda$ small. Let $(\overline{b}_i)_{i\in\omega}$ be an $N${\hyp}indiscernible sequence with $\tp(\overline{b}/N)=\tp(\overline{b}_i/N)$. Suppose that $Y(\overline{b})\notin \mu$. By the locally atomic property of $\mu$, there is a wide $N,\overline{b}${\hyp}type $q(\overline{b})\subseteq Y(\overline{b})$. By \cref{l:index lemma 1}(\ref{itm:index lemma 1:left}), there is $a$ such that $q(\overline{b})\subseteq a\cdot S$. Let $(\overline{b}'_i)_{i\in \omega}$  be an $N,a^*${\hyp}indiscernible sequence realising $\tp((\overline{b}_i)_{i\in\omega}/N)$. By invariance under left translations, $a^{-1}q(\overline{b})\notin \mu$. As $S$ is $N${\hyp}medium, $a^{-1}q(\overline{b}'_0)\cap a^{-1}q(\overline{b}'_1)\notin \mu$. Thus, by $N${\hyp}invariance and invariance under left translations, we get that $q(\overline{b}_0)\cap q(\overline{b}_1)\notin \mu$, concluding $Y(\overline{b}_0)\cap Y(\overline{b}_1)\notin \mu$. 
\end{proof}
\end{prop}
\subsection{Near{\hyp}subgroups}
Let $G$ be a piecewise $A${\hyp}hyperdefinable group and $\mu$ an ideal of $\bigwedge_{<\lambda}${\hyp}definable subsets invariant under left translations with $\kappa\geq \lambda>|\lang|+|A|$. A \emph{$\mu${\hyp}near{\hyp}subgroup} of $G$ over $A$ is a wide $\bigwedge_A${\hyp}definable symmetric set $X$ generating $G$ such that $\mu_{\mid X^3}$ is locally atomic and $A${\hyp}medium. We say that $X$ is a near{\hyp}subgroup if it is a $\mu${\hyp}near{\hyp}subgroup for some ideal.
\begin{theo} \label{t:stabilizer nearsubgroups}
Let $G$ be a piecewise $A${\hyp}hyperdefinable group, $\mu$ an ideal of $\bigwedge_{<\lambda}${\hyp}definable subsets and $X$ a $\mu${\hyp}near{\hyp}subgroup over $A$. Then, there is a wide $\bigwedge_{|A|+|\lang|}${\hyp}definable normal subgroup of small index $S$ contained in $X^4$ and contained in every $\bigwedge_A${\hyp}definable subgroup of small index. Furthermore, $S=(p\cdot p^{-1})^2$ and $ppp^{-1}=pS=yS$ for any $y\in p$, where $p\subseteq X$ is a wide type over an elementary substructure of size $|A|+|\lang|$.
\begin{proof} By L\"{o}wenheim{\hyp}Skolem Theorem \cite[Theorem 2.3.1]{tent2012course}, we find an elementary substructure $\mathfrak{N}\preceq \mathfrak{M}$ containing $A^*$ with $|N|=|A^*|+|\lang|<\lambda$. As $A^*\subseteq N$, $\mu_{\mid X^3}$ is $N${\hyp}medium. By the locally atomic property, there is a wide $N${\hyp}type  $p\subseteq X$. Applying the Stabilizer \cref{t:stabilizer theorem 2} (and \cref{o:remarks stabilizer}(\ref{itm:remarks stabilizer:locally atomic})), we conclude that there is a wide  $\bigwedge_N${\hyp}definable normal subgroup $S\trianglelefteq G$ of small index which does not have proper $\bigwedge_N${\hyp}definable subgroups of small index. Furthermore, $S=\Stab(p)=(p\cdot p^{-1})^2\subseteq X^4$ and $ppp^{-1}=pS=yS$ for any $y\in p$. As $A^*\subseteq N$, we conclude. 
\end{proof}
\end{theo}
\begin{theo}\label{t:general lie model}
Let $G=\underrightarrow{\lim}\, X^n$ be a piecewise $A${\hyp}hyperdefinable group generated by a near{\hyp}subgroup $X$ over $A$. Assume that $X^n$ is an approximate subgroup for some $n\in\N$. Then, $G$ has a connected Lie model $\pi_{H/K}:\ H\rightarrow L=\rfrac{H}{K}$ with $K\subseteq X^{2n+4}$ such that
\begin{enumerate}[label={\rm{(\alph*)}}, wide]
\item $H\cap X^{2n}$ is $\bigwedge_{|A|+|\lang|}${\hyp}definable and commensurable to $X^n$,
\item $\pi_{H/K}[H\cap X^{2n}]$ is a compact neighbourhood of the identity in $L$,
\item $H$ is generated by $H\cap X^{2n+4}$, and
\item $\pi_{H/K}$ is continuous and proper from the logic topology using $|A|+|\lang|$ many parameters. 
\end{enumerate}
\begin{proof} As $X$ is a near{\hyp}subgroup, by \cref{t:stabilizer nearsubgroups}, $G^{00}_N\subseteq X^4$ with $|N|\leq |A|+|\lang|$. As $X^n$ is an approximate subgroup, $\rfrac{X^n}{G^{00}_N}$ is an approximate subgroup. Thus, $\rfrac{X^{2n}}{G^{00}_N}$ is a neighbourhood of the identity by the Generic Set Lemma (\cref{t:generic set lemma}). Then, by \cref{t:logic existence of lie core}, we can find a connected Lie core $\pi_{H/K}:\ H\rightarrow L=\rfrac{H}{K}$ with $K\subseteq X^{2n}G^{00}_N\subseteq X^{2n+4}$. \begin{enumerate}[label={\rm{(\alph*)}}, wide]
\item Since $\rfrac{X^{n}}{G^{00}_A}$ is compact and $\rfrac{H}{G^{00}_A}$ is an open subgroup in the global logic topology, we conclude that $\left[\rfrac{X^n}{G^{00}_A}:\rfrac{H}{G^{00}_A}\right]$ is finite. As $G^{00}_A\leq H$, we get that $[X^n:H]$ is finite. Thus, by \cite[Lemma 2.2, Lemma 2.3]{machado2023closed}, $H\cap X^{2n}$ is an approximate subgroup commensurable to $X^{n}$.
\item As $\rfrac{X^{2n}}{G^{00}_N}$ is a compact neighbourhood of the identity and $\rfrac{H}{G^{00}_N}$ is an open subgroup in the global logic topology, we conclude that $\rfrac{H\cap X^{2n}}{G^{00}_N}$ is a compact neighbourhood of the identity. Recall that $\widetilde{\pi}_{H/K}:\ \rfrac{H}{G^{00}_N}\rightarrow L$ given by $\pi_{H/K}=\widetilde{\pi}_{H/K}\circ\pi_{G/G^{00}_N}$ is a continuous, closed, open and proper onto group homomorphism. Thus, we conclude that $\pi_{H/K}[H\cap X^{2n}]$ is a compact neighbourhood of the identity. 
\item Since $L$ is connected and $\pi_{H/K}[X^{2n}\cap H]$ is a neighbourhood of the identity, $\pi_{H/K}[H\cap X^{2n}]$ generates $L$. Therefore, we have that $\pi^{-1}_{H/K}[\pi_{H/K}[H\cap X^{2n}]]=(H\cap X^{2n})\cdot K$ generates $H$. Now, $K\subseteq H\cap X^{2n+4}$, so $(H\cap X^{2n})\cdot K\subseteq  (H\cap X^{2n+4})^2$, concluding that $H\cap X^{2n+4}$ generates $H$.
\item By \cref{p:reduce number of parameters}, the global logic topology of $\rfrac{G}{G^{00}_N}$ is given using $|A|+|\lang|$ many parameters, so $H\cap X^{2n}$ is $\bigwedge_{|A|+|\lang|}${\hyp}definable and $\pi_{H/K}$ is continuous and closed using $|A|+|\lang|$ many parameters.\qedhere
\end{enumerate}
\end{proof}
\end{theo}
Alternatively, using Gleason{\hyp}Yamabe{\hyp}Kreitlon \cref{t:gleason yamabe carolino} rather than Gleason{\hyp}Yamabe \cref{t:gleason yamabe}, we get the following variation which provides some extra control over some of the parameters.
\begin{theo}\label{t:general lie model carolino}
There are functions $c:\ \N\rightarrow \N$ and $d:\ \N\rightarrow\N$ such that the following holds for any $\kappa${\hyp}saturated $\lang${\hyp}structure $\mathfrak{M}$ and any set of parameters $A$ with $\kappa>|\lang|+|A|$:\smallskip

Let $G=\underrightarrow{\lim}\, X^n$ be a piecewise $A${\hyp}hyperdefinable group generated by a near{\hyp}subgroup $X$ over $A$. Assume that $X^n$ is a $k${\hyp}approximate subgroup for some $n$ and $k$. Then, $G$ has a Lie model $\pi_{H/K}:\ H\rightarrow L=\rfrac{H}{K}$ with $K\subseteq X^{12n+4}$ and $\dim(L)\leq d(k)$ such that
\begin{enumerate}[label={\rm{(\alph*)}}, wide]
\item $H\cap X^{2n}$ is $\bigwedge_{|A|+|\lang|}${\hyp}definable and $c(k)${\hyp}commensurable to $X^n$,
\item $\pi_{H/K}[H\cap X^{2n}]$ is a compact neighbourhood of the identity in $L$,
\item $H$ is generated by $H\cap X^{12n+4}$, and
\item $\pi_{H/K}$ is continuous and proper from the logic topology using $|A|+|\lang|$ many parameters. 
\end{enumerate}
\begin{proof} As $X$ is a near{\hyp}subgroup, by \cref{t:stabilizer nearsubgroups}, $G^{00}_N\subseteq X^4$ with $|N|\leq |A|+|\lang|$. As $X^n$ is a $k${\hyp}approximate subgroup, $\rfrac{X^n}{G^{00}_N}$ is an approximate subgroup. Thus, $\rfrac{X^{2n}}{G^{00}_N}$ is a neighbourhood of the identity by the Generic Set Lemma (\cref{t:generic set lemma}), and, in particular, $\rfrac{G}{G^{00}_N}$ is a locally compact topological group with the global logic topology. Thus, $\rfrac{X^n}{G^{00}_N}$ is contained in the interior, $\rfrac{U}{G^{00}_N}$, of $\rfrac{X^{3n}}{G^{00}_N}$ in the global logic topology. Now, we obviously have that $\rfrac{U^2}{G^{00}_N}\subseteq \rfrac{X^{6n}}{G^{00}_N}$, which is covered by $k^5$ left translates of $\rfrac{X^n}{G^{00}_N}\subseteq \rfrac{U}{G^{00}_N}$. Hence, $\rfrac{U}{G^{00}_N}$ is an open precompact $k^5${\hyp}approximate subgroup.

Let $c_0$ and $d_0$ be the functions provided by Gleason{\hyp}Yamabe{\hyp}Kreitlon \cref{t:gleason yamabe carolino}. Applying \cref{t:gleason yamabe carolino}, we get a Lie model $\pi_{H/K}:\ H\rightarrow L=\rfrac{H}{K}$ with $G^{00}_N\leq K\subseteq U^4\subseteq X^{12n}G^{00}_N\subseteq X^{12n+4}$ and $\dim(L)\leq d(k)\coloneqq d_0(k^5)$ such that $H\cap U^4$ generates $H$ and is $c_0(k^5)${\hyp}commensurable to $U$. Thus, $H\cap X^{12n+4}$ generates $H$ and $c_0(k^5)$ left translates of $H\cap X^{12n+4}$ cover $X^n$. As $X^n$ is a $k${\hyp}approximate subgroup, $k^{15}$ left translates of $X^n$ cover $X^{16n}\supseteq H\cap X^{12n+4}$. By \cite[Lemma 2.2]{machado2023closed}, we get that $k^{15}$ left translates of $H\cap X^{2n}$ cover $H\cap X^{12n+4}$, so $H\cap X^{2n}$ and $X^n$ are $c(k)${\hyp}commensurable, where $c(k)\coloneqq c_0(k^5)k^{15}$. \smallskip

Now, $\rfrac{H\cap X^{2n}}{G^{00}_N}$ is a compact neighbourhood of the identity in the global logic topology. Recall that $\widetilde{\pi}_{H/K}:\ \rfrac{H}{G^{00}_N}\rightarrow L$ given by $\pi_{H/K}=\widetilde{\pi}_{H/K}\circ\pi_{G/G^{00}_N}$ is a continuous, closed, open and proper onto group homomorphism. Thus, $\pi_{H/K}[H\cap X^{2n}]$ is a compact neighbourhood of the identity of $L$. 

By \cref{p:reduce number of parameters}, the global logic topology of $\rfrac{G}{G^{00}_N}$ is given using $|A|+|\lang|$ many parameters, so $H\cap X^{2n}$ is $\bigwedge_{|A|+|\lang|}${\hyp}definable and $\pi_{H/K}$ is continuous and proper using $|A|+|\lang|$ many parameters.
\end{proof}
\end{theo}
We conclude applying \cref{t:general lie model} to the case of rough approximate subgroups. Recall that a \emph{$T${\hyp}rough $k${\hyp}approximate subgroup} of a group $G$ is a subset $X\subseteq G$ such that $X^2\subseteq \Delta XT$ with $|\Delta|\leq k\in\N_{>0}$ and $1\in T\subseteq G$. \medskip

Let $G$ be a group and $X\subseteq G$ a symmetric subset. For some fixed $n$, assume $X^n$ is a $T_i${\hyp}rough $k${\hyp}approximate subgroup for a sequence $(T_i)_{i\in\N}$ of thickenings satisfying that
\begin{enumerate}[label={\rm{(\roman*)}}, wide]
\item it decreases in doubling scales, i.e. $T_{i+1}T^{-1}_{i+1}\subseteq T_i$ for each $i\in\N$, and
\item it is asymptotically normalised by $X$, i.e. $x^{-1}T_{i+1}x\subseteq T_i$ for each $x\in X$ and $i\in \N$.
\end{enumerate}

Assuming saturation and $\bigwedge${\hyp}definability, we have that $T=\bigcap T_i$ is an $\bigwedge${\hyp}definable subgroup of $G$ normalised by $X$ and $X^n$ is a $T${\hyp}rough $k${\hyp}approximate subgroup. 
\begin{theo}[Rough Lie Model Theorem] \label{t:rough lie model}
Let $G$ be an $A${\hyp}definable group, $T\leq G$ an $\bigwedge_{A}${\hyp}definable subgroup and $X\subseteq G$ a symmetric $\bigwedge_{A}${\hyp}definable subset. Write $\widetilde{G}$ for the subgroup generated by $X$. Assume that $X$ normalises $T$, $\rfrac{X}{\kern 0.1em T}$ is a near{\hyp}subgroup of $\rfrac{\widetilde{G}}{\kern 0.1em T}$ and $X^n$ is a $T${\hyp}rough approximate subgroup for some $n$. Then, $\widetilde{G}T\leq G$ has a connected Lie model $\pi_{H/K}:\ H\rightarrow L=\rfrac{H}{K}$ with $T\subseteq K\subseteq X^{2n+4}T$ such that
\begin{enumerate}[label={\rm{(\alph*)}}, wide]
\item $H\cap X^{2n}T$ is $\bigwedge_{|A|+|\lang|}${\hyp}definable and commensurable to $X^nT$,
\item $\pi_{H/K}[H\cap X^{2n}]$ is a compact neighbourhood of the identity in $L$,
\item $H$ is generated by $H\cap X^{2n+4}T$, and
\item $\pi_{H/K}$ is continuous and proper from the logic topology using $|A|+|\lang|$ many parameters.
\end{enumerate}
\end{theo} 
Alternatively, using \cref{t:general lie model carolino}:
\begin{theo}[Rough Lie model Theorem, version 2] \label{t:rough lie model carolino}
There are functions $c:\ \N\rightarrow \N$ and $d:\ \N\rightarrow\N$ such that the following holds for any $\kappa${\hyp}saturated $\lang${\hyp}structure $\mathfrak{M}$ and any set of parameters $A$ with $\kappa>|\lang|+|A|$:\smallskip

Let $G$ be an $A${\hyp}definable group, $T\leq G$ an $\bigwedge_{A}${\hyp}definable subgroup and $X\subseteq G$ a symmetric $\bigwedge_{A}${\hyp}definable subset. Write $\widetilde{G}$ for the subgroup generated by $X$. Assume that $X$ normalises $T$, $\rfrac{X}{\kern 0.1em T}$ is a near{\hyp}subgroup of $\rfrac{\widetilde{G}}{\kern 0.1em T}$ and $X^n$ is a $T${\hyp}rough $k${\hyp}approximate subgroup for some $n$ and $k$. Then, $\widetilde{G}T\leq G$ has a Lie model $\pi_{H/K}:\ H\rightarrow L=\rfrac{H}{K}$ with $T\subseteq K\subseteq X^{12n+4}T$ and $\dim(L)\leq d(k)$ such that
\begin{enumerate}[label={\rm{(\alph*)}}, wide]
\item $H\cap X^{2n}T$ is $\bigwedge_{|A|+|\lang|}${\hyp}definable and $c(k)${\hyp}commensurable to $X^nT$,
\item $\pi_{H/K}[H\cap X^{2n}]$ is a compact neighbourhood of the identity in $L$,
\item $H$ is generated by $H\cap X^{12n+4}T$, and
\item $\pi_{H/K}$ is continuous and proper from a logic topology using $|A|+|\lang|$ many parameters.
\end{enumerate}
\end{theo}

\bibliographystyle{alphaurl}
\bibliography{On_piecewise_hyperdefinable_groups}
\end{document}